\documentclass[a4j,11pt,dvipdfmx]{article}

\usepackage{graphicx}
\usepackage{amsmath} %
\usepackage{amssymb}
\usepackage{amsthm}
\usepackage{mathrsfs}

\usepackage{color}
\usepackage{amsthm} 
\usepackage{subfigure}
\usepackage{caption}
\usepackage{bm}
\usepackage{framed}
\usepackage{here}
\usepackage{url}
\usepackage{algorithm}
\usepackage{algorithmic}     
\usepackage{multirow}

\newtheorem{Thm}{Theorem}
\newtheorem{Prop}{Proposition}
\newtheorem{lemm}{Lemma}
\newtheorem{Cor}{Corollary}
\theoremstyle{definition}
\newtheorem{Rem}{Remark}
\newtheorem*{Condition}{Condition}

\newcommand{\indepe}{\mathop{\perp\!\!\!\perp}}
\newcommand{\norm}[1]{\left\|#1\right\|}



\allowdisplaybreaks[4]

\title{Two-sample test based on maximum variance discrepancy}
\author{Natsumi Makigusa\\
	\\
	Graduate School of Science and Engineering, Chiba University}
\date{}
\begin{document}
\maketitle
%
\begin{abstract}
In this article, we introduce a novel discrepancy called the maximum variance discrepancy for the purpose of measuring the difference between two distributions in Hilbert spaces that cannot be found via the maximum mean discrepancy. 
We also propose a two-sample goodness of fit test based on this discrepancy. 
We obtain the asymptotic null distribution of this two-sample test, which provides an efficient approximation method for the null distribution of the test.
\end{abstract}
\section{Introduction}
For probability distributions $P$ and $Q$, the test for the null hypothesis $H_0: P=Q$ against an alternative hypothesis $H_1:P \neq Q$ based on data $X_1, \dots, X_n \overset{i.i.d.}{\sim} P$ and $Y_1, \dots ,Y_m \overset{i.i.d.}{\sim} Q$ is known as a two-sample  test.
Such tests have applications in various areas. 
There is a huge body of literature on two-sample tests in Euclidean space, so we will not attempt a complete bibliography. 
In \cite{Gretton2007}, a two-sample test based on Maximum Mean Discrepancy (MMD) is proposed, where the MMD is defined by (\ref{Eq_MMD}) in Section \ref{MVD}.
The MMD for a reproducing kernel Hilbert space $H(k)$ associated with a positive definite kernel $k$ is defined as (\ref{Eq_MMD_inRKHS}) in Section \ref{MVD}.

In this paper, we propose a novel discrepancy between two distributions defined as
\[
T=\norm{V_{X \sim P}[k(\cdot,X)]-V_{Y \sim Q}[k(\cdot,Y)]}_{H(k)^{\otimes 2}},
\]
and we call this the Maximum Variance Discrepancy (MVD), where $V_{X \sim P}[k(\cdot,X)]$ is a covariance operator in $H(k)$.
The MVD is composed by replacing the kernel mean embedding in (\ref{Eq_MMD_inRKHS}) with a covariance operator; hence, it is natural to consider a two-sample test based on the MVD.
A related work can be found in \cite{Boente2018}, where a test for the equality of covariance operators in Hilbert spaces was proposed.

Our aim in this research is to clarify the properties
of the MVD test from two perspectives: an asymptotic investigation as $n,m \to \infty$, and its practical implementation. 
We first obtain the asymptotic distribution of a consistent estimator $\widehat{T}_{n,m}^2$ of $T^2$, under $H_0$.
We also derive the asymptotic distribution of $\widehat{T}_{n,m}^2$ under the alternative hypothesis $H_1$. 
Furthermore, we consider a sequence of local alternative distributions $Q_{nm}=(1-1/\sqrt{n+m})P+(1/\sqrt{n+m})Q$ for $P \neq Q$ and address the asymptotic distribution of $\widehat{T}^2_{n,m}$ under this sequence.
For practical purposes, a method to approximate the distribution of the test by $\widehat{T}^2_{n,m}$ under $H_0$ is developed.
The method is based on the eigenvalues of the centered Gram matrices associated with the dataset.
Those eigenvalues will be shown to be estimators of the weights appearing in the asymptotic null distribution of the test.
Hence, the method based on the eigenvalues is expected to provide a fine approximation of the distribution of the test.
However, this approximation does not actually work well.
Therefore, we further modify the method based on the eigenvalues, and the obtained method provides a better approximation.

The rest of this paper is structured as follows. 
Section \ref{MVD} introduces the framework of the two-sample test and defines the test statistics based on the MVD. 
In addition, the representation of test statistics based on the centered Gram matrices is described.
Section \ref{asymptotic null distribution} develops the asymptotics for the test by the MVD under $H_0$.
The test by the MVD under $H_1$ is addressed in Section \ref{asymptotic nonnull distribution}. 
Furthermore, the behavior of the test by the MVD under the local alternative hypothesis is clarified in Section \ref{asymptotic distribution under contiguous alternatives}.
Section \ref{Gram matrix spectrum} describes the estimation of the weights that appear in the asymptotic null distribution obtained in Section \ref{asymptotic null distribution}. 
Section \ref{implementation} examines the implementation of the MVD test with a Gaussian kernel in the Hilbert space $\mathcal{H}=\mathbb{R}^d$.
Section \ref{approximation of the null distribution} introduces the modification of the approximate distribution given in Section \ref{Gram matrix spectrum}. Section \ref{simulations} reports the results of simulations for the type I error and the power of the MVD and MMD tests.
Section \ref{Applications to real data sets} presents the results of applications to real data sets, including high-dimension low-sample-size data.
Conclusions are given in Section \ref{Conclusion}. 
All proofs of theoretical results are provided in Section \ref{section_Proof}.

\section{Maximum Variance Discrepancy}\label{MVD}
Let $\mathcal{H}$ be a separable Hilbert space and $(\mathcal{H},\mathcal{A})$ be a measurable space.
Let $\left<\cdot,\cdot\right>_{\mathcal{H}}$ be the inner product of $\mathcal{H}$ and $\|\cdot\|_{\mathcal{H}}=\sqrt{\left<\cdot,\cdot\right>_{\mathcal{H}}}$ be the associated norm.
Let $X_1,\dots,X_n \in \mathcal{H}$ and $Y_1,\dots,Y_m \in \mathcal{H}$ denote a sample of independent and identically distributed (i.i.d.) random variables drawn from unknown distributions $P$ and $Q$, respectively.
Our goal is to test whether the unknown distributions $P$ and $Q$ are equal.

Let us define the null hypothesis $H_0:P=Q$ and the alternative hypothesis $H_1:P \neq Q$.
Following \cite{Gretton2007}, the gap between two distributions $P$ and $Q$ on $\mathcal{H}$ is measured by:
\begin{equation}\label{Eq_MMD}
\text{MMD}(P,Q)= \sup_{f \in \mathcal{F}} |\mathbb{E}_{X \sim P}[f(X)]-\mathbb{E}_{Y \sim Q}[f(Y)]|,	
\end{equation}
where $\mathcal{F}$ is a class of real-valued functions on $\mathcal{H}$.
Regardless of $\mathcal{F}$, $\text{MMD}(P,Q)$ always defines a pseudo-metric on the space of probability distributions. 
Let $\mathcal{F}$ be the unit ball of a reproducing kernel Hilbert space $H(k)$ associated with a characteristic kernel $k: \mathcal{H} \times \mathcal{H} \to \mathbb{R}$ (see \cite{Aronszajn} and \cite{Fukumizu} for details) and assume that $\mathbb{E}_{X \sim P}[\sqrt{k(X,X)}] < \infty$ and $\mathbb{E}_{Y \sim Q}[\sqrt{k(Y,Y)}] < \infty$. 
Then, the MMD in $H(k)$ is defined as the distance between $P$ and $Q$ as follows:
\begin{equation}\label{Eq_MMD_inRKHS}
	\text{MMD}(P,Q)=\sup_{\norm{f}_{H(k)}=1}|\left<f,\mu_k(P)-\mu_k(Q)\right>_{H(k)}|=\norm{\mu_k(P)-\mu_k(Q)}_{H(k)},	
\end{equation}
where $\mu_k(P)=\mathbb{E}_{X \sim P}[k(\cdot,X)]$ and $\mu_k(Q)=\mathbb{E}_{Y \sim Q}[k(\cdot,Y)]$ are called kernel mean embeddings of $P$ and $Q$, respectively, in $H(k)$ (see \cite{Gretton2007}).
The MMD focuses on the difference between distributions $P$ and $Q$ depending on the difference between the means of $k(\cdot,X)$ and $k(\cdot,Y)$ in $H(k)$.
The motivation for this research is to focus on the difference between the distributions $P$ and $Q$ due to the difference between those variances in $H(k)$, based on a similar idea as the MMD.
Assume $\mathbb{E}_{X \sim P}[k(X,X)] < \infty$ and $\mathbb{E}_{Y \sim Q}[k(Y,Y)] < \infty$, then the variance $V_{X \sim P}[k(\cdot,X)]:H(k) \to H(k)$ is defined by
\[
V_{X \sim P}[k(\cdot,X)]=\mathbb{E}_{X \sim P}[(k(\cdot,X)-\mu_k(P))^{\otimes 2}] \in H(k)^{\otimes 2}.
\]
Here, for any $h,h' \in H(k)$, the tensor product $h \otimes h'$ is defined as the operator $H(k) \to H(k), x \mapsto \left<h',x\right>_{H(k)}h$,  $h^{\otimes 2}$ is defined as $h^{\otimes 2}=h \otimes h$, and $H(k)^{\otimes 2} =H(k) \otimes H(k)$ (see Section II.4 in \cite{reed1981functional} for details).
Let $V_{X \sim P}[k(\cdot,X)]=\Sigma_k(P)$ and $V_{Y \sim Q}[k(\cdot,Y)]=\Sigma_k(Q)$. Then, we define the MVD in $H(k)$ as
\[
\text{MVD}(P,Q)=\sup_{\norm{A}_{H(k)^{\otimes 2}}=1}|\left<A,\Sigma_k(P)-\Sigma_k(Q)\right>_{H(k)^{\otimes 2}}|=\norm{\Sigma_k(P)-\Sigma_k(Q)}_{H(k)^{\otimes 2}},
\]
which can be seen as a discrepancy between distributions $P$ and $Q$. 
The $T^2=\text{MVD}(P,Q)^2$ can be estimated by
\begin{equation}\label{Eq_widehat_Tnm}
\widehat{T}^2_{n,m}=\norm{\Sigma_k(\widehat{P})-\Sigma_k(\widehat{Q})}^2_{H(k)^{\otimes 2}}, 
\end{equation}
where
\[
\Sigma_k(\widehat{P})=\frac{1}{n}\sum_{i=1}^{n}(k(\cdot,X_i)-\mu_k(\widehat{P}))^{\otimes 2},~~\mu_k(\widehat{P})=\frac{1}{n}\sum_{i=1}^{n}k(\cdot,X_i)
\]
and
\[
\Sigma_k(\widehat{Q})=\frac{1}{m}\sum_{j=1}^{m}(k(\cdot,Y_j)-\mu_k(\widehat{Q}))^{\otimes 2},~~\mu_k(\widehat{Q})=\frac{1}{m}\sum_{j=1}^{m}k(\cdot,Y_j).
\]
Let the Gram matrices be $K_{X}=(k(X_i,X_s))_{1 \leq i,s \leq n}$, $K_{Y}=(k(Y_j,Y_t))_{1 \leq j,t \leq m}$, and $K_{XY}=(k(X_i,Y_j))_{\substack{1 \leq i \leq n\\ 1 \leq j \leq m}}$; the centering matrix be $P_n=I_n-(1/n)\underline{1}_{n}\underline{1}_{n}$; and the centered Gram matrices be $\widetilde{K}_{X}=P_n K_{X} P_n$, $\widetilde{K}_{Y} =P_m K_{Y}P_m$, and $\widetilde{K}_{XY} =P_nK_{XY}P_m$, where $I_n$ is the $n \times n$ identity matrix. 
This test statistic can be expanded as:
\[
\widehat{T}^2_{n,m}
=\frac{1}{n^2} \text{tr}(\widetilde{K}_X^2) -\frac{2}{nm} \text{tr}(\widetilde{K}_{XY}\widetilde{K}_{XY}^T) +\frac{1}{m^2} \text{tr}(\widetilde{K}_Y^2).
\]

We investigate the $\text{MMD}(P,Q)^2$ and $\text{MVD}(P,Q)^2$ when $\mathcal{H}=\mathbb{R}^d$, the kernel $k(\cdot,\cdot)$ is the Gaussian kernel:
\begin{equation}\label{Eq_Gaussign_kernel}
	k(\underline{t},\underline{s})=\exp(-\sigma \norm{\underline{t}-\underline{s}}^2_{\mathbb{R}^d}),~~\sigma >0,
\end{equation}
$P=N(\underline{0},I_d)$, and $Q=N(\underline{m},\Sigma)$.
Under this setting, straightforward calculations using the properties of Gaussian density yield the following result for MMD:
\begin{Prop}\label{Prop_MMD_N(0,1)_and_N(m,Sigma)}
	When $\mathcal{H}=\mathbb{R}^d$ and $k(\cdot,\cdot)$ is the Gaussian kernel in (\ref{Eq_Gaussign_kernel}), we have
	\begin{align*}
		&{\rm MMD}(N(\underline{0},I_d),N(\underline{m},\Sigma))^2\\
		&=(1+4\sigma)^{-d/2}+|I_d+4\sigma\Sigma|^{-1/2}
		-2|(1+2\sigma)I_d+2\sigma\Sigma|^{-1/2} \exp\left(
		-\sigma \underline{m}^T ((1+2\sigma)I_d +2\sigma\Sigma)^{-1} \underline{m}
		\right).
	\end{align*}
\end{Prop}
The result for MVD is also obtained as follows using the Gaussian density property as well as the result for MMD:
\begin{Prop}\label{Prop_MVD_N(0,1)_and_N(m,Sigma)}
	When $\mathcal{H}=\mathbb{R}^d$ and $k(\cdot,\cdot)$ is the Gaussian kernel in (\ref{Eq_Gaussign_kernel}), we have
	\begin{align*}
		&{\rm MVD}(N(\underline{0},I_d),N(\underline{m},\Sigma))^2\\
		&=(1+8\sigma)^{-d/2}-2(1+8\sigma+12\sigma^2)^{-d/2}+(1+4\sigma)^{-d}\\
		&~~~~~+|I_d+8\sigma \Sigma|^{-1/2}
		-2|I_d+8\sigma \Sigma+12\sigma^2 \Sigma^2|^{-1/2}
		+|I_d+4\sigma \Sigma|^{-1}\\
		&~~~~~-2|(1+4\sigma)I_d+4\sigma \Sigma|^{-1/2} \exp\left(
		-2\sigma \underline{m}^T\left(
		(1+4\sigma)I_d+4\sigma \Sigma
		\right)^{-1} \underline{m}
		\right) \\
		&~~~~~+2|I_d+2\sigma \Sigma|^{-1/2} |(1+4\sigma)I_d+2\sigma \Sigma|^{-1/2}
		\exp\left(
		-2\sigma \underline{m}^T \left(
		(1+4\sigma)I_d+2\sigma \Sigma
		\right)^{-1} \underline{m}
		\right)\\
		&~~~~~
		+2(1+2\sigma)^{-d/2} |(1+2\sigma)I_d+4\sigma \Sigma|^{-1/2}
		\exp\left(-2\sigma \underline{m}^T \left((1+2\sigma)I_d+4\sigma \Sigma\right)^{-1} \underline{m}\right)\\
		&~~~~~
		-2|(1+2\sigma)I_d+2\sigma \Sigma|^{-1} \exp\left(
		-2\sigma \underline{m}^T \left(
		(1+2\sigma)I_d+2\sigma \Sigma
		\right)^{-1} \underline{m}
		\right). 
	\end{align*}
\end{Prop}
In particular, $\text{MMD}(P,Q)$ and  $\text{MVD}(P,Q)$ for $P=N(\underline{0},I_d)$ and $Q=N(t\underline{1},sI_d)$ are derived by Propositions \ref{Prop_MMD_N(0,1)_and_N(m,Sigma)} and \ref{Prop_MVD_N(0,1)_and_N(m,Sigma)}.
The result is Corollary \ref{Cor_MMD_MVD}.
\begin{Cor}\label{Cor_MMD_MVD}
	When $\mathcal{H}=\mathbb{R}^d$ and $k(\cdot,\cdot)$ is the Gaussian kernel in (\ref{Eq_Gaussign_kernel}), we have
	\begin{align*}
		&{\rm MMD}(N(\underline{0},I_d),N(t\underline{1},sI_d))^2\\
		&=(1+4\sigma)^{-d/2}+(1+4\sigma s )^{-d/2}
		-2 (1+2\sigma+2\sigma s)^{-d/2} \exp\left(
		-\sigma t^2 d (1+2\sigma +2\sigma s)^{-1}  
		\right)
	\end{align*}
	and
	\begin{align*}
		&{\rm MVD}(N(\underline{0},I_d),N(t\underline{1},sI_d))^2\\
		&=(1+8\sigma)^{-d/2}-2(1+8\sigma+12\sigma^2)^{-d/2}+(1+4\sigma)^{-d}\\
		&~~~~~+(1+8\sigma s)^{-d/2}
		-2(1+8\sigma s+12\sigma^2 s^2)^{-d/2}
		+(1+4\sigma s)^{-d}\\
		&~~~~~-2(1+4\sigma+4\sigma s)^{-d/2}
		\exp\left(
		-2\sigma t^2 d (1+4\sigma +4\sigma s)^{-1}
		\right)\\
		&~~~~~
		+2(1+2\sigma s)^{-d/2}(1+4\sigma +2\sigma s)^{-d/2}
		\exp\left(
		-2\sigma t^2 d (1+4\sigma+2\sigma s)^{-1}
		\right) \\
		&~~~~~
		+2(1+2\sigma	)^{-d/2}(1+2\sigma +4\sigma s)^{-d/2}
		\exp\left(
		-2\sigma t^2 d (1+2\sigma+4\sigma s)^{-1}
		\right) \\
		&~~~~~
		-2(1+2\sigma+2\sigma s)^{-d}
		\exp\left(
		-2\sigma t^2 d (1+2\sigma+2\sigma s)^{-1}
		\right). 
	\end{align*}
\end{Cor}
We investigate the behavior of ${\rm MMD}(N(\underline{0},I_d),N(t\underline{1},sI_d))^2$ and ${\rm MVD}(N(\underline{0},I_d),N(t\underline{1},sI_d))^2$ for the difference of $(t,s)$ from (0,1) by using Corollary \ref{Cor_MMD_MVD}.
A sensitive reaction to the difference of $(t,s)$ from (0,1) means that it can sensitively react to differences between distributions from $N(\underline{0},I_d)$.
Using such a discrepancy in the framework of the test is expected to correctly reject $H_0$ under $H_1$.

More generally, kernel $k'(x,y)=\exp(C)k(x,y)$ based on a constant $C$ and a positive definite kernel $k(x,y)$ is also positive definite.
Then, $\text{MMD}_{k'}(P,Q)$ and $\text{MVD}_{k'}(P,Q)$ calculated by the kernel $k'$ hold $\text{MMD}_{k'}(P,Q)^2=\exp(C)\text{MMD}_{k}(P,Q)^2$ and $\text{MVD}_{k'}(P,Q)^2=\exp(2C)\text{MVD}_{k}(P,Q)^2$ for any distributions $P$ and $Q$ using $\text{MMD}_{k}(P,Q)$ and $\text{MVD}_{k}(P,Q)$ calculated by the kernel $k$.

\begin{figure}[H]
	\centering
	\includegraphics[width=1\linewidth]{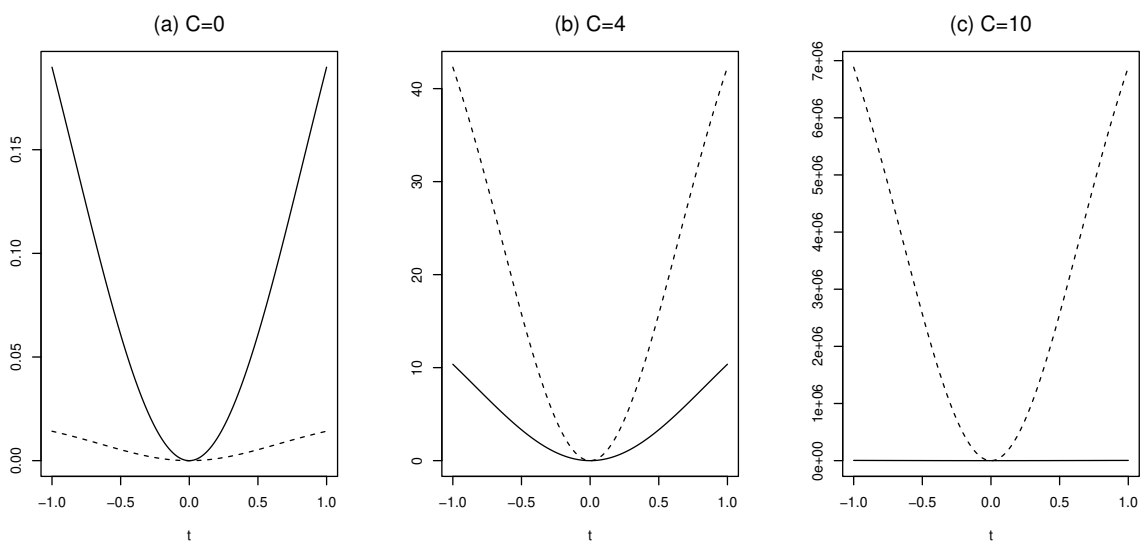}
	\caption{The $\text{MMD}_{k'}(N(\underline{0},I_d),N(t\underline{1},sI_d))^2$ (solid) and $\text{MVD}_{k'}(N(\underline{0},I_d),N(t\underline{1},sI_d))^2$ (dashed) for each $t$: $s=1$, $d=10$, $\sigma=0.1$, and (a) $C=0$, (b) $C=4$, and (c) $C=10$.}
	\label{fig:mmdmvdt}
\end{figure}

\begin{figure}[H]
	\centering
	\includegraphics[width=1\linewidth]{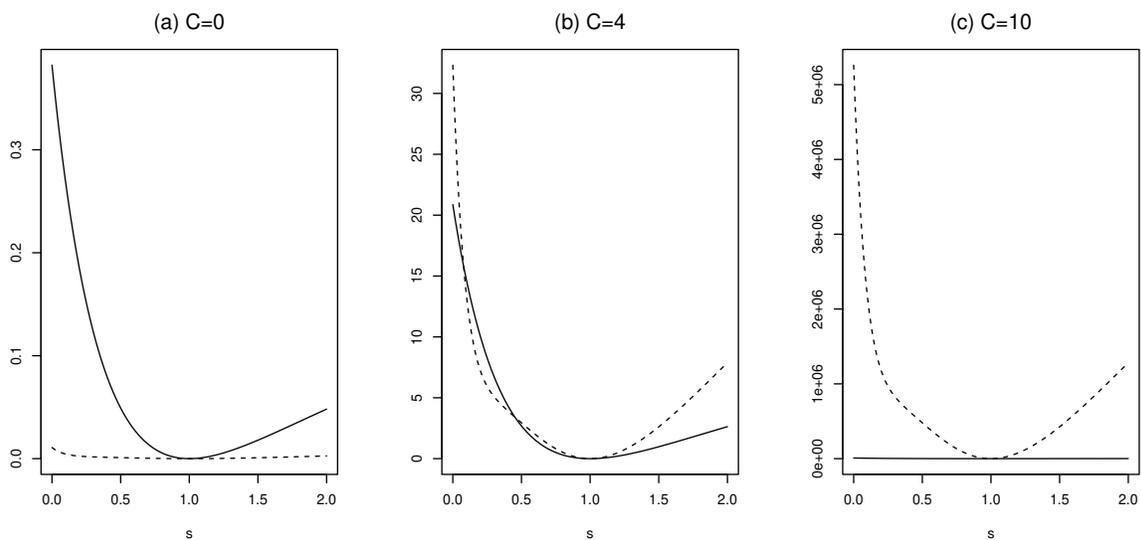}
	\caption{The $\text{MMD}_{k'}(N(\underline{0},I_d),N(t\underline{1},sI_d))^2$ (solid) and $\text{MVD}_{k'}(N(\underline{0},I_d),N(t\underline{1},sI_d))^2$ (dashed) for each $s$: $t=0$, $d=10$, $\sigma=0.1$, and (a) $C=0$, (b) $C=4$, and (c) $C=10$.}
	\label{fig:mmdmvds}
\end{figure}

The graph of $\text{MMD}_{k'}(P,Q)^2$ and $\text{MVD}_{k'}(P,Q)^2$ is displayed for each $t$ when $s=1$ in Figure \ref{fig:mmdmvdt} and for each $s$ when $t=0$ in Figure \ref{fig:mmdmvds}.
The kernel $k$ is a Gaussian kernel in (\ref{Eq_Gaussign_kernel}), and the parameters are $C=0, 4, 10$, $d=10$, and $\sigma=0.1$ in both Figures \ref{fig:mmdmvdt} and \ref{fig:mmdmvds}.
Figure \ref{fig:mmdmvdt} shows the $\text{MMD}_{k'}(P,Q)^2$ and $\text{MVD}_{k'}(P,Q)^2$ for the difference of the mean from the standard normal distribution.
In Figure \ref{fig:mmdmvdt} (a), $\text{MMD}_{k'}(P,Q)^2$ is larger than $\text{MVD}_{k'}(P,Q)^2$, but in Figures \ref{fig:mmdmvdt} (b) and (c), where the value of $C$ is increased, $\text{MVD}_{k'}(P,Q)^2$ is larger than $\text{MMD}_{k'}(P,Q)^2$ for each $t$.
In addition, Figure \ref{fig:mmdmvds} shows the reaction of the MMD and MVD to the difference of the covariance matrix from the standard normal distribution, and $\text{MVD}_{k'}(P,Q)^2$ is larger than $\text{MMD}_{k'}(P,Q)^2$ for each $s$ when $C$ is large.
This means that MVD is more sensitive to differences from the standard normal distribution than MMD for $ k'$ with large $C$.
The fact that there is a kernel $k'$ for which MVD is larger than MMD is a motivation for the two-sample test based on MVD in the next section.


\section{Test statistic for two-sample problem}\label{asymptotics}
We consider a two-sample test based on $T^2_{n,m}$ for $H_0:P=Q$ and $H_1: P \neq Q$, and $\widehat{T}^2_{n,m}$ is defined as a test statistic.
If $\widehat{T}^2_{n,m}$ is large, then the null hypothesis $H_0$ is rejected since $T^2_{n,m}$ is the difference between $P$ and $Q$.
The condition to derive the asymptotic distribution of this test statistic is as follows:
\begin{Condition}
	$\mathbb{E}_{X \sim P}[k(X,X)^2] < \infty$, $\mathbb{E}_{Y \sim Q}[k(Y,Y)^2] < \infty$ and $\lim_{n,m \to \infty} n/(n+m) \to \rho, ~0 < \rho <1$.
\end{Condition}
\subsection{Asymptotic null distribution}\label{asymptotic null distribution}
In this section, we derive an asymptotic distribution of $\widehat{T}^2_{n,m}$ under $H_0$. 
In what follows, the symbol $``\xrightarrow{\mathcal{D}}"$ designates convergence in distribution.
\begin{Thm}\label{Asymptotic_Null_Distribution}
	Suppose that Condition is satisfied.
	Then, under $H_0: P=Q$, as $n,m \to \infty$,
	\[
	(n+m)\widehat{T}^2_{n,m} \xrightarrow{\mathcal{D}} \frac{1}{\rho(1-\rho)} \sum_{\ell=1}^{\infty} \lambda_{\ell} Z_{\ell}^2,
	\]
	where $Z_{\ell} \overset{i.i.d.}{\sim}N(0,1)$ and $\lambda_{\ell}$ is an eigenvalue of $V_{X \sim P}[(k(\cdot,X)-\mu_k(P))^{\otimes 2}]$.
\end{Thm}
It is generally not easy to utilize such an asymptotic null distribution because it is an infinite sum and determining the weights included in the asymptotic distribution is itself a difficult problem.
For practical purposes, a method to approximate the distribution of the test by $\widehat{T}_{n,m}^2$ under $H_0$ is developed in Section \ref{implementation}.
The method is based on the eigenvalues of the centered Gram matrices associated with the data set in Section \ref{Gram matrix spectrum}.

\subsection{Asymptotic non-null distribution}\label{asymptotic nonnull distribution}
In this section, an asymptotic distribution of $\widehat{T}^2_{n,m}$ under $H_1$ is investigated.
\begin{Thm}\label{Asymptotic_Nonnull_Distribution}
	Suppose that Condition is satisfied.
	Then, under $H_1: P\neq Q$, as $n,m \to \infty$,
	\[
	\sqrt{n+m}(\widehat{T}^2_{n,m}-T^2) \xrightarrow{\mathcal{D}} N(0,4\rho^{-1}v^2_P+4(1-\rho)^{-1}v^2_Q),
	\]
	where 
	\[
	v_P^2= V_{X \sim P}\left[\left< \Sigma_k(P)-\Sigma_k(Q), (k(\cdot,X)-\mu_k(P))^{\otimes 2}-\Sigma_k(P)
	\right>_{H(k)^{\otimes 2}} \right]\\
	\]
	and
	\[
	v_Q^2=V_{Y \sim Q}\left[\left< \Sigma_k(P)-\Sigma_k(Q), (k(\cdot,Y)-\mu_k(Q))^{\otimes 2}-\Sigma_k(Q)
	\right>_{H(k)^{\otimes 2}}\right].
	\]
\end{Thm}
\begin{Rem}
	We see by Theorem \ref{Asymptotic_Nonnull_Distribution} that
	\[
	\frac{\sqrt{n+m}(\widehat{T}^2_{n,m}-T^2)}{v} \xrightarrow{\mathcal{D}} N(0,1),
	\]
	where $v=\sqrt{4\rho^{-1} v^2_{P}+4 (1-\rho)^{-1} v^2_{Q}}$.
	Thus, we can evaluate the power of the test by $(n+m)\widehat{T}^2_{n,m}$ as
	\begin{align*}
	\Pr((n+m) \widehat{T}^2_{n,m} \geq t_{\alpha}|H_1)
	&=\Pr((n+m)(\widehat{T}^2_{n,m}-T^2) \geq t_{\alpha}-(n+m)T^2|H_1)\\
	&=\Pr\left(
	\frac{\sqrt{n+m}(\widehat{T}^2_{n,m}-T^2)}{v} \geq \frac{t_{\alpha}}{\sqrt{n+m} v}-\frac{\sqrt{n+m} T^2}{v} \Bigg| H_1
	\right)\\
	&\approx 1-\Phi\left(
	\frac{t_{\alpha}}{\sqrt{n+m}v}-\frac{\sqrt{n+m}T^2}{v} 
	\right) \to 1
	\end{align*}
	as $n,m \to \infty$, where $t_{\alpha}$ is the $(1-\alpha)$-quantile of the distribution of  $(n+m)\widehat{T}^2_{n,m}$ under $H_0$, and $\Phi$ is the distribution function of the standard normal distribution.
	Therefore, this test is consistent.
\end{Rem}
\subsection{Asymptotic distribution under contiguous alternatives}\label{asymptotic distribution under contiguous alternatives}
In this section, we develop an asymptotic distribution of $\widehat{T}^2_{n,m}$ under a sequence of local alternative distributions $Q_{nm}=(1-1/\sqrt{n+m})P+(1/\sqrt{n+m})Q$ for $P \neq Q$.
\begin{Thm}\label{Asymptotic distribution under contiguous alternatives}
	Let $X_1,\dots,X_n \overset{i.i.d.}{\sim} P$ and $Y_1,\dots,Y_m \overset{i.i.d.}{\sim} Q_{nm}$.
	Suppose that Condition is satisfied.
	Let $h: \mathcal{H} \times \mathcal{H} \to \mathbb{R}$ be a kernel defined as
	\begin{equation}\label{Eq_h(x,y)}
		h(x,y)=\left<(k(\cdot,x)-\mu_k(P))^{\otimes 2}-\Sigma_k(P),(k(\cdot,y)-\mu_k(P))^{\otimes 2}-\Sigma_k(P)\right>_{H(k)^{\otimes 2}},~~x,y \in \mathcal{H}
	\end{equation}
	and
	\[
	h(\cdot,x)= (k(\cdot,x)-\mu_k(P))^{\otimes 2}-\Sigma_k(P) \in H(k)^{\otimes 2} \nonumber
	\]
	and let $S_k:L_2(\mathcal{H},P) \to L_2(\mathcal{H},P)$ be a self-adjoint operator defined as	
	\begin{equation}\label{Eq_S_k}
		S_kg(x)= \int_{\mathcal{H}} h(x,y)g(y) dP(y),~~g \in L_2(\mathcal{H},P)
	\end{equation}
	(see Sections VI.1, VI.3, and VI.6 in \cite{reed1981functional} for details). 
	Then, as $n,m \to \infty$
		\[
		(n+m)\norm{\Sigma_k(\widehat{P})-\Sigma_k(\widehat{Q}_{nm})}^2_{H(k)^{\otimes 2}} \xrightarrow{\mathcal{D}}  \frac{1}{\rho(1-\rho)}\sum_{\ell=1}^{\infty} \theta_{\ell}W_{\ell}^2,	
		\]
		where $W_{\ell} \sim N(\sqrt{\rho (1-\rho)}\cdot \zeta_{\ell}(P,Q)/\theta_{\ell},1),~W_{\ell} \indepe W_{\ell'} ~(\ell \neq \ell')$,
		\[
		\zeta_{\ell}(P,Q)=\int_{\mathcal{H}}\left<\Sigma_k(Q)-\Sigma_k(P)+(\mu_k(Q)-\mu_k(P))^{\otimes 2},h(\cdot,y)\right>_{H(k)^{\otimes 2}} \Psi_{\ell}(y)dP(y),
		\]
		and $\theta_{\ell}$ and $\Psi_{\ell}$ are, respectively, the eigenvalue of $S_k$ and the eigenfunction corresponding to $\theta_{\ell}$.
\end{Thm}
	The following proposition claims that the eigenvalues $\theta_{\ell}$ appearing in Theorem \ref{Asymptotic distribution under contiguous alternatives} are the same as the eigenvalues $\lambda_{\ell}$ appearing in Theorem \ref{Asymptotic_Null_Distribution}:
\begin{Prop}\label{Prop_eigenvalue_T1_and_T3}
	The eigenvalues of $V_{X \sim P}[(k(\cdot,X)-\mu_k(P))^{\otimes 2}]$ in Theorem \ref{Asymptotic_Null_Distribution} and $S_k$ in (\ref{Eq_S_k}) of Theorem \ref{Asymptotic distribution under contiguous alternatives} are the same.
\end{Prop}
From Theorems \ref{Asymptotic_Null_Distribution} and \ref{Asymptotic distribution under contiguous alternatives} and Proposition \ref{Prop_eigenvalue_T1_and_T3}, it can be seen that the local power of the test by $(n+m)\widehat{T}^2_{n,m}$ is dominated by the noncentrality parameters.
It follows that
\[
\zeta_{\ell} (P,Q) =\int_{\mathcal{H}} \left<\mathbb{E}_{Y \sim Q} [h(\cdot,y)] ,h(\cdot,x) \right> \Psi_{\ell}(y) dP(y)=\lambda_{\ell} \mathbb{E}_{Y \sim Q} [\Psi_{\ell}(Y)]
\]
by which we obtain 
\begin{align*}
\mathbb{E}\left[\frac{1}{\rho(1-\rho)} \sum_{\ell=1}^{\infty} \theta_{\ell} W_{\ell}^2\right] 
&=\frac{1}{\rho(1-\rho)}\sum_{\ell=1}^{\infty} \lambda_{\ell}  \left(1+\rho(1-\rho) \cdot \frac{\zeta_{\ell}(P,Q)^2}{\lambda_{\ell}^2}\right)\\
&= \frac{1}{\rho(1-\rho)}\sum_{\ell=1}^{\infty} \lambda_{\ell} \Big(1+\rho(1-\rho) \{\mathbb{E}_{Y \sim Q} [\Psi_{\ell}(Y)]\}^2\Big).
\end{align*}
In addition,  from Theorem 1 in \cite{Minh2006},  we have
\begin{align*}
\sum_{\ell=1}^{\infty} \lambda_{\ell}  \{\mathbb{E}_{Y \sim Q}[\Psi_{\ell}(Y)]\}^2
&=\sum_{\ell=1}^{\infty} \lambda_{\ell} \mathbb{E}_{Y \sim Q} [\Psi_{\ell}(Y)] \mathbb{E}_{Y' \sim Q} [\Psi_{\ell}(Y')]\\
&=\mathbb{E}_{Y,Y' \sim Q} [h(Y,Y')]\\
&=\norm{\Sigma_k(Q)-\Sigma_k(P)+(\mu_k(Q)-\mu_k(P))^{\otimes 2}}_{H(k)^{\otimes 2}}^2.
\end{align*}
Hence, the local power reveals not only the difference between $\Sigma_k(P)$ and $\Sigma_k(Q)$ but also that between $\mu_k(Q)$ and $\mu_k(P)$.
\subsection{Null distribution estimates using Gram matrix spectrum}\label{Gram matrix spectrum}
The asymptotic null distribution was obtained in Theorem \ref{Asymptotic_Null_Distribution}, but it is difficult to derive its weights.
The following theorem  shows that this weight can be estimated using the estimator of $V[(k(\cdot,X)-\mu_k(P))^{\otimes 2}]$.

\begin{Thm}\label{Asymptotic_distribution_Gram_matrix_spectrum}
Assume that $\mathbb{E}_{X \sim P}[k(X,X)^2]< \infty$.
Let $\{\lambda_{\ell}\}_{\ell=1}^{\infty}$ and $\{\widehat{\lambda}^{(n)}_{\ell}\}_{\ell=1}^{\infty}$ be the eigenvalues of $\Upsilon$ and $\widehat{\Upsilon}^{(n)}$, respectively, where 
\[
\Upsilon=V\left[(k(\cdot,X)-\mu_k(P))^{\otimes 2}\right]~~\text{and}~~\widehat{\Upsilon}^{(n)}=\frac{1}{n}\sum_{i=1}^{n} \left\{
(k(\cdot,X_i)-\mu_k(\widehat{P}))^{\otimes 2}-\Sigma_k(\widehat{P})\right\}^{\otimes 2}.
\]
Then, as $n \to \infty$
\[
\sum_{\ell=1}^{\infty} \widehat{\lambda}_{\ell}^{(n)} Z_{\ell}^2 \xrightarrow{\mathcal{D}} \sum_{\ell=1}^{\infty} \lambda_{\ell} Z_{\ell}^2,
\]
where $Z_{\ell} \overset{i.i.d.}{\sim} N(0,1)$.
\end{Thm}
In addition, Proposition \ref{Lem_eigenvaue_same} claims the eigenvalues of $\widehat{\Upsilon}^{(n)}$ and the Gram matrix are the same.
\begin{Prop}\label{Lem_eigenvaue_same}
	The $n \times n$ Gram matrix $H=(H_{ij})_{1 \leq i,j \leq n}$ is defined as 
	\[
	H_{ij}=\left<(k(\cdot,X_i)-\mu_k(\widehat{P}))^{\otimes 2} -\Sigma_k(\widehat{P}),(k(\cdot,X_j)-\mu_k(\widehat{P}))^{\otimes 2} -\Sigma_k(\widehat{P})\right>_{H(k)^{\otimes 2}}.
	\]
	Then, the eigenvalues of $\widehat{\Upsilon}^{(n)}$ and $H/n$ are the same.
\end{Prop}
\begin{Rem}
	By Proposition \ref{Lem_eigenvaue_same}, the critical value can be obtained by calculating $1/\{\rho(1-\rho)\}\sum_{\ell=1}^{n-1} \widehat{\lambda}_{\ell}^{(n)} Z_{\ell}^2$ using the eigenvalue $\widehat{\lambda}^{(n)}_{\ell},~\ell \in \{1,\dots,n-1\}$ of $H/n$.
	In addition, the matrix $H$ is expressed as 
	\begin{equation}\label{H_matrix_Gram}
	H=P_n(\widetilde{K}_{X}\odot\widetilde{K}_{X})P_n,
	\end{equation}
	where $\odot$ is the Hadamard product.
	 The $n \times n$ Gram matrix $K_X$ is a positive definite, but $H$ has eigenvalue 0 since $H \underline{1}=\underline{0}$.
\end{Rem}
\section{Implementation}\label{implementation}
This section proposes corrections to the asymptotic distribution for both the MVD and MMD tests, and describes the results of simulations of the type-I error and power for the modifications.
The MMD test is a two-sample test for $H_0$ and $H_1$ using the test statistic:
\[
\widehat{\Delta}_{n,m}^2=\norm{\frac{1}{n}\sum_{i=1}^{n}k(\cdot,X_i)-\frac{1}{m}\sum_{j=1}^{m}k(\cdot,Y_j)}_{H(k)}^2.
\]
The $\widehat{\Delta}^2_{n,m}$ of the asymptotic null distribution is the infinite sum of the weighted chi-square distribution, which is the  same as $\widehat{T}^2_{n,m}$ in (\ref{Eq_widehat_Tnm}).
The approximate distribution can be obtained by estimating the eigenvalues by the centered Gram matrix (see \cite{Gretton2009a} for details).
\subsection{Approximation of the null distribution} \label{approximation of the null distribution}
In this section, we discuss methods to approximate the null distributions of the MVD and MMD tests.
The asymptotic null distribution of the MVD test was obtained in Theorem \ref{Asymptotic_Null_Distribution} as an infinite sum of weighted chi-squared random variables with one degree of freedom, and according to Theorem \ref{Asymptotic_distribution_Gram_matrix_spectrum}, those weights $\lambda_{\ell}~(\ell \geq 1)$ can be estimated by the eigenvalues $\widehat{\lambda}_{\ell}^{(n)}$ of the matrix $H/n$.  
Similar results was obtained for MMD by \cite{Gretton2009a}.
However, this approximate distribution based on this estimated eigenvalue does not actually work well.
In fact, by comparing the simulated exact null distribution with this approximate distribution based on estimated eigenvalues, it can be seen that variance of the approximate distribution is larger than that of the simulated exact null distribution.
We modify the approximate distribution $1/\{\rho(1-\rho)\}\sum_{\ell=1}^{n-1} \widehat{\lambda}_{\ell}^{(n)} Z_{\ell}^2$ by obtaining the variance of this simulated exact null distribution.
The variance of the exact null distribution $V[(n+m)\widehat{T}^2_{n,m}]$ is obtained as the following proposition:
\begin{Prop}\label{Variance_MVD_P=Q}
	Assume that $\mathbb{E}_{X \sim P}[k(X,X)^2] < \infty$.
	Then under $H_0$,
	\[
	V[(n+m)\widehat{T}^2_{n,m}]=\frac{2(n+m)^4}{n^2 m^2} \norm{\Upsilon}^2_{H(k)^{\otimes 4}} +O\left(\frac{1}{n}\right)+O\left(\frac{1}{m}\right).
	\]
\end{Prop}
Proposition \ref{Variance_MVD_P=Q} leads to 
\begin{equation}\label{Eq_nm_kl_MVD}
V[(n+m)\widehat{T}^2_{n,m}] \approx \frac{(n+m)^4}{n^2m^2} \cdot \frac{k^2 \ell^2}{(k+\ell)^4} V[(k+\ell)\widehat{T}^2_{k,\ell}],
\end{equation}
with respect to $V[(n+m)\widehat{T}^2_{n,m}]$ and $V[(k+\ell)\widehat{T}^2_{k,\ell}],~k,\ell \in \mathbb{N}$.
If we can estimate the variance $V[(k+\ell) \widehat{T}_{k,\ell}^2]$ at $k$ and $\ell$ that is less than $n$ and $m$, respectively, we can estimate $V[(n+m)\widehat{T}^2_{n,m}]$ by using (\ref{Eq_nm_kl_MVD}).
In addition, the method of estimating $V[(k+\ell)\widehat{T}_{k,\ell}^2]$ by choosing $(k,\ell)$ from $(n,m)$ without replacement is known as subsampling.

The following proposition for MMD shows a similar result to MVD:
\begin{Prop}\label{Variance_MMD_P=Q}
	Assume that $\mathbb{E}_{X \sim P}[k(X,X)] < \infty$.
	Then, under $H_0$
	\[
	V[(n+m)\widehat{\Delta}^2_{n,m}]=\frac{2(n+m)^4}{n^2m^2}\norm{\Sigma_k(P)}^2_{H(k)^{\otimes 2}}+O\left(
	\frac{1}{n}\right)
	+O\left(\frac{1}{m}\right).
	\]
\end{Prop}
\subsubsection{Subsampling method}\label{Bootstrap_method}
We used the subsampling method to estimate $V[(k+\ell)\widehat{T}^2_{k,\ell}]$ and $V[(k+\ell)\widehat{\Delta}^2_{k,\ell}]$ (see Section 2.2 in \cite{politis1999} for details).
In order to obtain two samples under the null hypothesis, we divide $X_1, \dots , X_n$ into $X_1, \dots, X_{n_1}$ and $X_{n_1+1},\dots, X_n$.
Then, we randomly select $X^*_1(i), \dots, X^*_k(i)$ and $Y^*_1(i), \dots, Y^*_\ell(i)$ from $X_1, \dots, X_{n_1}$ and $X_{n_1+1},\dots, X_n$ without replacement, which repaet in each iteration $i \in \{1,\dots, I\}$.
These randomly selected values generate the replications of the test statistic
\[
\widehat{T}_{k,\ell}^2(i)=  \widehat{T}_{k,\ell}^2(X^*_1(i),\dots,X^*_k(i); Y^*_1(i),\dots,Y^*_\ell(i))
\]
for iterations $i \in \{1,\dots, I\}$.
The generated test statistics $(k+\ell) \widehat{T}_{k,\ell}^2(i)$ in $I$ iterations estimate $V[(k+\ell)\widehat{T}^2_{k,\ell}]$ by calculating the unbiased sample variance:
\[
V[(k+\ell)\widehat{T}^2_{k,\ell}]_{\text{sub}}=\frac{1}{I-1}\sum_{j=1}^{I} \left\{
(k+\ell){\widehat{T}^2}_{k,\ell}(j)-(k+\ell)\overline{\widehat{T}^2}_{k,\ell}
\right\}^2,
\]
where $\overline{\widehat{T}^2}_{k,\ell}=(1/I) \sum_{i=1}^{I} {\widehat{T}^2}_{k,\ell}(i)$.
According to (\ref{Eq_nm_kl_MVD}), $V[(n+m)\widehat{T}^2_{n,m}]$ is estimated by 
\begin{equation}\label{Eq_MVD_V_nm}
V[(n+m)\widehat{T}^2_{n,m}]_{\text{sub}}=\frac{(n+m)^4}{n^2 m^2} \frac{k^2 \ell^2}{(k+\ell)^4}V[(k+\ell)]\widehat{T}^2_{k,\ell}]_{\text{sub}}.
\end{equation}
We also estimate $V[(n+m)\widehat{\Delta}^2_{n,m}]$ by using 
\begin{equation}\label{Eq_MMD_V_nm}
V[(n+m)\widehat{\Delta}^2_{n,m}]_{\text{sub}}=\frac{(n+m)^4}{n^2 m^2} \frac{k^2 \ell^2}{(k+\ell)^4}V[(k+\ell)\widehat{\Delta}^2_{k,\ell}]_{\text{sub}}
\end{equation}
from Proposition \ref{Variance_MMD_P=Q}, where 
\[
\widehat{\Delta}_{k,\ell}^2(i)=  \widehat{\Delta}_{k,\ell}^2(X^*_1(i),\dots,X^*_k(i); Y^*_1(i),\dots,Y^*_\ell(i))
\]
for $i \in \{1,\dots,I\}$,
\[
V[(k+\ell)\widehat{\Delta}^2_{k,\ell}]_{\text{sub}}=\frac{1}{I-1}\sum_{j=1}^{I} \left\{
(k+\ell){\widehat{\Delta}^2}_{k,\ell}(j)-(k+\ell)\overline{\widehat{\Delta}^2}_{k,\ell}
\right\}^2
\]
and $\overline{\widehat{\Delta}^2}_{k,\ell}=(1/I) \sum_{i=1}^{I} {\widehat{\Delta}^2}_{k,\ell}(i)$.

The columns of $(n+m) \widehat{T}^2_{n,m}$ in Table \ref{MVD_variance} and $(n+m) \widehat{\Delta}^2_{n,m}$ in Table \ref{MMD_variance} are variances of $(n+m) \widehat{T}^2_{n,m}$ and $(n+m) \widehat{\Delta}^2_{n,m}$, which are estimated by a simulation of 10,000 iterations with $X_1,\dots,X_n \overset{i.i.d.}{\sim} N(\underline{0},I_d)$ and $Y_1,\dots,Y_m \overset{i.i.d.}{\sim} N(\underline{0},I_d)$ for each $\sigma, d$ and $(n,m)$.
The subsampling variances $V[(n+m)\widehat{T}^2_{n,m}]_{\text{sub}}$ in (\ref{Eq_MVD_V_nm}) and $V[(n+m)\widehat{\Delta}^2_{n,m}]_{\text{sub}}$ in (\ref{Eq_MMD_V_nm}) with $X_1,\dots,X_n \overset{i.i.d.}{\sim} N(\underline{0},I_d)$ are given in the columns labeled ``Subsampling" for $(k,\ell)$.
Tables \ref{MVD_variance} and \ref{MMD_variance} show that subsampling variances $V[(n+m)\widehat{T}^2_{n,m}]_{\text{sub}}$ and $V[(n+m)\widehat{\Delta}^2_{n,m}]_{\text{sub}}$ estimate the exact variances well.
However, these variances tend to underestimate the exact variances.

\begin{figure}[H]
	\begin{minipage}{0.5\hsize}
		\centering
		\includegraphics[width=1\linewidth]{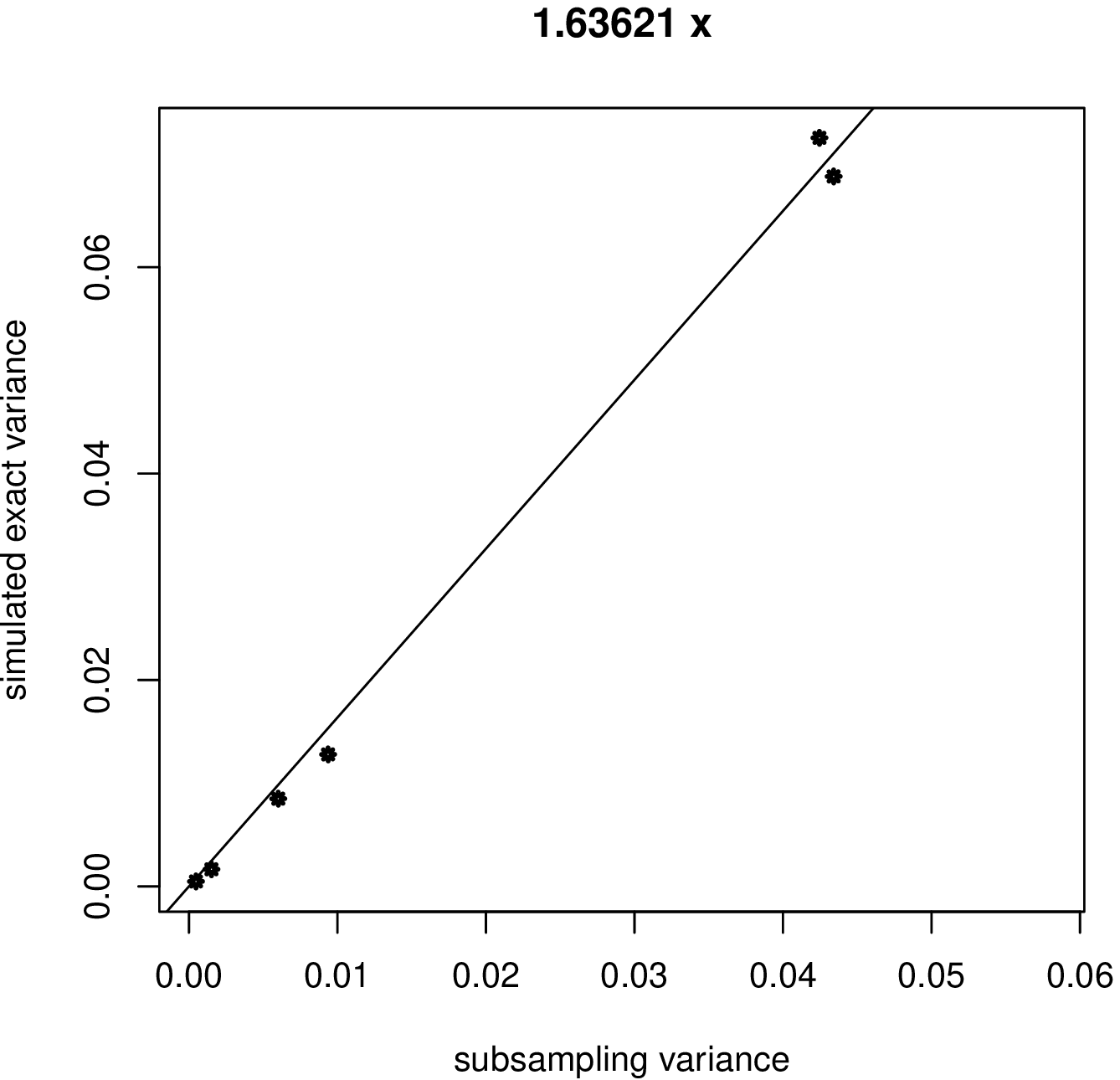}
		\vspace{-1cm}
		\caption*{(a)}
	\end{minipage}	
	\begin{minipage}{0.5\hsize}
		\centering
		\includegraphics[width=1\linewidth]{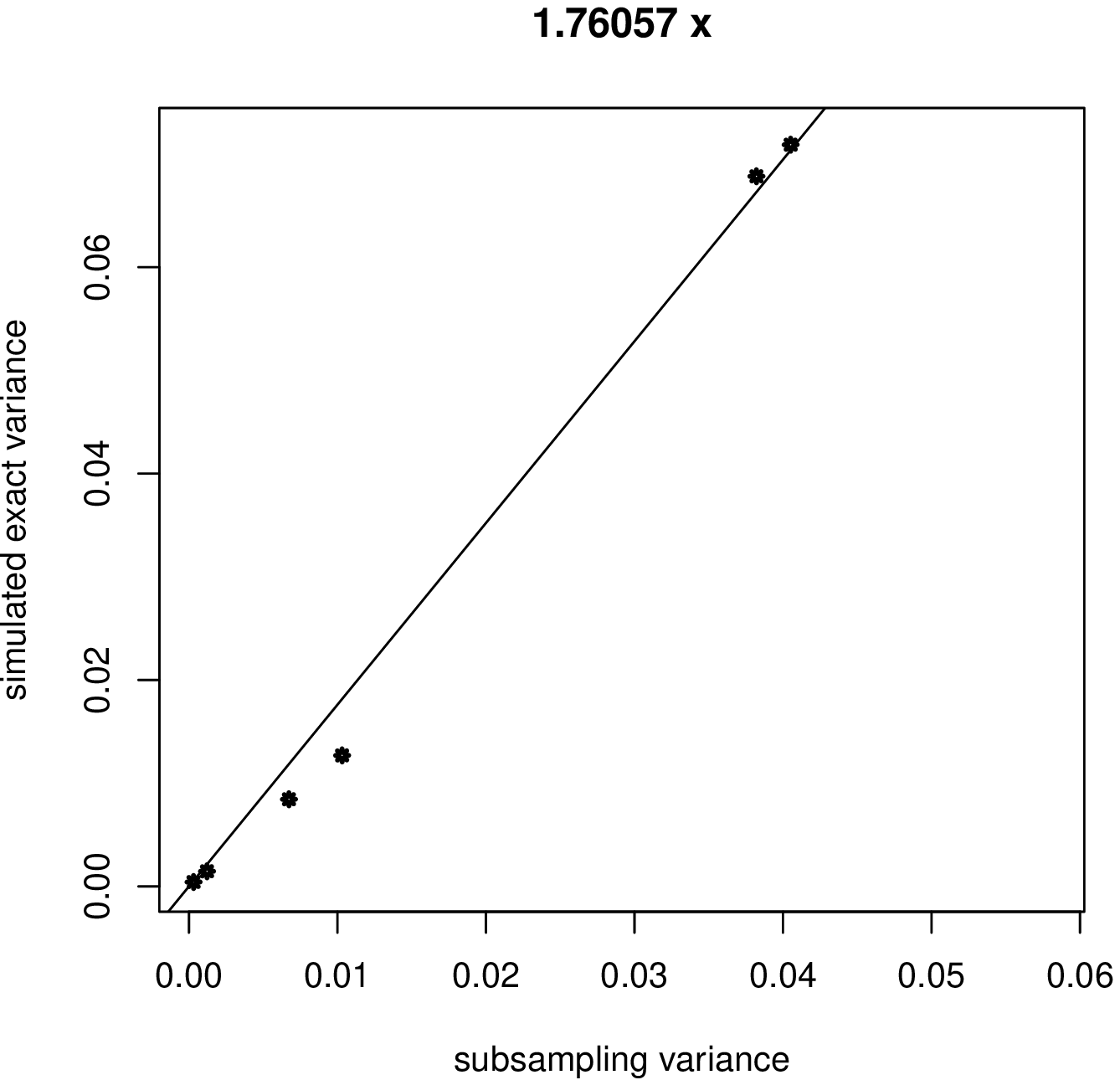}
		\vspace{-1cm}
		\caption*{(b)}
	\end{minipage}	
	\begin{minipage}{0.5\hsize}
		\centering
		\includegraphics[width=1\linewidth]{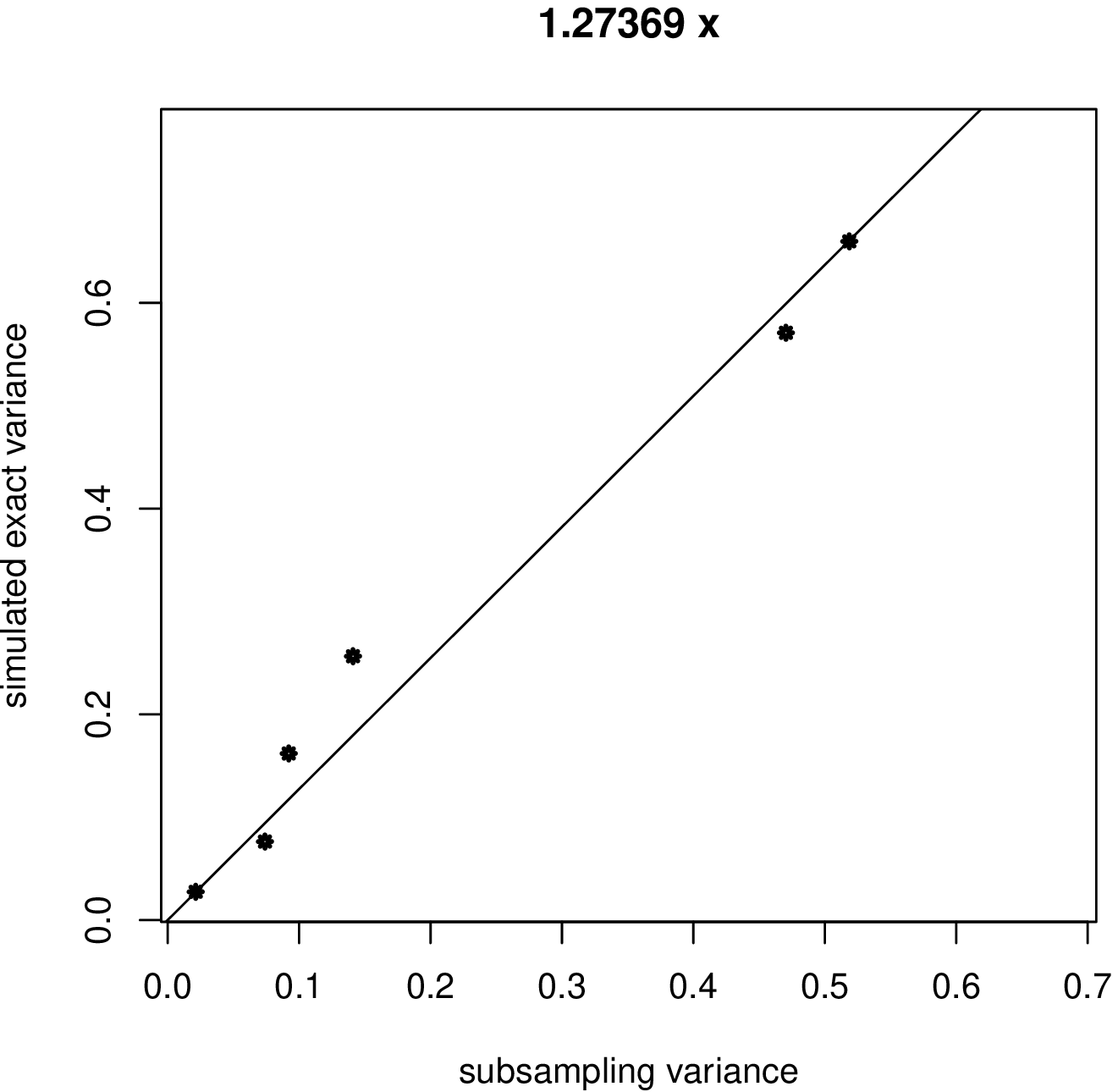}
		\vspace{-1cm}
		\caption*{(c)}
	\end{minipage}	
	\begin{minipage}{0.5\hsize}
		\centering
		\includegraphics[width=1\linewidth]{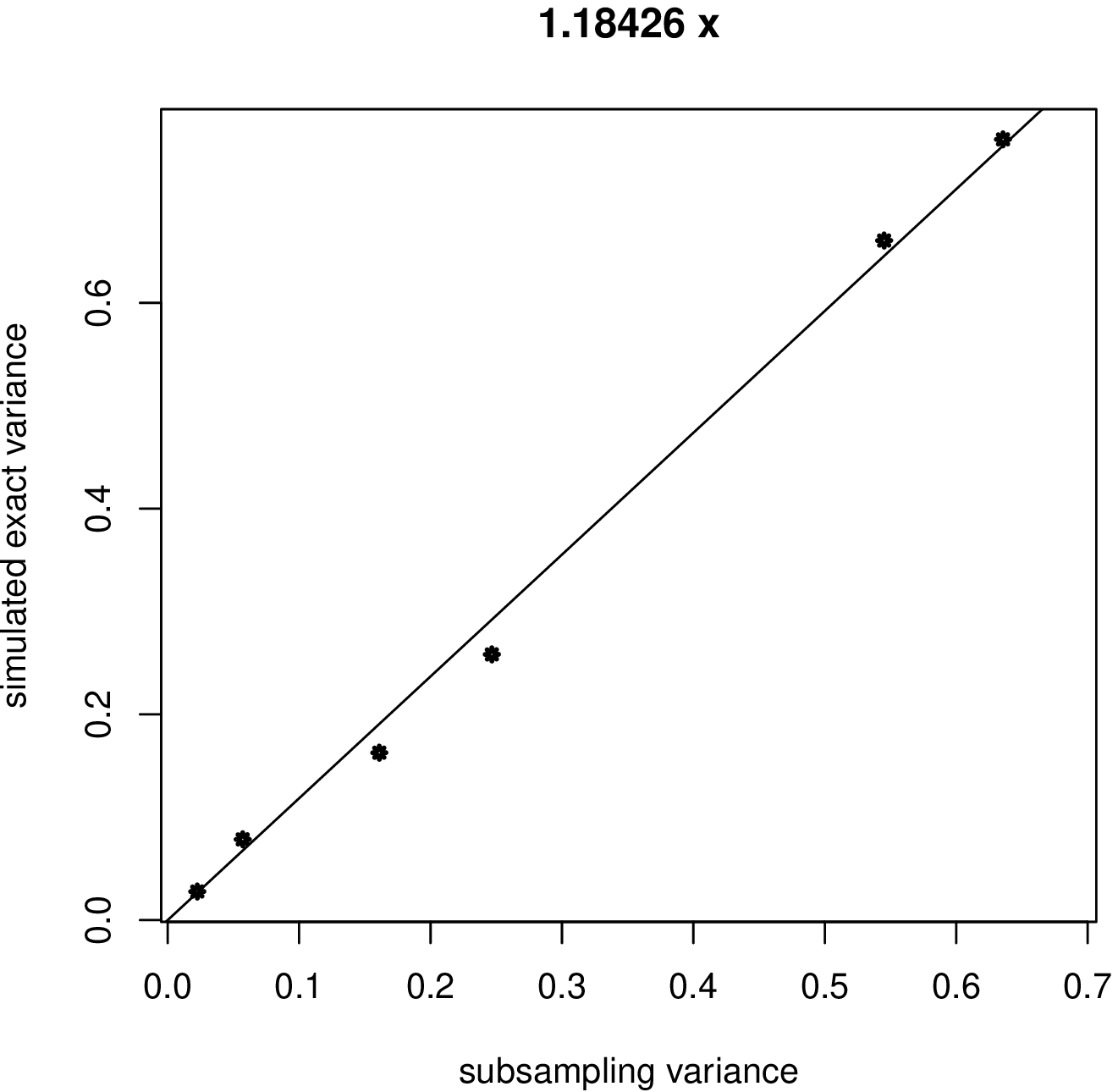}
		\vspace{-1cm}
		\caption*{(d)}
	\end{minipage}	
  \caption{The results of linear regression of the simulated variance of $(n+m)\widehat{T}^2_{n,m}$ and the subsampling variance $V[(n+m)\widehat{\Delta}^2_{n,m}]_{\text{sub}}$, and the results for the MMD.
  (a) MVD $(n, m) = (200,200),~ (k,\ell)=(50,50)$.
  (b) MVD $(n, m) = (500,500),~(k,\ell)=(125,125)$.
  (c) MMD $(n, m) = (200,200),~ (k,\ell)=(50,50)$.
  (d) MMD $(n, m) = (500,500),~(k,\ell)=(125,125)$.
  }
  \label{MVD_slope}
\end{figure}

We investigate how much smaller $V[(n+m)\widehat{T}^2_{n,m}]$ is than the variance of $(n+m)\widehat{T}^2_{n,m}$ by performing a linear regression of $V[(n+m)\widehat{T}^2_{n,m}]$ and the variance of $(n+m)\widehat{T}^2_{n,m}$ with an intercept equal to 0, with the same for the MMD.
The results are shown in Figure \ref{MVD_slope}; (a) and (b) show results for the MVD and (c) and (d) show results for the MMD, which are respectively $(n,m)=(200,200)$ and $(k,\ell)=(50,50)$ and $(n,m)=(500,500)$ and $(k,\ell)=(125,125)$ cases.
In Figure \ref{MVD_slope}, the $x$ axis is $V[(n+m)\widehat{T}^2_{n,m}]_{\text{sub}}$ or $V[(n+m)\widehat{\Delta}^2_{n,m}]_{\text{sub}}$ and the $y$ axis is variance of $(n+m)\widehat{T}^2_{n,m}$ or $(n+m)\widehat{\Delta}^2_{n,m}$ for each $\sigma,~d~,(n, m)$ and $(k, \ell)$ in Table \ref{MVD_variance} or Table \ref{MMD_variance}.
The line in Figure \ref{MVD_slope} is a regression line found by the least-squares method in the form $y = ax + \varepsilon$.
It can be seen from Figure \ref{MVD_slope} that when $V[(n+m)\widehat{T}^2_{n,m}]_{\text{sub}}$ and $V[(n+m)\widehat{\Delta}^2_{n,m}]_{\text{sub}}$ are multiplied by the term associated regression coefficient, they approach the variances $(n+m)\widehat{T}^2_{n,m}$ and $(n+m)\widehat{\Delta}^2_{n,m}$.
The coefficient of linear regression  with intercept 0 is written in the row labeled ``slope of the line" in Tables \ref{MVD_variance} and \ref{MMD_variance}.

\begin{table}[H]
	\caption{The variance of $(n+m)\widehat{T}^2_{n,m}$ under $P=Q=N(\underline{0},I_d)$ and $V[(n+m)\widehat{T}^2_{n,m}]_{\text{sub}}$ : $I = 1,000$, $n_1=n/2$, and $X_1,\dots,X_n \overset{i.i.d.}{\sim} N(\underline{0},I_d)$.
	}
	\label{MVD_variance}
	\centering
	\begin{tabular}{ccc|c|ccc} 
		$\sigma$&$d$&$(n,m)$&{$(n+m)\widehat{T}^2_{n,m}$} &\multicolumn{3}{c}{Subsampling $(k,\ell)$}\\ 
		\hline 
		\multicolumn{3}{c|}{ }&&$(n/4,n/4)$&$(n/6,n/6)$&$(n/8,n/8)$\\ 
		\hline 
		$d^{-3/4}$&5&(200,200)& 0.06880 &0.04341&0.05168&0.04902\\ 
		$d^{-3/4}$&5&(500,500)& 0.06881&0.03821&0.04897&0.04921\\ 
		$d^{-7/8}$&5&(200,200)& 0.07254&0.04246&0.05138&0.05798\\ 
		$d^{-7/8}$&5&(500,500)& 0.07188&0.04052&0.05500&0.05593\\ 
		$d^{-3/4}$&10&(200,200)&0.00850&0.00602&0.00812&0.00898\\ 
		$d^{-3/4}$&10&(500,500)&0.00845&0.00674&0.00751&0.00753\\ 
		$d^{-7/8}$&10&(200,200)&0.01280 &0.00937&0.01224&0.01377\\ 
		$d^{-7/8}$&10&(500,500)& 0.01270&0.01032&0.01251&0.01255\\ 
		$d^{-3/4}$&20&(200,200)& 0.00049&0.00048&0.00070&0.00094\\ 
		$d^{-3/4}$&20&(500,500)& 0.00043&0.00031&0.00046&0.00060\\ 
		$d^{-7/8}$&20&(200,200)& 0.00166&0.00152&0.00261&0.00330\\ 
		$d^{-7/8}$&20&(500,500)&0.00147&0.00122&0.00165&0.00204\\ 
		\hline
		\multicolumn{2}{c}{ }&(200,200)&1&1.63621&1.35769  &1.29601  \\
		\multicolumn{2}{c}{slope of the line}&(500,500)&1&
		1.76057 &1.33845 &1.3232 \\
		\multicolumn{2}{c}{}&both&1&1.69348&1.34798  &1.30928   
	\end{tabular} 
\end{table}
\begin{table}[H]
	\caption{The variance of $(n+m)\widehat{\Delta}^2_{n,m}$ under $P=Q=N(\underline{0},I_d)$ and $V[(n+m)\widehat{\Delta}^2_{n,m}]_{\text{sub}}$ : $I = 1,000$, $n_1=n/2$, and $X_1,\dots,X_n \overset{i.i.d.}{\sim} N(\underline{0},I_d)$.}
	\label{MMD_variance}
	\centering
	\begin{tabular}{ccc|c|ccc} 
		$\sigma$&$d$&$(n,m)$&{$(n+m)\widehat{\Delta}^2_{n,m}$} &\multicolumn{3}{c}{Subsampling $(k,\ell)$}\\ 
		\hline 
		\multicolumn{3}{c|}{ }& &$(n/4,n/4)$&$(n/6,n/6)$&$(n/8,n/8)$\\ 
		\hline 
		$d^{-3/4}$&5&(200,200)& 0.57100 &0.47044&0.50315 &0.57047 \\ 
		$d^{-3/4}$&5&(500,500)&  0.66068 &0.54518& 0.63054&0.58848 \\ 
		$d^{-7/8}$&5&(200,200)&0.65987  & 0.51867 &0.57258 & 0.54567\\ 
		$d^{-7/8}$&5&(500,500)&0.75903&0.63563 &0.68024 &0.65349 \\ 
		$d^{-3/4}$&10&(200,200)& 0.16205 &0.09213 & 0.12017 & 0.12940 \\ 
		$d^{-3/4}$&10&(500,500)&0.16279 &0.16106 &0.16809 &0.18334 \\ 
		$d^{-7/8}$&10&(200,200)& 0.25656 &0.14104 &0.18435 &0.20457\\ 
		$d^{-7/8}$&10&(500,500)&0.25836&0.24670 &0.25567 &0.26402\\ 
		$d^{-3/4}$&20&(200,200)&0.02757&0.02135&0.02632 &0.02814\\ 
		$d^{-3/4}$&20&(500,500)&0.02784&0.02255&0.02407 &0.02615 \\ 
		$d^{-7/8}$&20&(200,200)&0.07642&0.07404& 0.08121& 0.08744\\ 
		$d^{-7/8}$&20&(500,500)&0.07856&0.05714&0.07492&0.07320\\ 
		\hline
		\multicolumn{2}{c}{ }&(200,200)&1&1.27369 &1.16037&1.11083   \\
		\multicolumn{2}{c}{slope of the line}&(500,500)&1& 1.18426  &1.07578 &1.12080\\
		\multicolumn{2}{c}{ }&both&1&1.21990  &1.10951  &1.11643
	\end{tabular} 
\end{table}

\begin{Rem}
	Since the subsampling variance $V[(k+\ell)\widehat{T}^2_{k,\ell}]_{\text{sub}}$ is an unbiased sample variance of $I$ times, we get
	\begin{align*}
	&\mathbb{E}[V[(k+\ell)\widehat{T}^2_{k,\ell}]_{\text{sub}}]\\
	&=\frac{1}{I-1}\sum_{j=1}^{I} \mathbb{E}[\left\{
	(k+\ell){\widehat{T}^2}_{k,\ell}(j)-(k+\ell)\overline{\widehat{T}^2}_{k,\ell}
	\right\}^2]\\
	&=\frac{1}{I-1}\sum_{j=1}^{I} \mathbb{E}\left[\left\{
	(k+\ell){\widehat{T}^2}_{k,\ell}(j)-\mathbb{E}[(k+\ell)\widehat{T}^2_{k,\ell}]+\mathbb{E}[(k+\ell)\widehat{T}^2_{k,\ell}]-(k+\ell)\overline{\widehat{T}^2}_{k,\ell}
	\right\}^2\right]\\
	&=\frac{1}{I-1}\sum_{j=1}^{I}\Big[
	\mathbb{E}\left[\left\{
	(k+\ell){\widehat{T}^2}_{k,\ell}(j)-\mathbb{E}[(k+\ell)\widehat{T}^2_{k,\ell}]
	\right\}^2\right]
	+\mathbb{E}\left[
	\left\{
	\mathbb{E}[(k+\ell)\widehat{T}^2_{k,\ell}]-(k+\ell)\overline{\widehat{T}^2}_{k,\ell}
	\right\}^2
	\right]\\
	&~~~~~+2\mathbb{E}\left[\left\{
	(k+\ell){\widehat{T}^2}_{k,\ell}(j)-\mathbb{E}[(k+\ell)\widehat{T}^2_{k,\ell}]
	\right\}
	\left\{
	\mathbb{E}[(k+\ell)\widehat{T}^2_{k,\ell}]-(k+\ell)\overline{\widehat{T}^2}_{k,\ell}
	\right\}
	\right]
	\Big] \\
	&=\frac{1}{I-1} \sum_{j=1}^{I}\Big\{
	V[(k+\ell)\widehat{T}^2_{k,\ell}]
	+\frac{1}{I^2} \sum_{i,s=1}^{I}Cov\left((k+\ell)\widehat{T}^2_{k,\ell}(i),(k+\ell)\widehat{T}^2_{k,\ell}(s)\right)\\
	&~~~~~
	-\frac{2}{I} \sum_{i=1}^{I}Cov\left((k+\ell)\widehat{T}^2_{k,\ell}(j),(k+\ell)\widehat{T}^2_{k,\ell}(i)\right)
	\Big\}\\
	&=\frac{I}{I-1} V[(k+\ell)\widehat{T}^2_{k,\ell}]-\frac{1}{I(I-1)}\sum_{i,j=1}^{I}Cov\left((k+\ell)\widehat{T}^2_{k,\ell}(i),(k+\ell)\widehat{T}^2_{k,\ell}(j)\right)\\
	&=V[(k+\ell)\widehat{T}^2_{k,\ell}]-\frac{2}{I(I-1)} \sum_{i < j}Cov\left((k+\ell)\widehat{T}^2_{k,\ell}(i),(k+\ell)\widehat{T}^2_{k,\ell}(j)\right).
	\end{align*}
	However, it is not easy to estimate $Cov\left((k+\ell)\widehat{T}^2_{k,\ell}(i),(k+\ell)\widehat{T}^2_{k,\ell}(j)\right)$. 
	This fact caused us to modify subsampling variances $V[(n+m)\widehat{T}^2_{n,m}]_{\text{sub}}$ in (\ref{Eq_MVD_V_nm}) and $V[(n+m)\widehat{\Delta}^2_{n,m}]_{\text{sub}}$ in (\ref{Eq_MMD_V_nm}) to $(1+\tau)V[(n+m)\widehat{T}^2_{n,m}]_{\text{sub}}$ and $(1+\tau)V[(n+m)\widehat{\Delta}^2_{n,m}]_{\text{sub}}$ using $\tau >0$, which is determined from the regression coefficient in the row ``slope of the line" in Tables \ref{MVD_variance} and \ref{MMD_variance}.
\end{Rem}

\subsubsection{Modification of approximation distribution}\label{Wdash}
\begin{figure}[H]
	\begin{minipage}{0.5\hsize}
		\centering
		\includegraphics[width=1\linewidth]{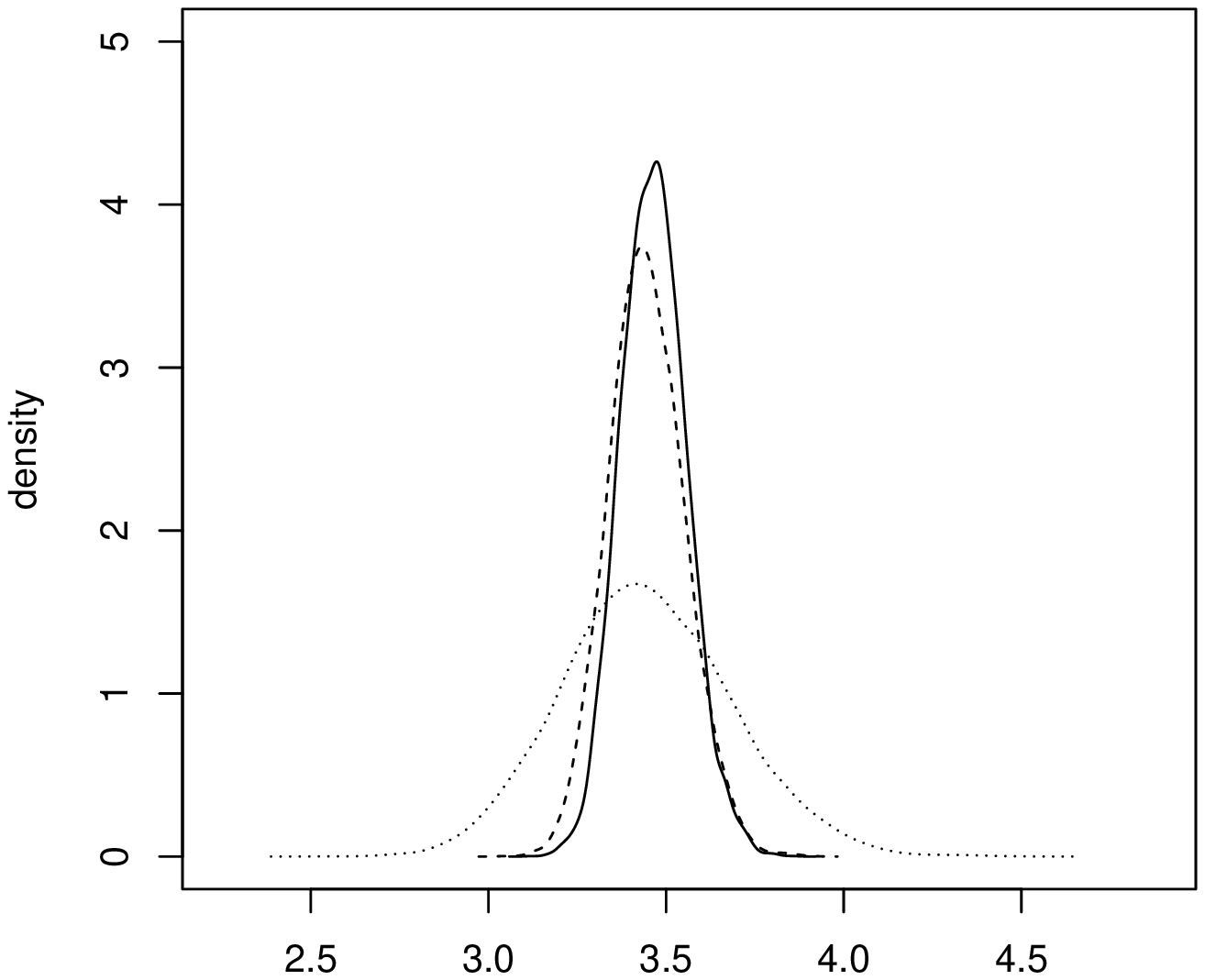}
	\end{minipage}	
	\begin{minipage}{0.5\hsize}
		\centering
		\includegraphics[width=1\linewidth]{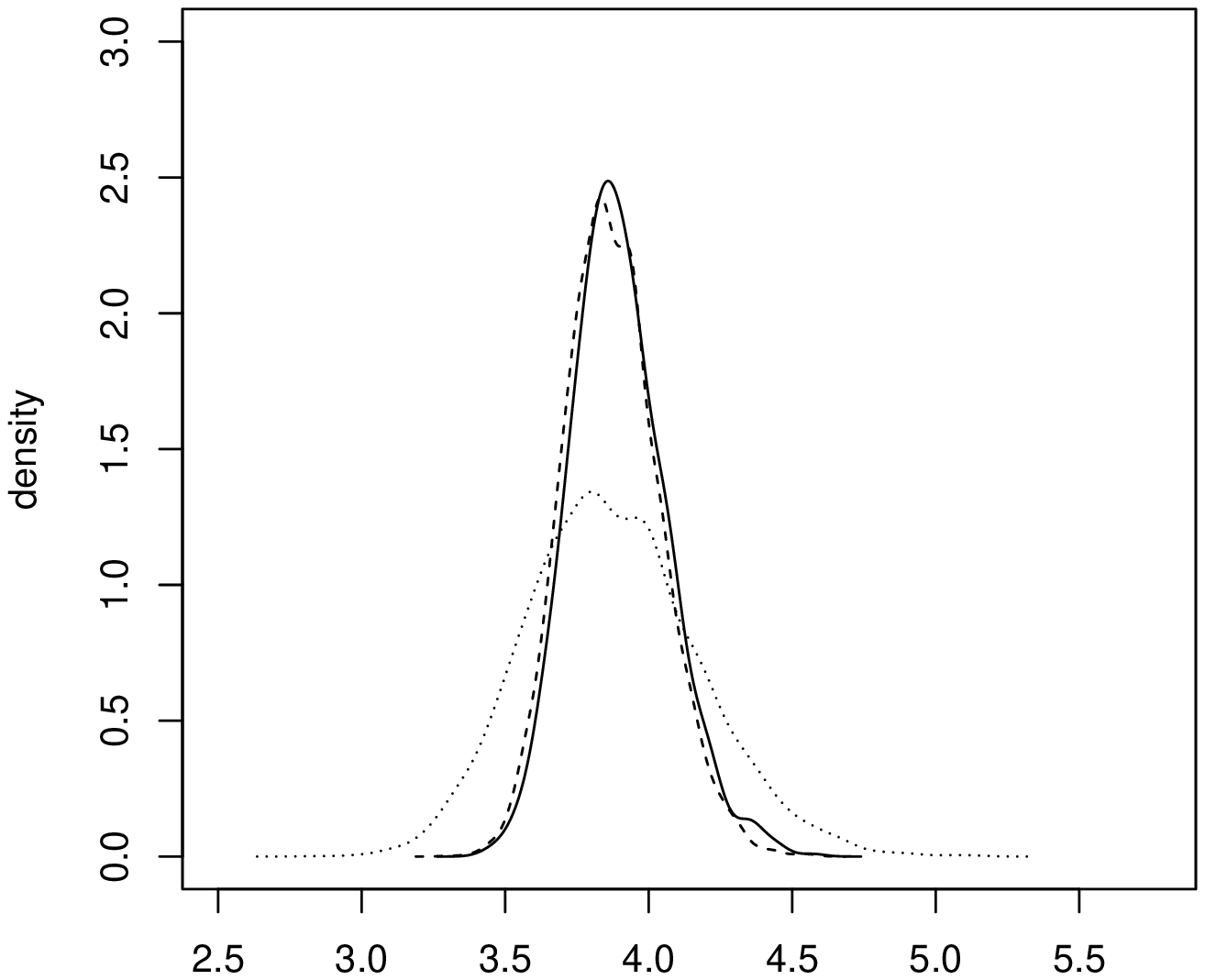}
	\end{minipage}	
	\caption{Density estimates of the distributions for simulated exact null distribution (solid), $W'_n$ (dashed), and $1/\{\rho(1-\rho)\} \sum_{\ell'=1}^{n-1} \widehat{\lambda}^{(n)}_{\ell} Z_{\ell}^2$ (dotted), with $(n,m)=(500, 500),~ \sigma =d^{-3/4}$, and $(k,\ell)=(125, 125)$.
	Left panel: Comparison of MVD results with $d=10,~\tau_{\text{MVD}}=0.69348$.
	Right panel: Comparison of MMD results with $d=20,~\tau_{\text{MMD}}=0.21990$.
	}
	\label{Wdash_comp}
\end{figure}

Since the subsampling variance $V[(n+m)\widehat{T}^2_{n,m}]_{\text{sub}}$ underestimates $V[(n+m)\widehat{T}^2_{n,m}]$, as seen in Table \ref{MVD_variance} and \ref{MMD_variance}, we use positive $\tau >0$ to estimate $V[(n+m)\widehat{T}^2_{n,m}]$ by $(1+\tau) V[(n+m)\widehat{T}^2_{n,m}]_{\text{sub}}$ in (\ref{Eq_MVD_V_nm}).
This underestimation is the same for the MMD test and $V[(n+m)\widehat{\Delta}^2_{n,m}]$ is estimated by $(1 + \tau) V[(n+m)\widehat{\Delta}^2_{n,m}]_{\text{sub}}$ in (\ref{Eq_MMD_V_nm}) with positive $ \tau> 0 $.
Our approximation of the null distribution is based on a modification of the large variance of $1/\{\rho(1-\rho)\} \sum_{\ell=1}^{\infty} \widehat{\lambda}_{\ell}^{(n)} Z_{\ell}^2$ to  $(1+\tau)V[(n+m)\widehat{T}^2_{n,m}]_{\text{sub}}$.
The method aims to approximate the exact null distribution by using 
\begin{equation}\label{Eq_W'}
W_n'= \xi_n/\{\rho(1-\rho)\}\sum_{\ell=1}^{n-1} \widehat{\lambda}_{\ell}^{(n)} Z_{\ell}^2+c_n.
\end{equation}
The parameters $\xi_{n}$ and $c_{n}$ are determined so that the means of $W'_n$ and $1/\{\rho(1-\rho)\}\sum_{\ell=1}^{n-1} \widehat{\lambda}_{\ell}^{(n)}$ are equal and the variance of $W'_n$ is equal to $(1+\tau) V[(n+m)\widehat{T}^2_{n,m}]_{\text{sub}}$, which can be established by
\[
\mathbb{E}[W'_n]=\frac{1}{\rho(1-\rho)} \sum_{\ell=1}^{\infty} \widehat{\lambda}_{\ell}
\]
and
\[
V[W'_n]= (1+\tau) V[(n+m)\widehat{T}^2_{n,m}]_{\text{sub}}.
\]
This approximation method can be similarly discussed for the MMD test using $(1+\tau) V[(n+m)\widehat{\Delta}^2_{n,m}]_{\text{sub}}$.
In this paper, the parameter $\tau> 0$ is determined using the value of ``slope of the line" in Tables \ref{MVD_variance} and \ref{MMD_variance}.
Figure \ref{Wdash_comp} shows that this $W'_n$ can approximate the simulated exact distribution better than $1/\{\rho(1-\rho)\} \sum_{\ell'=1}^{n-1} \widehat{\lambda}^{(n)}_{\ell} Z_{\ell}^2$.
Algorithm shows how to obtain the critical value of the MVD test using this modification.
The algorithm for the MMD test can be obtained by changing $H$ and $\widehat{T}^2$ in Algorithm to $\widetilde{K}_{X}$ and $\widehat{\Delta}^2$.

\begin{algorithm}                      
	\caption{Calculation of critical value for the MVD test.}         
	\begin{algorithmic}                  
		\REQUIRE $X_1,\dots,X_n, Y_1,\dots,Y_m \in \mathcal{H},~k:\mathcal{H} \times \mathcal{H} \to \mathbb{R}$ (kernel), $0 < \alpha < 1$ (significance label) and $(k,\ell) \in \{1,\dots,[n/2]\}, \tau >0$ (parameters).
		For example, $\tau$ is selected as shown in the values of ``slope of the line" in Tables \ref{MVD_variance} and \ref{MMD_variance}.
		\\
		1. Compute the eigenvalues $\widehat{\lambda}_{\ell}^{(n)}$ of $H$ in (\ref{H_matrix_Gram}) and obtain $1/\{\rho(1-\rho)\}\sum_{\ell=1}^{n-1} \widehat{\lambda}_{\ell}^{(n)} Z_{\ell}^2(j)$ by random element $Z_1^{(j)},\dots,Z_{n-1}^{(j)} \overset{i.i.d.}{\sim} N(0,1),~j \in \{1,\dots,J\} $.\\
		2. (a) 
		Obtain copies $(k+\ell)\widehat{T}^2_{k,\ell}(i),~i \in 
		\{1,\dots,I\}$ of $(k+\ell)\widehat{T}^2_{k,\ell}$ under $H_0$ by the subsampling method.\\		
		~~~~(b) Compute subsampling variance $V[(n+m) \widehat{T}^2_{n,m}]_{\text{sub}}$ from $(k+\ell)\widehat{T}^2_{k,\ell}(i)$.\\
		3. Compute $W_n'$ in (\ref{Eq_W'}) by $1/\{\rho(1-\rho)\}\sum_{\ell=1}^{n-1} \widehat{\lambda}_{\ell}^{(n)} Z_{\ell}^2(j)$ and $V[(n+m) \widehat{T}^2_{n,m}]_{\text{sub}}$.
		\ENSURE We obtain the critical value $t_{\alpha}(W'_n)$ as the $J (1-\alpha)$-th from the top sorted in ascending order.
	\end{algorithmic}
\end{algorithm}

\subsection{Simulations}\label{simulations}
In this section, we investigate the performance of $(n+m)\widehat{T}^2_{n,m}$ under a specific null hypothesis and specific alternative hypotheses when $\mathcal{H} = \mathbb{R}^d$ and $k(\cdot, \cdot)$ is the Gaussian kernel in (\ref{Eq_Gaussign_kernel}).
In particular, a Monte Carlo simulation is performed to observe the type-I error and the power of the MVD and MMD tests.
Two cases are implemented: a uniform distribution $Q_1$ and an exponential distribution $Q_2$, with $P=N(0,1)$, all of which have means and variances 0 and 1, respectively.  
The critical values are determined based on $W'_n$ in Section \ref{Wdash} from a normal distribution.
The type-I error of $(n+m)\widehat{T}^2_{n,m}$ can be obtained by counting the number of times $(n+m)\widehat{T}^2_{n,m}$ exceeds the critical value in 1,000 iterations under the null hypothesis.
Next, the estimated power of $(n+m)\widehat{T}^2_{n,m}$ is similarly obtained by counting the number of times $(n+m)\widehat{T}^2_{n,m}$ exceeds the critical value under each alternative distribution in 1,000 iterations.
We execute the above for $(n,m)=(200,200)$ and $(500,500)$ and $d = 5, 10$, and $20$. 
It is known that the selection of the value of $\sigma$ involved in the Gaussian kernel affects the performance. 
We utilize $\sigma$ depending on dimension $d$. 
The significance level is $\alpha = 0.05$. 
The critical values are determined on the basis of $W'_n$ in Section \ref{Wdash} from a normal distribution. 
The type-I error and estimated power can be obtained by counting how many times $(n+m)\widehat{T}^2_{n,m}$ exceeds the critical values in 1,000 iterations under $P=Q$ and $P \neq Q$.
With $n_1=n/2$ and $(k,\ell)=(n/8,n/8)$, $\tau_{\text{MVD}}$ for MVD is $\tau_{\text{MVD}}=0.30928$ and $\tau_{\text{MMD}}$ for MMD is $\tau_{\text{MMD}}=0.11643$ by ``slope of the line" in Tables \ref{MVD_variance} and \ref{MMD_variance}.
The following can be seen from Table \ref{Type I_and_power}:
\begin{itemize}
	\item Table \ref{Type I_and_power} shows that the probabilities of type-I error at $d = 5$ and 10 are near the significance level of $\alpha=0.05$ for both the MVD and MMD.
	\item The probability of type-I error at $d = 20$ exceeds the significance level of $\alpha=0.05$ for the MVD, but decreases as $(n, m)$ increases.
	\item It can be seen that the critical value by $W'_n$ of the MVD tends to be estimated as less than that point of the null distribution.
    \item In hypothesis $P=Q_1$, it can be seen that the MVD has a higher power than the MMD.
	\item It can be seen that the MVD and MMD have higher powers for hypothesis $P=Q_2$ than hypothesis $P=Q_1$ and it is difficult to distinguish between the normal distribution and the uniform distribution by the MVD and MMD for a Gaussian kernel.
	\item Note that the critical value changes depending on the distribution of the null hypothesis.
\end{itemize}
\begin{table}[H]
	\centering
	\caption{Type-I error and power of the test by $(n+m)\widehat{T}^2_{n,m}$ for each sample size and each parameter $\sigma$.
	}
	\label{Type I_and_power}
	\begin{tabular}{ccc|ccc|ccc} 
		\hline
		$\sigma$&$d$&$(n,m)$&\multicolumn{3}{c|}{MVD}&\multicolumn{3}{c}{MMD}\\
		\hline
		\multicolumn{3}{c|}{ }&Type-I error&$Q_1$ &$Q_2$&Type-I error&$Q_1$ &$Q_2$\\ 
		\hline
		$d^{-3/4}$&5&(200,200)&0.060&0.797&1&0.047&0.401&1\\ 
		$d^{-3/4}$&5&(500,500)&0.072&1&1&0.063&0.877&1\\ 
		$d^{-7/8}$&5&(200,200)&0.056&0.728&1&0.052&0.305&1\\ 
		$d^{-7/8}$&5&(500,500)&0.067&0.999&1&0.053&0.735&1\\  
		$d^{-3/4}$&10&(200,200)&0.073&0.612&1&0.086&0.342&1\\ 
		$d^{-3/4}$&10&(500,500)&0.047&0.991&1&0.040&0.630&1\\ 
		$d^{-7/8}$&10&(200,200)&0.054&0.482&1&0.086&0.235&1\\ 
		$d^{-7/8}$&10&(500,500)&0.034&0.955&1&0.044&0.363&1\\ 
		$d^{-3/4}$&20&(200,200)&0.279&0.816&1&0.082&0.239&1\\  
		$d^{-3/4}$&20&(500,500)&0.099&0.948&1&0.068&0.477&1\\ 
		$d^{-7/8}$&20&(200,200)&0.060&0.332&1&0.047&0.113&0.989\\  
		$d^{-7/8}$&20&(500,500)&0.034&0.728&1&0.069&0.240&1\\ 
		\hline
	\end{tabular} 
\end{table} 
\section{Application to real datasets} \label{Applications to real data sets}
The MVD test was applied to some real data sets.
The significance level was $\alpha=0.05$ and the critical value $t_{0.05}(H)$ was obtained through 10,000 iterations of $1/\{\rho(1-\rho)\}\sum_{\ell=1}^{n-1} \widehat{\lambda}^{(n)}_{\ell} Z_{\ell}^2$ based on the eigenvalues of the matrix $H/n$.
We also calculated the critical value $t_{0.05}(W'_n)$ of the approximate distribution $W'_n$ according to Algorithm in Section \ref{Wdash}.
The $t_{0.05}(\widetilde{K}_{X})$ calculates the critical value for the MMD test from the distribution obtained based on Theorem 1 in \cite{Gretton2009a} through 10,000 iterations.
\subsection{USPS data}
The USPS dataset consists of handwritten digits represented by a $16\times16$ grayscale matrix (\url{https://www.csie.ntu.edu.tw/~cjlin/libsvmtools/datasets/multiclass.html#usps}). 
The sizes of each sample are shown in Table \ref{USPSdata}.
\begin{table}[H]
	\centering
	\caption{Sample sizes of the USPS data.}
	\label{USPSdata}
	\begin{tabular}{ccccccccccc}
		\hline 
	 	index&0&1&2&3&4&5&6&7&8&9\\ 
		\hline 
		sample size&359&264&198&166&200&160&170&147&166&177\\ 
		\hline 
	\end{tabular} 
\end{table}
Each group was divided into two sets of sample size 70 and the MVD test was applied to each set.
Table \ref{MVD_MMD_USPS} shows the results of applying the MVD and MMD tests to this USPS dataset.
Parameters $ \sigma = d^{-3/4}, d^{-7/8} $, $ n_1 = 35 $, $k$ and $\ell = 18$ are adopted and $ \tau = 0.69348$ for the MVD and $\tau=0.21990$ for the MMD were utilized from the slope of the line in Tables \ref{MVD_variance} and \ref{MMD_variance}.
In each cell, the values of $(n + m)\widehat{T}^2_{n,m}$ and $(n+m)\widehat{\Delta}^2_{n,m}$ for each number are written, with the values of $(n+m)\widehat{\Delta}^2_{n,m}$ in parentheses.

From Table \ref{MVD_MMD_USPS}, there is a tendency that different groups will be rejected and that the same groups are not rejected by the MVD test.
For $P=$ USPS 2 and $Q=$ USPS 2, the value of $(n+m)\widehat{T}^2_{n,m}$ is 2.953, which is larger than $t_{0.05}(W'_n)=2.946$ but smaller than $t_{0.05}(H)=3.416$.
On the other hand, for the MMD test, the value of $(n+m)\widehat{\Delta}^2_{n,m}$ is 5.014, which is larger than both $t_{0.05}(W'_n)=4.488$ and $t_{0.05}(H)=4.698$.
By modifying the distribution, there is a tendency to reject the null hypothesis.

\begin{table}[H]
	\centering
	\caption{Values of $(n + m)\widehat{T}^2_{n,m}$ and $(n+m)\widehat{\Delta}^2_{n,m}$ for $\sigma=d^{-3/4}$, $n_1 = 35 $, $ k, \ell = 18 $, $\tau_{\text{MVD}} = 0.69348$, and $\tau_{\text{MMD}}=0.21990$.
	}
	\label{MVD_MMD_USPS}
	\small
	\begin{tabular}{c|cccccccccc}
		\hline
		& 0 & 1 & 2 & 3 & 4 &5 & 6 & 7 & 8 & 9 \\ 
		\hline
		\hline 
		$t_{0.05}(H)$&3.056&0.680&3.416&2.742&2.747&3.319&2.705&2.079&2.884&1.941\\ 
		\hline
		$t_{0.05}(W'_n)$&2.755&0.572&2.940&2.317&2.349&2.785&2.336&1.823&2.431&1.719\\ 
		\hline
		\multirow{2}{*}{0}& 2.241 & \multirow{2}{*}{6.328} & \multirow{2}{*}{6.803}  &\multirow{2}{*}{6.834}  & \multirow{2}{*}{7.390 }&\multirow{2}{*}{6.589}  &\multirow{2}{*}{7.117}  &\multirow{2}{*}{7.513}  &\multirow{2}{*}{6.930}  &\multirow{2}{*}{4.573}  \\ 
		&(2.740)&&&&&&&&&\\
		\multirow{2}{*}{1}& \multirow{2}{*}{(124.6)}& 0.287 &\multirow{2}{*}{4.269 } &\multirow{2}{*}{4.290}  &\multirow{2}{*}{4.356}  &\multirow{2}{*}{4.393}  &\multirow{2}{*}{5.132}  &\multirow{2}{*}{4.265}  &\multirow{2}{*}{4.233}  &\multirow{2}{*}{4.170}  \\  
		&&(0.585)&&&&&&&&\\
		\multirow{2}{*}{2}&\multirow{2}{*}{(34.38)} & \multirow{2}{*}{(94.14)}& 2.953 &\multirow{2}{*}{4.730}  &\multirow{2}{*}{5.253}  &\multirow{2}{*}{5.053}  &\multirow{2}{*}{5.880}  &\multirow{2}{*}{5.264}  &\multirow{2}{*}{4.732}  &\multirow{2}{*}{5.165}  \\ 
		&&&(5.014)&&&&&&&\\ 
		\multirow{2}{*}{3}& \multirow{2}{*}{(42.61)} & \multirow{2}{*}{(105.0)} & \multirow{2}{*}{(26.78)}& 2.383&\multirow{2}{*}{5.248}  &\multirow{2}{*}{4.239} &\multirow{2}{*}{6.022}& \multirow{2}{*}{5.242} &\multirow{2}{*}{4.575}  &\multirow{2}{*}{5.022} \\  
		&&& & (3.345)&&&&&&\\
		\multirow{2}{*}{4}&\multirow{2}{*}{(55.81)}&\multirow{2}{*}{(93.27)}  &\multirow{2}{*}{(36.46)}  & \multirow{2}{*}{(45.34)} & 2.067 &\multirow{2}{*}{5.259}  &\multirow{2}{*}{6.237}  &\multirow{2}{*}{4.930}  &\multirow{2}{*}{4.849}  &\multirow{2}{*}{3.717} \\  
		&&&&& (2.745) &&&&&\\
		\multirow{2}{*}{5}&  \multirow{2}{*}{(30.65)}  & \multirow{2}{*}{(95.83)} & \multirow{2}{*}{(24.64)} & \multirow{2}{*}{(18.91)} & \multirow{2}{*}{(35.32)} & 2.761 &\multirow{2}{*}{5.757}  &\multirow{2}{*}{5.434}  &\multirow{2}{*}{4.814} &\multirow{2}{*}{5.106} \\ 
		&&&&&& (3.822)&&&&\\ 
		\multirow{2}{*}{6}&\multirow{2}{*}{(39.15)} & \multirow{2}{*}{(102.6)} & \multirow{2}{*}{(30.11)} & \multirow{2}{*}{(47.36)} & \multirow{2}{*}{(48.80)} & \multirow{2}{*}{(29.49)} & 2.261 &\multirow{2}{*}{6.527} &\multirow{2}{*}{5.946}&\multirow{2}{*}{6.344}\\
		&&&&&&& (5.643) &&&\\  
		\multirow{2}{*}{7}&\multirow{2}{*}{(72.41)} & \multirow{2}{*}{(111.3)} & \multirow{2}{*}{(50.29)} & \multirow{2}{*}{(52.91)} & \multirow{2}{*}{(45.54)} & \multirow{2}{*}{(51.29)} & \multirow{2}{*}{(74.30)} & 1.560 &\multirow{2}{*}{5.142}&\multirow{2}{*}{4.062}\\  
		&&&&&&&&(1.785)&&\\
		\multirow{2}{*}{8}& \multirow{2}{*}{(44.46)} & \multirow{2}{*}{(86.90)} & \multirow{2}{*}{(25.20)} & \multirow{2}{*}{(28.77)} & \multirow{2}{*}{(31.13)} & \multirow{2}{*}{(25.01)} & \multirow{2}{*}{(40.92)} & \multirow{2}{*}{(51.27)} & 2.055 &\multirow{2}{*}{4.352} \\ 
		&&&&&&&&& (2.666)&\\ 
		\multirow{2}{*}{9}&\multirow{2}{*}{(71.81)}& \multirow{2}{*}{(95.38)} & \multirow{2}{*}{(51.48)} & \multirow{2}{*}{(49.58)} & \multirow{2}{*}{(25.87)} & \multirow{2}{*}{(46.76)} & \multirow{2}{*}{(70.81)} & \multirow{2}{*}{(31.19)} & \multirow{2}{*}{(33.73)} & 1.677\\  
		&&&&&&&&&& (2.336)\\
		\hline
		\hline
		$t_{0.05}(\widetilde{K}_{X})$&(4.983)&(2.191)&(4.698)&(4.352)& (4.343)&(4.691)&(4.550)&(3.917)&(4.424)&(3.767)\\ 
		\hline 
		$t_{0.05}(W'_n)$&(5.228)&(1.749)&(4.502)&(3.928)&(3.961)&(4.085)&(4.510)&(4.104)&(4.245)&(4.455)\\ 
	\end{tabular}
\end{table}

\subsection{MNIST data}
The MNIST dataset consists of images of $28 \times 28=784$ pixels in size (\url{http://yann.lecun.com/exdb/mnist}). 
The sizes of each sample are shown in Table \ref{MNISTdata}.

\begin{table}[H]
	\centering
	\caption{Sample sizes of the MNIST data.}
	\label{MNISTdata}
	\begin{tabular}{ccccccccccc}
		\hline
		index& 0 & 1 & 2 & 3 &4  & 5 & 6 & 7 & 8 & 9 \\ 
		\hline 
		sample size&5,923  & 6,742 & 5,958 & 6,131 & 5,842 & 5,421 &5,918  & 6,265 & 5,851 &5,949  \\ 
		\hline 
	\end{tabular}
\end{table} 
The MNIST data are divided into two sets of sample size 2,000 and the MVD and MMD tests are applied.
Table \ref{MVD_MMD_MNIST} shows the results of applying the MVD  and MMD tests to the MNIST data.
The approximate distribution $W'_n$ is calculated with $n_1=1,000$, $k, \ell=500$, $\tau_{\text{MVD}}=0.69348$, and $\tau_{\text{MMD}}=0.21990$.
As in Table \ref{MVD_MMD_USPS}, the values of $(n+m)\widehat{T}^2_{n,m}$ and $(n+m)\widehat{\Delta}^2_{n,m}$ are written in each cell, with the values of $(n+m)\widehat{\Delta}^2_{n,m}$ in parentheses.
In Table \ref{MVD_MMD_MNIST}, $(n+m)\widehat{T}_{n,m}^2$ tends to take a larger value than both $t_{0.05}(H)$ and $t_{0.05}(W'_n)$.
This result is the same for the MMD test.
The MVD and MMD tests tend to reject the null hypothesis with the modifications in Section \ref{Wdash}.

\begin{table}[H]
	\centering
	\caption{Values of $(n + m)\widehat{T}^2_{n,m}$ and $(n+m)\widehat{\Delta}^2_{n,m}$ for $\sigma=d^{-2}$, $n_1 = 1,000 $, $ k, \ell = 500 $, $\tau_{\text{MVD}} = 0.69348$, and $\tau_{\text{MMD}}=0.21990$.
	}
	\label{MVD_MMD_MNIST}
	\small
		\begin{tabular}{c|cccccccccc}
			\hline
			& 0 & 1 & 2 & 3 & 4 &5 & 6 & 7 & 8 & 9 \\ 
			\hline
			\hline 
			$t_{0.05}(H)$&4.194&3.883&4.202&4.196&4.189&4.195&4.182&4.163&4.197&4.171\\ 
			\hline 
			$t_{0.05}(W'_n)$&3.993&3.937&3.993&3.989&3.990&3.993&3.989&3.990&3.989&3.988\\  
			\hline
			\multirow{2}{*}{0}& 3.999 & \multirow{2}{*}{34.86} & \multirow{2}{*}{4.207} & \multirow{2}{*}{4.268} & \multirow{2}{*}{4.394} &  \multirow{2}{*}{4.304} & \multirow{2}{*}{4.599} & \multirow{2}{*}{5.240} & \multirow{2}{*}{4.256} & \multirow{2}{*}{4.861} \\ 
			&(4.092)&&&&&&&&&\\
			\multirow{2}{*}{1}& \multirow{2}{*}{(217.4)} &5.379 & \multirow{2}{*}{34.68} & \multirow{2}{*}{34.75} & \multirow{2}{*}{34.86} & \multirow{2}{*}{34.80} & \multirow{2}{*}{35.08} & \multirow{2}{*}{35.66} & \multirow{2}{*}{34.73} & \multirow{2}{*}{35.33}  \\  
			&&(15.42)&&&&&&&&\\
			\multirow{2}{*}{2}&\multirow{2}{*}{(10.77)} & \multirow{2}{*}{(210.5)}& 4.001 & \multirow{2}{*}{4.118} & \multirow{2}{*}{4.245} & \multirow{2}{*}{4.156} & \multirow{2}{*}{4.451} & \multirow{2}{*}{5.088} & \multirow{2}{*}{4.106} & \multirow{2}{*}{4.711}\\ 
			&&&(4.131)&&&&&&&\\ 
			\multirow{2}{*}{3}& \multirow{2}{*}{(13.55)} & \multirow{2}{*}{(213.4)} & \multirow{2}{*}{(9,806)}& 4.007 & \multirow{2}{*}{4.306} & \multirow{2}{*}{4.211} & \multirow{2}{*}{4.512} & \multirow{2}{*}{5.146} & \multirow{2}{*}{4.164} & \multirow{2}{*}{4.769}\\ 
			&&& & (4.208)&&&&&&\\
			\multirow{2}{*}{4}&  \multirow{2}{*}{(16.56)}&\multirow{2}{*}{(216.1)}  &\multirow{2}{*}{(12.83)}& \multirow{2}{*}{(15.59)}&4.020 & \multirow{2}{*}{4.341} & \multirow{2}{*}{4.637} & \multirow{2}{*}{5.262} & \multirow{2}{*}{4.292} & \multirow{2}{*}{4.805}\\
			&&&&& (4.482) &&&&&\\
			\multirow{2}{*}{5}& \multirow{2}{*}{(13.35)}  & \multirow{2}{*}{(213.9)} & \multirow{2}{*}{(10.10)} & \multirow{2}{*}{(11.57)} & \multirow{2}{*}{(15.12)} & 
			4.018 & \multirow{2}{*}{4.546} & \multirow{2}{*}{5.188} & \multirow{2}{*}{4.201} & \multirow{2}{*}{4.806}\\  
			&&&&&& (4.357)&&&&\\ 
			\multirow{2}{*}{6}&  \multirow{2}{*}{(19.39)} & \multirow{2}{*}{(219.5)} & \multirow{2}{*}{(16.06)} & \multirow{2}{*}{(18.90)} & \multirow{2}{*}{(21.28)} & \multirow{2}{*}{(18.29)} & 4.031 & \multirow{2}{*}{5.485} & \multirow{2}{*}{4.499} & \multirow{2}{*}{5.105} \\
			&&&&&&& (4.573) &&&\\  
			\multirow{2}{*}{7}&\multirow{2}{*}{(28.85)} & \multirow{2}{*}{(225.2)} & \multirow{2}{*}{(24.85)} & \multirow{2}{*}{(27.49)} &\multirow{2}{*}{(28.28)} & \multirow{2}{*}{(27.56)} & \multirow{2}{*}{(34.24)} &  4.067 &  \multirow{2}{*}{5.138} & \multirow{2}{*}{5.625}\\  
			&&&&&&&&(4.625)&&\\
			\multirow{2}{*}{8}& \multirow{2}{*}{(13.12)} & \multirow{2}{*}{(211.3)} & \multirow{2}{*}{(9.261)} & \multirow{2}{*}{(11.37)} & \multirow{2}{*}{(14.72)} & \multirow{2}{*}{(11.51)} & \multirow{2}{*}{(18.19)} & \multirow{2}{*}{(26.97)}& 4.005 &\multirow{2}{*}{4.755} \\ 
			&&&&&&&&& (4.190)&\\ 
			\multirow{2}{*}{9}& \multirow{2}{*}{(24.80)}& \multirow{2}{*}{(223.2)} & \multirow{2}{*}{(21.11)} & \multirow{2}{*}{(23.23)} & \multirow{2}{*}{(19.05)} & \multirow{2}{*}{(22.80)} & \multirow{2}{*}{(29.84)} & \multirow{2}{*}{(29.91)} & \multirow{2}{*}{(22.07)} & 4.046\\  
			&&&&&&&&&& (4.686)\\
			\hline
			\hline
			$t_{0.05}(\widetilde{K}_{X})$&(4.210)&(4.718)&(4.208)&(4.206)&(4.210)&(4.209)&(4.218)&(4.238)&(4.207)&(4.224)\\ 
			\hline 
			$t_{0.05}(W'_n)$& (4.058)&(5.137)&(4.024)&(4.038)&(4.077)&(4.107)&(4.092)&(4.137)&(4.039)&(4.159)\\  
	\end{tabular}
\end{table}

\subsection{Colon data}
The Colon dataset contains gene expression data from DNA microarray experiments of colon tissue samples with $d$ = 2,000 and $n = 62$ (see \cite{Alon} for details). 
Among the 62 samples, 40 are tumor tissues and 22 are normal tissues. 
Tables \ref{MVD_tumor_normal} and \ref{MMD_tumor_normal} show the results of the MVD and MMD tests for $P=$ tumor and $Q=$ normal.
The ``tumor vs. normal" column shows the values of $(n+m)\widehat{T}^2_{n,m}$ and $(n+m)\widehat{\Delta}^2_{n,m}$ for $P =$ tumor and $Q =$ normal.
The ``normal" and ``tumor" columns show, $t_ {0.05} (W')$ and $t_ {0.05} (H)$ calculated from respectively the normal tissues and tumor tissues datasets.

For the MVD, $(n+m)\widehat{T}^2_{n,m}$ does not exceed $t_{0.05}(H)$, but $(n+m)\widehat{T}^2_{n,m}$ exceeds $t_{0.05}(W'_n)$ by modifying the approximate distribution.
By contrast, for the MMD, $(n+m)\widehat{\Delta}^2_{n,m}$ exceeds both $t_{0.05}(H)$ and $t_{0.05}(W_n')$ without modifying the approximate distribution.

\begin{table}[H]
	\centering
	\caption{Values of $(n + m)\widehat{T}^2_{n,m}$ and critical values for normal ($n_1= 11,~k, \ell = 6$) and tumor ($n_1 = 20,~k, \ell = 10$), with $\tau_{\text{MVD}} = 0.69348$.}
	\label{MVD_tumor_normal}
	\begin{tabular}{cc|cc|cc}
		\hline
		&&\multicolumn{2}{c|}{normal}&\multicolumn{2}{c}{tumor}\\
		$\sigma$&tumor vs. normal & $t_{0.05}(W'_n)$&$t_{0.05}(H)$ &$t_{0.05}(W'_n)$&$t_{0.05}(H)$\\
		\hline
		$d^{-3/4}$& 3.867&3.536&5.280&3.728&5.050\\
		$d^{-7/8}$& 2.291&2.097&2.907&2.258&2.846\\
		$d^{-1}$& 0.684&0.660&0.879&0.757&0.906\\
		\hline
	\end{tabular}
\end{table}

\begin{table}[H]
	\centering
	\caption{Values of $(n + m)\widehat{\Delta}^2_{n,m}$ and critical values for normal ($n_1=11,~k, \ell=6$) and tumor ($n_1=20,~k, \ell=10$), with $\tau_{\text{MMD}}=0.21990$.}
	\label{MMD_tumor_normal}
	\begin{tabular}{cc|cc|cc}
		\hline
		&&\multicolumn{2}{c|}{normal}&\multicolumn{2}{c}{tumor}\\
		$\sigma$&tumor vs. normal & $t_{0.05}(W'_n)$&$t_{0.05}(H)$ &$t_{0.05}(W'_n)$&$t_{0.05}(H)$\\
		\hline
		$d^{-3/4}$& 6.695&4.456&6.201&4.618& 5.713\\
		$d^{-7/8}$& 8.787&3.974&4.827&3.945&4.491\\
		$d^{-1}$& 6.282&2.439&2.754&2.412&2.634\\
		\hline
	\end{tabular}
\end{table}

Next, tumor $ (\text{sample size} = 40) $ was divided into $P=$ tumor 1 $ (n = 20) $ and $Q=$ tumor 2 $ (m = 20) $, and two-sample tests by the MVD and MMD were applied.
The results are shown in Tables \ref{MVD_tumor_tumor} and \ref{MMD_tumor_tumor}, with the values for $(n+m)\widehat{T}_{n,m}^2$ and $(n+m)\widehat{\Delta}_{n,m}^2$ in the column ``tumor 1 vs. tumor 2".
In Table \ref{MVD_tumor_tumor}, when $\sigma=d^{-3/4}$, $(n+m)\widehat{T}_{n,m}^2$ exceeds $t_{0.05}(W'_n)$, but in other cases $(n+m)\widehat{T}_{n,m}^2$ does not exceed $t_{0.05}(W'_n)$ and $(n+m)\widehat{T}_{n,m}^2$ does not exceed $t_{0.05}(H)$ and $t_{0.05}(W_n')$.
Table \ref{MMD_tumor_tumor} shows that, for all $\sigma$, $(n+m)\widehat{\Delta}^2_{n,m}$ exceeds $t_{0.05}(W'_n)$ for the MMD test, but $(n+m)\widehat{\Delta}^2_{n,m}$ does not exceed $t_{0.05}(H)$.

\begin{table}[H]
	\centering
	\caption{Values of $(n + m)\widehat{T}^2_{n,m}$ and critical values for tumor 1 ($n_1=10,~k, \ell=5$) and tumor 2 ($n_1=10,~k, \ell=5$), with $\tau=0.69348$.}
	\label{MVD_tumor_tumor}
	\begin{tabular}{cc|cc|cc}
		\hline
		&&\multicolumn{2}{c|}{tumor 1}&\multicolumn{2}{c}{tumor 2}\\
		$\sigma$&tumor 1 vs. tumor 2 & $t_{0.05}(W'_n)$&$t_{0.05}(H)$ &$t_{0.05}(W'_n)$&$t_{0.05}(H)$\\
		\hline
		$d^{-3/4}$& 3.379&3.180&4.596&3.245&4.942\\
		$d^{-7/8}$& 1.858&1.915&2.502&2.085&2.875\\
		$d^{-1}$&0.558&0.629&0.800&0.727&0.921\\
		\hline
	\end{tabular}
\end{table}

\begin{table}[H]
	\centering
	\caption{Values of $(n + m)\widehat{T}^2_{n,m}$ and critical values for tumor 1 ($n_1=10,~k, \ell=5$) and tumor 2 ($n_1=10,~k, \ell=5$), with $\tau=0.21990$.
	}
	\label{MMD_tumor_tumor}
	\begin{tabular}{cc|cc|cc}
		\hline
		&&\multicolumn{2}{c|}{tumor 1}&\multicolumn{2}{c}{tumor 2}\\
		$\sigma$&tumor 1 vs. tumor 2 & $t_{0.05}(W'_n)$&$t_{0.05}(H)$ &$t_{0.05}(W'_n)$&$t_{0.05}(H)$\\
		\hline
		$d^{-3/4}$& 4.627&4.064&5.621&3.961&5.782\\
		$d^{-7/8}$& 4.123&3.305&4.206&3.426&4.551\\
		$d^{-1}$&2.453&1.942&2.377&2.102&2.656\\
		\hline
	\end{tabular}
\end{table}
\section{Conclusion}\label{Conclusion}
We defined a Maximum Variance Discrepancy (MVD) with a similar idea to the Maximum Mean Discrepnacy (MMD) in Section \ref{MVD}.
We derived the asymptotic null distribution for the MVD test in Section \ref{asymptotic null distribution}.
This was the infinite sum of the weighted chi-square distributions.
In Section \ref{asymptotic nonnull distribution}, we derived an asymptotic nonnull distribution for the MVD test, which was a normal distribution.
The asymptotic normality of the test under the alternative hypothesis showed that the two-sample test by the MVD has consistency.
Furthermore, we developed an asymptotic distribution for the test under a sequence of local alternatives in Section \ref{asymptotic distribution under contiguous alternatives}. 
This was the infinite sum of weighted noncentral chi-squared distributions.
We constructed an estimator of asymptotic null distributed weights based on the Gram matrix in Section \ref{Gram matrix spectrum}.
The approximate distribution of the null distribution by these estimated weights does not work well, so we modified it in Section \ref{approximation of the null distribution}.
In the simulation of the power reported, we found that the power of the two-sample test by the MVD was larger than that of the MMD.
We confirmed in Section \ref{Applications to real data sets} that the two-sample test by the MVD works for real data-sets.
\section{Proofs}\label{section_Proof}
\begin{lemm}\label{asymptotic equivalent}
	Suppose that Condition is satisfied.
	Then, as $n,m \to \infty$,
	\[
	\sqrt{n+m}\left(
	\norm{\Sigma_k(\widehat{P})-\Sigma_k(\widehat{Q})}_{H(k)^{\otimes 2}}^2-\norm{\widetilde{\Sigma}_k(P)-\widetilde{\Sigma}_k(Q)}_{H(k)^{\otimes 2}}^2
	\right) \xrightarrow{P}0,
	\]
	where 
	\[
	\widetilde{\Sigma}_k(P)=\frac{1}{n}\sum_{i=1}^{n}(k(\cdot,X_i)-\mu_k(P))^{\otimes 2}~~\text{and}~~\widetilde{\Sigma}_k(Q)=\frac{1}{m}\sum_{j=1}^{m}(k(\cdot,Y_j)-\mu_k(Q))^{\otimes 2}.
	\]
\end{lemm}
\begin{lemm}\label{Lemma_three_evaluate}
	Let $Y_1,\dots,Y_m \overset{i.i.d.}{\sim} Q_{nm}$.
	Suppose that Condition is satisfied.
	Then, as $n,m \to \infty$, following evaluates are obtained 
	\begin{itemize}
		\item[(i)]  \[
		\norm{\mu_k(\widehat{Q}_{nm})-\mu_k(Q_{nm})}_{H(k)} =O_p\left(\frac{1}{\sqrt{n+m}}\right),
		\]
		\item[(ii)]
		\[
		\norm{\Sigma_k(\widehat{Q}_{nm})-\Sigma_k(Q_{nm})}_{H(k)^{\otimes 2}}=O_p\left(\frac{1}{\sqrt{n+m}}\right),
		\]
		\item[(iii)]
		\[
		\norm{\widetilde{\Sigma}_k(Q_{nm})-\Sigma_k(Q_{nm})}_{H(k)^{\otimes 2}}=O_p\left(\frac{1}{\sqrt{n+m}}\right).
		\]
	\end{itemize}
\end{lemm}
\begin{lemm}\label{asymptotic_equivalent_contiguous}
	Let $X_1,\dots,X_n \overset{i.i.d.}{\sim} P$ and $Y_1,\dots,Y_m \overset{i.i.d.}{\sim} Q_{nm}$.
	Suppose that Condition is satisfied.
	Then, as $n,m \to \infty$,
	\[
	(n+m)\left(
	\norm{\Sigma_k(\widehat{P})-\Sigma_k(\widehat{Q}_{nm})}^2_{H(k)^{\otimes 2}}-\norm{\widetilde{\Sigma}_k(P)-\widetilde{\Sigma}_k(Q_{nm})}^2_{H(k)^{\otimes 2}}
	\right) \xrightarrow{P}0.
	\]
\end{lemm}
\subsection{Proof of Proposition \ref{Prop_MMD_N(0,1)_and_N(m,Sigma)}}
The kernel mean embedding $\mu_k(N(\underline{m},\Sigma))$ with the Gaussian kernel in (\ref{Eq_Gaussign_kernel}) is obtained 
\begin{equation}\label{Eq_expression_KME}
\mu_k(N(\underline{m},\Sigma))=|I_d+2\sigma \Sigma|^{-1/2}\exp\left(
-\sigma(\cdot-\underline{m})^T(I_d+2\sigma \Sigma)^{-1}(\cdot-\underline{m})
\right)
\end{equation}
and the norm of that is derived
\begin{equation}\label{Eq_norm_of_kernel_mean_embedding}
	\norm{\mu_k(N(\underline{m},\Sigma))}^2_{H(k)}=|I_d+4\sigma \Sigma|^{-1/2}
\end{equation}
by Proposition 4.2 in \cite{Kellner_and_Celisse}.
We use the property of Gaussian density 
\begin{equation}\label{Eq_property_of_Gaussian_density}
	\phi_{\Sigma_1}(\underline{x}-\underline{m}_1)\phi_{\Sigma_2}(\underline{x}-\underline{m}_2)=\phi_{\Sigma_1+\Sigma_2}(\underline{m}_1-\underline{m}_2) \phi_{(\Sigma_1^{-1}+\Sigma_2^{-1})^{-1}}(\underline{x}-\underline{m}^*),
\end{equation} 
where 
\[
\underline{m}^*=(\Sigma^{-1}_1+\Sigma^{-1}_2)^{-1}(\Sigma^{-1}_1\underline{m}_2+\Sigma^{-1}_2\underline{m}_1)
\]
and $\phi_{\Sigma}(\cdot-\underline{m})$ designates the density of $N(\underline{m},\Sigma)$, see e.g. Appendix C in \cite{Wand_and_Jones}.
The property of Gaussian density (\ref{Eq_property_of_Gaussian_density}) is used repeatedly to calculate $\mathbb{E}_{\substack{\underline{X} \sim N(\underline{\mu},\Sigma)\\ \underline{X}' \sim N(\underline{m}_0,\Sigma_0)}}[k(\underline{X},\underline{X}')]$, and we get 
\begin{align}
&\mathbb{E}_{\substack{\underline{X} \sim N(\underline{\mu},\Sigma)\\ \underline{X}' \sim N(\underline{m}_0,\Sigma_0)}}[k(\underline{X},\underline{X}')] \nonumber\\
&=\int_{\mathbb{R}^d} \int_{\mathbb{R}^d} \exp(-\sigma \norm{\underline{x}-\underline{y}}^2_{\mathbb{R^d}}) dN(\underline{\mu},\Sigma)(\underline{x}) dN(\underline{m}_0,\Sigma_0)(\underline{y})\nonumber\\
&=\left(\frac{\pi}{\sigma}\right)^{d/2}\int_{\mathbb{R}^d} \int_{\mathbb{R}^d} \phi_{\frac{1}{2\sigma}I_d}(\underline{x}-\underline{y}) dN(\underline{\mu},\Sigma)(\underline{x}) dN(\underline{m}_0,\Sigma_0)(\underline{y})\nonumber\\
&=\left(\frac{\pi}{\sigma}\right)^{d/2}\int_{\mathbb{R}^d} \int_{\mathbb{R}^d} \phi_{\frac{1}{2\sigma}I_d}(\underline{x}-\underline{y}) \phi_{\Sigma} (\underline{x}-\underline{\mu}) d\underline{x} dN(\underline{m}_0,\Sigma_0)(\underline{y})\nonumber\\
&=\left(\frac{\pi}{\sigma}\right)^{d/2}\int_{\mathbb{R}^d} \int_{\mathbb{R}^d} \phi_{\frac{1}{2\sigma}I_d+\Sigma}(\underline{y}-\underline{\mu}) \phi_{(2\sigma I_d+\Sigma^{-1})^{-1}} (\underline{x}-\underline{m}_1^*) d\underline{x} dN(\underline{m}_0,\Sigma_0)(\underline{y})\nonumber\\
&=\left(\frac{\pi}{\sigma}\right)^{d/2}\int_{\mathbb{R}^d} \phi_{\frac{1}{2\sigma}I_d+\Sigma}(\underline{y}-\underline{\mu})dN(\underline{m}_0,\Sigma_0) (\underline{y})  \int_{\mathbb{R}^d}\phi_{(2\sigma I_d+\Sigma^{-1})^{-1}} (\underline{x}-\underline{m}_1^*)  d\underline{x}\nonumber\\
&=\left(\frac{\pi}{\sigma}\right)^{d/2}\int_{\mathbb{R}^d} \phi_{\frac{1}{2\sigma}I_d+\Sigma}(\underline{y}-\underline{\mu}) \phi_{\Sigma_0}(\underline{y}-\underline{m}_0)d\underline{y}\nonumber\\
&=\left(\frac{\pi}{\sigma}\right)^{d/2}\int_{\mathbb{R}^d} \phi_{\frac{1}{2\sigma}I_d+\Sigma+\Sigma_0}(\underline{\mu}-\underline{m}_0) \phi_{((\frac{1}{2\sigma}I_d+\Sigma)^{-1}+\Sigma_0^{-1})^{-1}}(\underline{y}-\underline{m}_2^*) d\underline{y}\nonumber\\
&=\left(\frac{\pi}{\sigma}\right)^{d/2}\phi_{\frac{1}{2\sigma}I_d+\Sigma+\Sigma_0}(\underline{\mu}-\underline{m}_0) \int_{\mathbb{R}^d} \phi_{((\frac{1}{2\sigma}I_d+\Sigma)^{-1}+\Sigma_0^{-1})^{-1}}(\underline{y}-\underline{m}_2^*) d\underline{y}\nonumber\\
&=\left(\frac{\pi}{\sigma}\right)^{d/2}\phi_{\frac{1}{2\sigma}I_d+\Sigma+\Sigma_0}(\underline{\mu}-\underline{m}_0),  \label{Eq_KME_inner_product}
\end{align}
where 
\begin{align*}
&m_1^*=(2\sigma I_d+\Sigma^{-1})^{-1}(2\sigma \underline{y}+\Sigma^{-1}\underline{\mu}),\\
&m_2^*=\left\{\left(\frac{1}{2\sigma}I_d+\Sigma \right)^{-1} +\Sigma_0^{-1}\right\}^{-1} \left\{\left(\frac{1}{2\sigma}I_d+\Sigma \right)^{-1}\underline{\mu}+\Sigma_0^{-1} \underline{m}_0\right\}.
\end{align*}
Using this results (\ref{Eq_norm_of_kernel_mean_embedding}) and (\ref{Eq_KME_inner_product}), $\text{MMD}(N(\underline{0},I_d),N(\underline{m},\Sigma))^2$ is obtained as 
\begin{align*}
	&\text{MMD}(N(\underline{0},I_d),N(\underline{m},\Sigma))^2\\
	&=\norm{\mu_k(N(\underline{0},I_d))-\mu_k(N(\underline{m},\Sigma))}^2_{H(k)}\\
	&=\norm{\mu_k(N(\underline{0},I_d))}^2_{H(k)} +\norm{\mu_k(N(\underline{m},\Sigma))}^2_{H(k)} -2\left<\mu_k(N(\underline{0},I_d)),\mu_k(N(\underline{m},\Sigma))\right>_{H(k)}\\
	&=\norm{\mu_k(N(\underline{0},I_d))}^2_{H(k)} +\norm{\mu_k(N(\underline{m},\Sigma))}^2_{H(k)} -2\mathbb{E}_{\substack{\underline{X} \sim N(\underline{0},I_d)\\ \underline{X}' \sim N(\underline{m},\Sigma)}}[k(\underline{X},\underline{X}')]\\
	&=|I_d+4\sigma I_d|^{-1/2}+|I_d+4\sigma \Sigma|^{-1/2} 
	-2\left(\frac{\pi}{\sigma}\right)^{d/2} \phi_{\frac{1}{2\sigma} I_d+I_d+\Sigma}(\underline{m})\\
	&=(1+4\sigma)^{-d/2} +|I_d+4\sigma \Sigma|^{-1/2} \\
	&~~~~~
	-2\left(\frac{\pi}{\sigma}\right)^{d/2} \left|\frac{\pi}{\sigma}\left((1+2\sigma)I_d+2\sigma \Sigma\right)\right|^{-1/2} \exp\left(
	-\sigma \underline{m}^T \left((1+2\sigma) I_d+2\sigma \Sigma\right)^{-1} \underline{m}
	\right)\\
	&=|I_d+4\sigma I_d|^{-1/2}+|I_d+4\sigma \Sigma|^{-1/2} 
	-2\left(\frac{\pi}{\sigma}\right)^{d/2} \phi_{\frac{1}{2\sigma} I_d+I_d+\Sigma}(\underline{m})\\
	&=(1+4\sigma)^{-d/2} +|I_d+4\sigma \Sigma|^{-1/2} \\
	&~~~~~
	-2 \left|(1+2\sigma)I_d+2\sigma \Sigma\right|^{-1/2} \exp\left(
	-\sigma \underline{m}^T \left((1+2\sigma) I_d+2\sigma \Sigma\right)^{-1} \underline{m}
	\right).
\end{align*}
\subsection{Proof of Proposition \ref{Prop_MVD_N(0,1)_and_N(m,Sigma)}}
In the following proof, the $\text{MVD}(N(\underline{0},I_d),N(\underline{m},\Sigma))^2$ when using the Gaussian kernel in (\ref{Eq_Gaussign_kernel}) is calculated by repeatedly using (\ref{Eq_property_of_Gaussian_density}).
From the expansion of the norm
\begin{align}
	&\text{MVD}(N(\underline{0},I_d),N(\underline{m},\Sigma))^2 \nonumber\\
	&=\norm{\Sigma_k(N(\underline{0},I_d))-\Sigma_k(N(\underline{m},\Sigma))}^2_{H(k)^{\otimes 2}}\nonumber\\
	&=\norm{\Sigma_k(N(\underline{0},I_d))}^2_{H(k)^{\otimes 2}}
	+\norm{\Sigma_k(N(\underline{m},\Sigma))}^2_{H(k)^{\otimes 2}}
	-2\left<\Sigma_k(N(\underline{0},I_d)),\Sigma_k(N(\underline{m},\Sigma))\right>_{H(k)^{\otimes 2}}, \label{Eq_MVD_NandN}
\end{align}
it is sufficient for us to calculate $\left<\Sigma_k(N(\underline{\mu},\Sigma)),\Sigma_k(N(\underline{m}_0,\Sigma_0))\right>_{H(k)^{\otimes 2}}$.
The definition of $\Sigma_k(P)$ and the tensor product $h^{\otimes 2}$ lead to 
\begin{align}
	&\left<\Sigma_k(N(\underline{\mu},\Sigma)),\Sigma_k(N(\underline{m}_0,\Sigma_0))\right>_{H(k)^{\otimes 2}} \nonumber\\
	&=\left<\mathbb{E}_{\underline{X} \sim N(\underline{\mu},\Sigma)}\left[
	\left(k(\cdot,\underline{X})-\mathbb{E}_{\underline{X} \sim N(\underline{\mu},\Sigma)}[k(\cdot,\underline{X})]\right)^{\otimes 2}
	\right],\right. \nonumber\\
	&~~~~~\left. \mathbb{E}_{\underline{Y} \sim N(\underline{m}_0,\Sigma_0)}\left[
	\left(k(\cdot,\underline{Y})-\mathbb{E}_{\underline{Y} \sim N(\underline{m}_0,\Sigma_0)}[k(\cdot,\underline{Y})]\right)^{\otimes 2}
	\right]\right>_{H(k)^{\otimes 2}} \nonumber\\
	&=\mathbb{E}_{\substack{\underline{X} \sim N(\underline{\mu},\Sigma) \nonumber\\ \underline{Y} \sim N(\underline{m}_0,\Sigma_0)}}\left[
	\left<
	\left(k(\cdot,\underline{X})-\mathbb{E}_{\underline{X} \sim N(\underline{\mu},\Sigma)}[k(\cdot,\underline{X})]\right)^{\otimes 2}, \left(k(\cdot,\underline{Y})-\mathbb{E}_{\underline{Y} \sim N(\underline{m}_0,\Sigma_0)}[k(\cdot,\underline{Y})]\right)^{\otimes 2}
	\right>_{H(k)^{\otimes 2}}
	\right] \nonumber\\
	&=\mathbb{E}_{\substack{\underline{X} \sim N(\underline{\mu},\Sigma)\\ \underline{Y} \sim N(\underline{m}_0,\Sigma_0)}}\left[
	\left<
	k(\cdot,\underline{X})^{\otimes 2}-\left\{\mathbb{E}_{\underline{X} \sim N(\underline{\mu},\Sigma)}[k(\cdot,\underline{X})] \right\}^{\otimes 2}, k(\cdot,\underline{Y})^{\otimes 2}-\left\{\mathbb{E}_{\underline{Y} \sim N(\underline{m}_0,\Sigma_0)}[k(\cdot,\underline{Y})]\right\} ^{\otimes 2}
	\right>_{H(k)^{\otimes 2}}
	\right] \nonumber\\
	&=\mathbb{E}_{\substack{\underline{X} \sim N(\underline{\mu},\Sigma)\\ \underline{Y} \sim N(\underline{m}_0,\Sigma_0)}}\left[
	\left<
	k(\cdot,\underline{X})^{\otimes 2}, k(\cdot,\underline{Y})^{\otimes 2}
	\right>_{H(k)^{\otimes 2}}
	-\left<
	k(\cdot,\underline{X})^{\otimes 2}, \left\{\mathbb{E}_{\underline{Y} \sim N(\underline{m}_0,\Sigma_0)}[k(\cdot,\underline{Y})]\right\} ^{\otimes 2}
	\right>_{H(k)^{\otimes 2}} \right. \nonumber\\
	&~~~~~
	-\left<
	\left\{\mathbb{E}_{\underline{X} \sim N(\underline{\mu},\Sigma)}[k(\cdot,\underline{X})] \right\}^{\otimes 2}, k(\cdot,\underline{Y})^{\otimes 2}
	\right>_{H(k)^{\otimes 2}} \nonumber\\
	&~~~~~\left.
	+\left<
	\left\{\mathbb{E}_{\underline{X} \sim N(\underline{\mu},\Sigma)}[k(\cdot,\underline{X})] \right\}^{\otimes 2}, \left\{\mathbb{E}_{\underline{Y} \sim N(\underline{m}_0,\Sigma_0)}[k(\cdot,\underline{Y})]\right\} ^{\otimes 2}
	\right>_{H(k)^{\otimes 2}}
	\right] \nonumber\\
	&=\mathbb{E}_{\substack{\underline{X} \sim N(\underline{\mu},\Sigma)\\ \underline{Y} \sim N(\underline{m}_0,\Sigma_0)}}\left[
	\left<
	k(\cdot,\underline{X}), k(\cdot,\underline{Y})
	\right>_{H(k)}^2
	-\left<
	k(\cdot,\underline{X}), \mathbb{E}_{\underline{Y} \sim N(\underline{m}_0,\Sigma_0)}[k(\cdot,\underline{Y})]
	\right>_{H(k)}^2 \right. \nonumber\\
	&~~~~~\left.
	-\left<
	\mathbb{E}_{\underline{X} \sim N(\underline{\mu},\Sigma)}[k(\cdot,\underline{X})] , k(\cdot,\underline{Y})
	\right>_{H(k)}^2
	+\left<
	\mathbb{E}_{\underline{X} \sim N(\underline{\mu},\Sigma)}[k(\cdot,\underline{X})] , \mathbb{E}_{\underline{Y} \sim N(\underline{m}_0,\Sigma_0)}[k(\cdot,\underline{Y})]
	\right>_{H(k)}^2
	\right] \nonumber\\
	&=\mathbb{E}_{\substack{\underline{X} \sim N(\underline{\mu},\Sigma)\\ \underline{Y} \sim N(\underline{m}_0,\Sigma_0)}}\left[
	k(\underline{X},\underline{Y})^2
	-\left\{\mathbb{E}_{\underline{Y} \sim N(\underline{m}_0,\Sigma_0)}[k(\underline{X},\underline{Y})]\right\}^2
	-\left\{
	\mathbb{E}_{\underline{X} \sim N(\underline{\mu},\Sigma)}[k(\underline{X},\underline{Y})]
	\right\}^2 \right. \nonumber\\
	&~~~~~\left.
	+\left\{
	\mathbb{E}_{\substack{\underline{X} \sim N(\underline{\mu},\Sigma) \\ \underline{Y} \sim N(\underline{m}_0,\Sigma_0)}}[k(\underline{X},\underline{Y})]
	\right\}^2
	\right] \nonumber\\
	&=\mathbb{E}_{\substack{\underline{X} \sim N(\underline{\mu},\Sigma)\\ \underline{Y} \sim N(\underline{m}_0,\Sigma_0)}}[k(\underline{X},\underline{Y})^2]
	-\mathbb{E}_{\underline{X} \sim N(\underline{\mu},\Sigma)}\left[ \left\{
	\mathbb{E}_{\underline{Y} \sim N(\underline{m}_0,\Sigma_0)}[k(\underline{X},\underline{Y})] \right\}^2
	\right] \nonumber
	\\
	&~~~~~
	-\mathbb{E}_{\underline{Y} \sim N(\underline{m}_0,\Sigma_0)}\left[ \left\{
	\mathbb{E}_{\underline{X} \sim N(\underline{\mu},\Sigma)}[k(\underline{X},\underline{Y})] \right\}^2
	\right]
	+\left\{\mathbb{E}_{\substack{\underline{X} \sim N(\underline{\mu},\Sigma)\\ \underline{Y} \sim N(\underline{m}_0,\Sigma_0)}}[k(\underline{X},\underline{Y})] \right\}^2 \nonumber\\
	&=:I_1-I_2-I_3+I_4. \label{Eq_prop2_I}
\end{align}
We calculate each of these terms.
The first term $I_1$ is derived as 
\begin{align}
I_1
&=\mathbb{E}_{\substack{\underline{X} \sim N(\underline{\mu},\Sigma) \\ \underline{Y} \sim N(\underline{m}_0,\Sigma_0)}}[k(\underline{X},\underline{Y})^2] \nonumber\\
&=\int_{\mathbb{R}^d}  \int_{\mathbb{R}^d} \exp\left(-2\sigma \norm{\underline{x}-\underline{y}}^2_{\mathbb{R}^d}\right)d N(\underline{m}_0,\Sigma_0)(\underline{x}) dN(\underline{\mu},\Sigma)(\underline{y}) \nonumber\\
&=\left(\frac{\pi}{2\sigma}\right)^{d/2} \int_{\mathbb{R}^d}\int_{\mathbb{R}^d} \phi_{\frac{1}{4\sigma}I_d}(\underline{x}-\underline{y})d N(\underline{m}_0,\Sigma_0)(\underline{x})dN(\underline{\mu},\Sigma)(\underline{y}) \nonumber\\
&=\left(\frac{\pi}{2\sigma}\right)^{d/2} \int_{\mathbb{R}^d }\int_{\mathbb{R}^d} \phi_{\frac{1}{4\sigma}I_d}(\underline{x}-\underline{y}) \phi_{\Sigma_0}(\underline{x}-\underline{m}_0)d\underline{x} dN(\underline{\mu},\Sigma)(\underline{y}) \nonumber\\
&=\left(\frac{\pi}{2\sigma}\right)^{d/2} \int_{\mathbb{R}^d}  \phi_{\frac{1}{4\sigma}I_d+\Sigma_0}(\underline{y}-\underline{m}_0) dN(\underline{\mu},\Sigma)(\underline{y}) \nonumber\\
&=\left(\frac{\pi}{2\sigma}\right)^{d/2} \int_{\mathcal{\mathbb{R}}} \phi_{\frac{1}{4\sigma}I_d}(\underline{x}-\underline{y}) \phi_{\Sigma_0}(\underline{y}-\underline{m}_0)d\underline{y} \nonumber\\
&=\left(\frac{\pi}{2\sigma}\right)^{d/2} \phi_{\frac{1}{4\sigma}I_d+\Sigma_0}(\underline{x}-\underline{m}_0).  \label{Eq_prop2_I_1}
\end{align}
by repeatedly using (\ref{Eq_property_of_Gaussian_density}).
By using the expression of kernel mean embedding in (\ref{Eq_expression_KME}) and the property of (\ref{Eq_property_of_Gaussian_density}), we obtain the second term 
\begin{align}
	I_2
	&=\mathbb{E}_{\underline{X} \sim N(\underline{\mu},\Sigma)}\left[ \left\{
	\mathbb{E}_{\underline{Y} \sim N(\underline{m}_0,\Sigma_0)}[k(\underline{X},\underline{Y})] \right\}^2
	\right] \nonumber\\
	&=\int_{\mathbb{R}^d} \left\{
	|I_d+2\sigma \Sigma_0|^{-1/2} \exp\left(
	-\sigma(\underline{x}-\underline{m}_0)^T (I_d+2\sigma \Sigma_0)^{-1}(\underline{x}-\underline{m}_0)
	\right)
	\right\}^2 dN(\underline{\mu},\Sigma)(\underline{x})\nonumber\\
	&=|I_d+2\sigma \Sigma_0|^{-1} \int_{\mathbb{R}^d} \exp\left(
	-2\sigma(\underline{x}-\underline{m}_0)^T (I_d+2\sigma \Sigma_0)^{-1}(\underline{x}-\underline{m}_0)
	\right) dN(\underline{\mu},\Sigma)(\underline{x}) \nonumber\\
	&=\left(\frac{\pi}{2\sigma}\right)^{d/2}|I_d+2\sigma \Sigma_0|^{-1/2} \int_{\mathbb{R}^d} \phi_{\frac{1}{4\sigma}(I_d+2\sigma \Sigma_0)} (\underline{x}-\underline{m}_0) dN(\underline{\mu},\Sigma)(\underline{x}) \nonumber\\
	&=\left(\frac{\pi}{2\sigma}\right)^{d/2}|I_d+2\sigma \Sigma_0|^{-1/2} \int_{\mathbb{R}^d} \phi_{\frac{1}{4\sigma}(I_d+2\sigma \Sigma_0)} (\underline{x}-\underline{m}_0) \phi_{\Sigma}(\underline{x}-\underline{\mu}) d\underline{x} \nonumber\\
	&=\left(\frac{\pi}{2\sigma}\right)^{d/2}|I_d+2\sigma \Sigma_0|^{-1/2} \phi_{\frac{1}{4\sigma}(I_d+2\sigma \Sigma_0+4\sigma \Sigma)}(\underline{m}_0-\underline{\mu}). \label{Eq_prop2_I_2}
\end{align}
The third term $I_3$ is derived 
\begin{align}
I_3
&=\mathbb{E}_{\underline{Y} \sim N(\underline{m}_0,\Sigma_0)}\left[ \left\{
\mathbb{E}_{\underline{X} \sim N(\underline{\mu},\Sigma)}[k(\underline{X},\underline{Y})] \right\}^2
\right] \nonumber \\
&=\left(\frac{\pi}{2\sigma}\right)^{d/2}|I_d+2\sigma \Sigma|^{-1/2} \phi_{\frac{1}{4\sigma}(I_d+2\sigma \Sigma+4\sigma \Sigma_0)}(\underline{m}_0-\underline{\mu}) \label{Eq_prop2_I_3}
\end{align}
by the same calculation as $I_2$.
Finally, the fourth term $I_4$ is calculated as follows 
\begin{align}
	I_4
	&=\left\{\mathbb{E}_{\substack{\underline{X} \sim N(\underline{\mu},\Sigma)\\ \underline{Y} \sim N(\underline{m}_0,\Sigma_0)}}[k(\underline{X},\underline{Y})] \right\}^2 \nonumber\\
	&=\left(\frac{\pi}{\sigma}\right)^d\left\{
	\phi_{\frac{1}{2\sigma}I_d+\Sigma+\Sigma_0}(\underline{m}_0-\underline{\mu})
	\right\}^2 \nonumber\\
	&=\left(\frac{\pi}{\sigma}\right)^d\left\{ \left|2\pi \left(\frac{1}{2\sigma}I_d+\Sigma+\Sigma_0\right) \right|^{-1/2} \exp\left(-\sigma(\underline{m}_0-\underline{\mu})^T(I_d+2\sigma \Sigma+2\sigma \Sigma_0)^{-1} (\underline{m}_0-\underline{\mu})\right)\right\}^2 \nonumber\\
	&=|I_d+2\sigma \Sigma+2\sigma \Sigma_0|^{-1} \exp\left(-2\sigma(\underline{m}_0-\underline{\mu})^T(I_d+2\sigma \Sigma+2\sigma \Sigma_0)^{-1} (\underline{m}_0-\underline{\mu})\right) \nonumber\\
	&=|I_d+2\sigma \Sigma+2\sigma \Sigma_0|^{-1} \left|\frac{\pi}{2 \sigma}(I_d+2\sigma \Sigma+2\sigma \Sigma_0)\right|^{1/2} \phi_{\frac{1}{4\sigma}(I_d+2\sigma \Sigma+2\sigma \Sigma_0)}(\underline{m}_0-\underline{\mu}) \nonumber\\
	&=\left(\frac{\pi}{2\sigma}\right)^{d/2} |I_d+2\sigma \Sigma+2\sigma \Sigma_0|^{-1/2} \phi_{\frac{1}{4\sigma}(I_d+2\sigma \Sigma+2\sigma \Sigma_0)}(\underline{m}_0-\underline{\mu}) \label{Eq_prop2_I_4}
\end{align}
by using (\ref{Eq_KME_inner_product}).
Hence, combining  (\ref{Eq_prop2_I}) and (\ref{Eq_prop2_I_1})-(\ref{Eq_prop2_I_4}) yields 
\begin{align}
	&\left<\Sigma_k(N(\underline{\mu},\Sigma)),\Sigma_k(N(\underline{m}_0,\Sigma_0))\right>_{H(k)^{\otimes 2}} \nonumber\\
	&=\left(\frac{\pi}{2\sigma}\right)^{d/2} \phi_{\frac{1}{4\sigma}I_d+\Sigma_0+\Sigma}(\underline{m}_0-\underline{\mu})
	-\left(\frac{\pi}{2\sigma}\right)^{d/2}|I_d+2\sigma \Sigma_0|^{-1/2} \phi_{\frac{1}{4\sigma}(I_d+2\sigma \Sigma_0+4\sigma \Sigma)}(\underline{m}_0-\underline{\mu}) \nonumber\\
	&~~~~~
	-\left(\frac{\pi}{2\sigma}\right)^{d/2}|I_d+2\sigma \Sigma|^{-1/2} \phi_{\frac{1}{4\sigma}(I_d+2\sigma \Sigma+4\sigma \Sigma_0)}(\underline{m}_0-\underline{\mu})\nonumber\\
	&~~~~~
	+\left(\frac{\pi}{2\sigma}\right)^{d/2} |I_d+2\sigma \Sigma+2\sigma \Sigma_0|^{-1/2} \phi_{\frac{1}{4\sigma}(I_d+2\sigma \Sigma+2\sigma \Sigma_0)}(\underline{m}_0-\underline{\mu}). \label{Eq_KVE_inner_product}
\end{align}
The following results are obtained by using (\ref{Eq_KVE_inner_product}):
\begin{align}
	&\norm{\Sigma_k(N(\underline{0},I_d))}^2_{H(k)^{\otimes 2}} \nonumber\\
	&=\left(\frac{\pi}{2\sigma}\right)^{d/2}\phi_{\frac{1}{4\sigma}(1+8\sigma)I_d}(\underline{0})
	-\left(\frac{\pi}{2\sigma}\right)^{d/2}|(1+2\sigma) I_d|^{-1/2} \phi_{\frac{1}{4\sigma}(1+6\sigma)I_d}(\underline{0}) \nonumber\\
	&~~~~~-\left(\frac{\pi}{2\sigma}\right)^{d/2} |(1+2\sigma)I_d|^{-1/2} \phi_{\frac{1}{4\sigma}(1+6\sigma)I_d}(\underline{0})
	+\left(\frac{\pi}{2\sigma}\right)^{d/2} |(1+4\sigma)I_d|^{-1/2} \phi_{\frac{1}{4\sigma}(1+4\sigma)I_d}(\underline{0})\nonumber\\
	&=(1+8\sigma)^{-d/2}-2(1+8\sigma+12\sigma^2)^{-d/2}+(1+4\sigma)^{-d}, \label{Eq_prop2_term1}\\
	&\norm{\Sigma_k(N(\underline{m},\Sigma))}^2_{H(k)^{\otimes 2}} \nonumber\\
	&=\left(\frac{\pi}{2\sigma}\right)^{d/2}\phi_{\frac{1}{4\sigma}(I_d+8\sigma \Sigma)}(\underline{0})
	-\left(\frac{\pi}{2\sigma}\right)^{d/2}|I_d+2\sigma \Sigma|^{-1/2}\phi_{\frac{1}{4\sigma}(I_d+6\sigma \Sigma)}(\underline{0}) \nonumber\\
	&~~~~~-\left(\frac{\pi}{2\sigma}\right)^{d/2}|I_d+2\sigma \Sigma|^{-1/2}\phi_{\frac{1}{4\sigma}(I_d+6\sigma \Sigma)}(\underline{0})
	+\left(\frac{\pi}{2\sigma}\right)^{d/2}|I_d+4\sigma \Sigma|^{-1/2} \phi_{\frac{1}{4\sigma}(I_d+4\sigma \Sigma)}(\underline{0}) \nonumber\\
	&=|I_d+8\sigma \Sigma|^{-1/2}
	-2|I_d+8\sigma \Sigma+12\sigma^2 \Sigma^2|^{-1/2}
	+|I_d+4\sigma \Sigma|^{-1}, \label{Eq_prop2_term2}
\end{align}
and 
\begin{align}
	&\left<\Sigma_k(N(\underline{0},I_d)),\Sigma_k(N(\underline{m},\Sigma))\right>_{H(k)^{\otimes 2}} \nonumber\\
	&=|(1+4\sigma)I_d+4\sigma \Sigma|^{-1/2} \exp\left(
	-2\sigma \underline{m}^T\left(
	(1+4\sigma)I_d+4\sigma \Sigma
	\right)^{-1} \underline{m}
	\right) \nonumber\\
	&~~~~~-|I_d+2\sigma \Sigma|^{-1/2} |(1+4\sigma)I_d+2\sigma \Sigma|^{-1/2}
	\exp\left(
	-2\sigma \underline{m}^T \left(
	(1+4\sigma)I_d+2\sigma \Sigma
	\right)^{-1} \underline{m}
	\right)\nonumber\\
	&~~~~~
	-(1+2\sigma)^{-d/2} |(1+2\sigma)I_d+4\sigma \Sigma|^{-1/2}
	\exp\left(-2\sigma \underline{m}^T \left((1+2\sigma)I_d+4\sigma \Sigma\right)^{-1} \underline{m}\right) \nonumber\\
	&~~~~~
	+|(1+2\sigma)I_d+2\sigma \Sigma|^{-1} \exp\left(
	-2\sigma \underline{m}^T \left(
	(1+2\sigma)I_d+2\sigma \Sigma
	\right)^{-1} \underline{m}
	\right).  \label{Eq_prop2_term3}
\end{align}
Therefore, we give 
\begin{align*}
	&\text{MVD}(N(\underline{0},I_d),N(\underline{m},\Sigma))^2\\
	&=(1+8\sigma)^{-d/2}-2(1+8\sigma+12\sigma^2)^{-d/2}+(1+4\sigma)^{-d}\\
	&~~~~~+|I_d+8\sigma \Sigma|^{-1/2}
	-2|I_d+8\sigma \Sigma+12\sigma^2 \Sigma^2|^{-1/2}
	+|I_d+4\sigma \Sigma|^{-1}\\
	&~~~~~-2|(1+4\sigma)I_d+4\sigma \Sigma|^{-1/2} \exp\left(
	-2\sigma \underline{m}^T\left(
	(1+4\sigma)I_d+4\sigma \Sigma
	\right)^{-1} \underline{m}
	\right) \\
	&~~~~~+2|I_d+2\sigma \Sigma|^{-1/2} |(1+4\sigma)I_d+2\sigma \Sigma|^{-1/2}
	\exp\left(
	-2\sigma \underline{m}^T \left(
	(1+4\sigma)I_d+2\sigma \Sigma
	\right)^{-1} \underline{m}
	\right)\\
	&~~~~~
	+2(1+2\sigma)^{-d/2} |(1+2\sigma)I_d+4\sigma \Sigma|^{-1/2}
	\exp\left(-2\sigma \underline{m}^T \left((1+2\sigma)I_d+4\sigma \Sigma\right)^{-1} \underline{m}\right)\\
	&~~~~~
	-2|(1+2\sigma)I_d+2\sigma \Sigma|^{-1} \exp\left(
	-2\sigma \underline{m}^T \left(
	(1+2\sigma)I_d+2\sigma \Sigma
	\right)^{-1} \underline{m}
	\right) 
\end{align*}
from substituting the formulas (\ref{Eq_prop2_term1})-(\ref{Eq_prop2_term3}) to (\ref{Eq_MVD_NandN}).
\subsection{Proof of Theorem \ref{Asymptotic_Null_Distribution}}
Theorem \ref{Asymptotic_Null_Distribution} is shown by regarding $k(\cdot,X_1),\dots,k(\cdot,X_n)$ and $k(\cdot,Y_1),\dots,k(\cdot,Y_m)$ as the data in Corollary 1 of \cite{Boente2018}.
\subsection{Proof of Theorem \ref{Asymptotic_Nonnull_Distribution}}
By Lemma \ref{asymptotic equivalent}, it suffices to derive the asymptotic distribution of 
\[
\sqrt{n+m} \left\{
\norm{\widetilde{\Sigma}_k(P)-\widetilde{\Sigma}_k(Q)}^2_{H(k)^{\otimes 2}} -\norm{\Sigma_k(P)-\Sigma_k(Q)}^2_{H(k)^{\otimes 2}}
\right\}.
\]
Let us expand the following quantity
\begin{align*}
&\sqrt{n+m} \left\{
\norm{\widetilde{\Sigma}_k(P)-\widetilde{\Sigma}_k(Q)}^2_{H(k)^{\otimes 2}} -\norm{\Sigma_k(P)-\Sigma_k(Q)}^2_{H(k)^{\otimes 2}}
\right\}\\
&=\sqrt{n+m} \left< \widetilde{\Sigma}_k(P) -\widetilde{\Sigma}_k(Q) +\Sigma_k(P)-\Sigma_k(Q) , \widetilde{\Sigma}_k(P) - \widetilde{\Sigma}_k(Q) -\Sigma_k(P)+\Sigma_k(Q)\right>_{H(k)^{\otimes 2}}\\
&=\sqrt{n+m} \left<\widetilde{\Sigma}_k(P)-\Sigma_k(P)-\{\widetilde{\Sigma}_k(Q) -\Sigma_k(Q)\}, \right.\\
&~~~~~\left.2\{\Sigma_k(P)-\Sigma_k(Q)\}+\widetilde{\Sigma}_k(P)-\Sigma_k(P)-\{\widetilde{\Sigma}_k(Q)-\Sigma_k(Q)\}\right>_{H(k)^{\otimes 2}}\\
&=2\sqrt{n+m} \left<\Sigma_k(P) -\Sigma_k(Q), \widetilde{\Sigma}_k(P)-\Sigma_k(P)-\{\widetilde{\Sigma}_k(Q) -\Sigma_k(Q)\} \right>_{H(k)^{\otimes 2}}\\
&~~~~~+\sqrt{n+m} \norm{\widetilde{\Sigma}_k(P)-\Sigma_k(P)-\{\widetilde{\Sigma}_k(Q) -\Sigma_k(Q)\}}^2_{H(k)^{\otimes 2}}\\
&=\sqrt{\frac{n+m}{n}}\frac{2}{\sqrt{n}} \sum_{i=1}^{n}\left< \Sigma_k(P)-\Sigma_k(Q), (k(\cdot,X_i)-\mu_k(P))^{\otimes 2}-\Sigma_k(P)
\right>_{H(k)^{\otimes 2}} \\
&~~~~~
-\sqrt{\frac{n+m}{m}}\frac{2}{\sqrt{m}} \sum_{j=1}^{m}\left< \Sigma_k(P)-\Sigma_k(Q), (k(\cdot,Y_j)-\mu_k(Q))^{\otimes 2}-\Sigma_k(Q)
\right>_{H(k)^{\otimes 2}}+O_p\left(\frac{1}{\sqrt{n+m}}\right),
\end{align*}
which converges in distribution to$N(0,4\rho^{-1} v_P^2+4(1-\rho)^{-1} v^2_{Q})$ by the central limit theorem.
\subsection{Proof of Lemma \ref{asymptotic equivalent}}
A direct calculation gives
\begin{align*}
	&\sqrt{n+m}\left(\norm{\Sigma_k(\widehat{P})-\Sigma_k(\widehat{Q})}^2_{H(k)^{\otimes 2}}
	-\norm{\widetilde{\Sigma}_k(P)-\widetilde{\Sigma}_k(Q)}^2_{H(k)^{\otimes 2}} \right)\\
	&=\left(
	\norm{\Sigma_k(\widehat{P})-\Sigma_k(\widehat{Q})}_{H(k)^{\otimes 2}}
	+\norm{\widetilde{\Sigma}_k(P)-\widetilde{\Sigma}_k(Q)}_{H(k)^{\otimes 2}}
	\right)\\
	&~~~~~\times \sqrt{n+m}\left(
	\norm{\Sigma_k(\widehat{P})-\Sigma_k(\widehat{Q})}_{H(k)^{\otimes 2}}
	-\norm{\widetilde{\Sigma}_k(P)-\widetilde{\Sigma}_k(Q)}_{H(k)^{\otimes 2}}
	\right)\\
	&\leq \left(
	\norm{\Sigma_k(\widehat{P})-\Sigma_k(\widehat{Q})}_{H(k)^{\otimes 2}}
	+\norm{\widetilde{\Sigma}_k(P)-\widetilde{\Sigma}_k(Q)}_{H(k)^{\otimes 2}}
	\right)\\
	&~~~~~\times \sqrt{n+m}\norm{\Sigma_k(\widehat{P})-\widetilde{\Sigma}_k(P)-\left(
		\Sigma_k(\widehat{Q})-\widetilde{\Sigma}_k(Q)
		\right)}_{H(k)^{\otimes 2}}.
\end{align*}
From direct expansion $\Sigma_k(\widehat{P})=\widetilde{\Sigma}_k(P) -(\mu_k(P)-\mu_k(\widehat{P}))^{\otimes 2}$, we have
\begin{align*}
	&\sqrt{n+m}\left(\norm{\Sigma_k(\widehat{P})-\Sigma_k(\widehat{Q})}^2_{H(k)^{\otimes 2}}
	-\norm{\widetilde{\Sigma}_k(P)-\widetilde{\Sigma}_k(Q)}^2_{H(k)^{\otimes 2}} \right)\\
	&\leq \left(
	\norm{\Sigma_k(\widehat{P})-\Sigma_k(\widehat{Q})}_{H(k)^{\otimes 2}}
	+\norm{\widetilde{\Sigma}_k(P)-\widetilde{\Sigma}_k(Q)}_{H(k)^{\otimes 2}}
	\right)\\
	&~~~~~\times \sqrt{n+m}\norm{
		(\mu_k(P)-\mu_k(\widehat{P}))^{\otimes 2}-(\mu_k(Q)-\mu_k(\widehat{Q}))^{\otimes 2}
	}_{H(k)^{\otimes 2}}\\
	&\leq \left(
	\norm{\Sigma_k(\widehat{P})-\Sigma_k(\widehat{Q})}_{H(k)^{\otimes 2}}
	+\norm{\widetilde{\Sigma}_k(P)-\widetilde{\Sigma}_k(Q)}_{H(k)^{\otimes 2}}
	\right)\\
	&~~~~~\times \sqrt{n+m}\left(\norm{
		\mu_k(P)-\mu_k(\widehat{P})}^2_{H(k)}+\norm{\mu_k(Q)-\mu_k(\widehat{Q}) 
	}^2_{H(k)} \right)\\
	&=\left(
	\norm{\Sigma_k(\widehat{P})-\Sigma_k(\widehat{Q})}_{H(k)^{\otimes 2}}
	+\norm{\widetilde{\Sigma}_k(P)-\widetilde{\Sigma}_k(Q)}_{H(k)^{\otimes 2}}
	\right)\\
	&~~~~~\times \left( \frac{\sqrt{n+m}}{n}\cdot n\norm{
		\mu_k(P)-\mu_k(\widehat{P})}^2_{H(k)}+\frac{\sqrt{n+m}}{m}\cdot m\norm{\mu_k(Q)-\mu_k(\widehat{Q}) 
	}^2_{H(k)} \right)\\
	&\xrightarrow{P}0,
\end{align*}
as $n,m \to \infty$.
\subsection{Proof of Theorem \ref{Asymptotic distribution under contiguous alternatives}}
From Lemma \ref{asymptotic_equivalent_contiguous} it is sufficient for us to derive the asymptotic distribution of $(n+m)\norm{\widetilde{\Sigma}_k(P)-\widetilde{\Sigma}_k(Q_{nm})}^2_{H(k)^{\otimes 2}}$.
It follows from direct calculations that
\begin{align}
&(n+m)\norm{\widetilde{\Sigma}_k(P)-\widetilde{\Sigma}_k(Q_{nm})}^2_{H(k)^{\otimes 2}} \nonumber\\
&=(n+m)\norm{\frac{1}{n}\sum_{i=1}^{n}(k(\cdot,X_i)-\mu_k(P))^{\otimes 2}
-\frac{1}{m}\sum_{j=1}^{m}(k(\cdot,Y_j)-\mu_k(Q_{nm}))^{\otimes 2}}^2_{H(k)^{\otimes 2}} \nonumber\\
&=(n+m)\left\|\frac{1}{n}\sum_{i=1}^{n}(k(\cdot,X_i)-\mu_k(P))^{\otimes 2}-\Sigma_k(P) \right. \nonumber\\
&~~~~~\left.-\frac{1}{m}\sum_{j=1}^{m}(k(\cdot,Y_j)-\mu_k(P)+\mu_k(P)-\mu_k(Q_{nm}))^{\otimes 2}+\Sigma_k(P) \right\|^2_{H(k)^{\otimes 2}} \nonumber\\
&=(n+m)\Bigg\|\frac{1}{n}\sum_{i=1}^{n}(k(\cdot,X_i)-\mu_k(P))^{\otimes 2}-\Sigma_k(P)-\left\{
\frac{1}{m}\sum_{j=1}^{m}(k(\cdot,Y_j)-\mu_k(P))^{\otimes 2}-\Sigma_k(P)
\right\} \nonumber\\
&~~~~~
-\frac{1}{m}\sum_{j=1}^{m}(k(\cdot,Y_j)-\mu_k(P))\otimes (\mu_k(P)-\mu_k(Q_{nm}))
-\frac{1}{m}\sum_{j=1}^{m} (\mu_k(P)-\mu_k(Q_{nm}))\otimes (k(\cdot,Y_j)-\mu_k(P)) \nonumber\\
&~~~~~-(\mu_k(P)-\mu_k(Q_{nm}))^{\otimes 2} \Bigg\|^2_{H(k)^{\otimes 2}} \nonumber\\
&=(n+m)\Bigg\|\frac{1}{n}\sum_{i=1}^{n}(k(\cdot,X_i)-\mu_k(P))^{\otimes 2}-\Sigma_k(P)-\left\{
\frac{1}{m}\sum_{j=1}^{m}(k(\cdot,Y_j)-\mu_k(P))^{\otimes 2}-\Sigma_k(P)
\right\}+\frac{A}{\sqrt{n+m}} \Bigg\|^2_{H(k)^{\otimes 2}} \nonumber\\
&=(n+m)\Bigg\|\frac{1}{n}\sum_{i=1}^{n}(k(\cdot,X_i)-\mu_k(P))^{\otimes 2}-\Sigma_k(P)-\left\{
\frac{1}{m}\sum_{j=1}^{m}(k(\cdot,Y_j)-\mu_k(P))^{\otimes 2}-\Sigma_k(P)
\right\}\Bigg\|^2_{H(k)^{\otimes 2}} \nonumber\\
&~~~~~+2\sqrt{n+m}\left<\frac{1}{n}\sum_{i=1}^{n}(k(\cdot,X_i)-\mu_k(P))^{\otimes 2}-\Sigma_k(P) -\left\{
\frac{1}{m}\sum_{j=1}^{m}(k(\cdot,Y_j)-\mu_k(P))^{\otimes 2}-\Sigma_k(P)
\right\},A \right>_{H(k)^{\otimes 2}} \nonumber\\
&~~~~~+\norm{A}^2_{H(k)^{\otimes 2}}, \label{Eq_Theorem3_h}
\end{align}
where 
\[
A=(\mu_k(\widehat{Q}_{nm})-\mu_k(P))\otimes(\mu_k(P)-\mu_k(Q))+(\mu_k(P)-\mu_k(Q))\otimes (\mu_k(\widehat{Q}_{nm})-\mu_k(P))-\frac{1}{\sqrt{n+m}}(\mu_k(P)-\mu_k(Q))^{\otimes 2}.
\]
In addition, we see that
\begin{align}
\norm{\mu_k(\widehat{Q}_{nm})-\mu_k(P)}_{H(k)}
&=\norm{\mu_k(\widehat{Q}_{nm})-\mu_k(Q_{nm})+\mu_k(Q_{nm})-\mu_k(P)}_{H(k)} \nonumber\\
&\leq \norm{\mu_k(\widehat{Q}_{nm})-\mu_k(Q_{nm})}_{H(k)} 
+\norm{\mu_k(Q_{nm})-\mu_k(P)}_{H(k)} \nonumber\\
&=O_p\left(\frac{1}{\sqrt{n+m}}\right) \label{Eq_mu_Qnm_P_order}
\end{align}
by (i) in Lemma \ref{Lemma_three_evaluate}. 
Thus, (\ref{Eq_mu_Qnm_P_order}) leads that
\begin{align}
\norm{A}_{H(k)^{\otimes 2}} 
&\leq 2\norm{\mu_k(\widehat{Q}_{nm})-\mu_k(P)}_{H(k)}  \norm{\mu_k(P)-\mu_k(Q)}_{H(k) }  
+\frac{1}{\sqrt{n+m}} \norm{\mu_k(P)-\mu_k(Q)}^2_{H(k)} \nonumber\\
&=O_p\left(\frac{1}{\sqrt{n+m}}\right). \label{Eq_A_order}
\end{align}
Further, following result is obtained by (i) and (ii) in Lemma \ref{Lemma_three_evaluate} 
\begin{align}
&\norm{\frac{1}{m} \sum_{j=1}^{m}(k(\cdot,Y_j)-\mu_k(P))^{\otimes 2}-\Sigma_k(P)}_{H(k)^{\otimes 2}} \nonumber\\
&=\norm{\frac{1}{m} \sum_{j=1}^{m}(k(\cdot,Y_j)-\mu_k(P))^{\otimes 2}-\Sigma_k(Q_{nm})+\Sigma_k(Q_{nm})-\Sigma_k(P)}_{H(k)^{\otimes 2}} \nonumber\\
&\leq \norm{\frac{1}{m} \sum_{j=1}^{m}(k(\cdot,Y_j)-\mu_k(P))^{\otimes 2}-\Sigma_k(Q_{nm})}_{H(k)^{\otimes 2}} +\norm{\Sigma_k(Q_{nm})-\Sigma_k(P)}_{H(k)^{\otimes 2}}\nonumber\\
&=\norm{\frac{1}{m} \sum_{j=1}^{m}(k(\cdot,Y_j)-\mu_k(P))^{\otimes 2}-\Sigma_k(Q_{nm})}_{H(k)^{\otimes 2}} +O\left(\frac{1}{\sqrt{n+m}}\right) \nonumber\\
&=\norm{\frac{1}{m} \sum_{j=1}^{m}(k(\cdot,Y_j)-\mu_k(Q_{nm})+\mu_k(Q_{nm})-\mu_k(P))^{\otimes 2}-\Sigma_k(Q_{nm})}_{H(k)^{\otimes 2}} +O\left(\frac{1}{\sqrt{n+m}}\right) \nonumber\\
&=\Bigg\|\frac{1}{m} \sum_{j=1}^{m}(k(\cdot,Y_j)-\mu_k(Q_{nm}))^{\otimes 2}-\Sigma_k(Q_{nm})
+(\mu_k(\widehat{Q}_{nm})-\mu_k(Q_{nm})) \otimes (\mu_k(Q_{nm})-\mu_k(P)) \nonumber\\
&~~~~~
+(\mu_k(Q_{nm})-\mu_k(P)) \otimes  (\mu_k(\widehat{Q}_{nm})-\mu_k(Q_{nm}))
+(\mu_k(Q_{nm})-\mu_k(P))^{\otimes 2} \Bigg\|_{H(k)^{\otimes 2}}  +O\left(\frac{1}{\sqrt{n+m}}\right) \nonumber\\
&\leq \norm{\widetilde{\Sigma}_k(Q_{nm})-\Sigma_k(Q_{nm})}_{H(k)} 
+\frac{2}{\sqrt{n+m}}\norm{\mu_k(\widehat{Q}_{nm})-\mu_k(Q_{nm})}_{H(k)}  \norm{\mu_k(Q)-\mu_k(P)}_{H(k)^{\otimes 2}} \nonumber\\
&~~~~~
+\frac{1}{n+m}\norm{\mu_k(Q)-\mu_k(P)}^2_{H(k)}+O\left(\frac{1}{\sqrt{n+m}}\right) \nonumber\\
&=O_p\left(\frac{1}{\sqrt{n+m}}\right). \label{Eq_Theorem3_order}
\end{align}
These results (\ref{Eq_A_order}) and (\ref{Eq_Theorem3_order}) yield that
\begin{align}
&\sqrt{n+m}\left<\frac{1}{n}\sum_{i=1}^{n}(k(\cdot,X_i)-\mu_k(P))^{\otimes 2}-\Sigma_k(P) -\left\{
\frac{1}{m}\sum_{j=1}^{m}(k(\cdot,Y_j)-\mu_k(P))^{\otimes 2}-\Sigma_k(P)
\right\},A \right>_{H(k)^{\otimes 2}}\nonumber \\
&\leq \sqrt{n+m}  \norm{\frac{1}{n}\sum_{i=1}^{n}(k(\cdot,X_i)-\mu_k(P))^{\otimes 2}-\Sigma_k(P) -\left\{
\frac{1}{m}\sum_{j=1}^{m}(k(\cdot,Y_j)-\mu_k(P))^{\otimes 2}-\Sigma_k(P)
\right\}}_{H(k)^{\otimes 2}} \nonumber\\
&~~~~~\times  \norm{A}_{H(k)^{\otimes 2}} \nonumber\\
&= \sqrt{n+m} \Bigg\{
\norm{\frac{1}{n}\sum_{i=1}^{n}(k(\cdot,X_i)-\mu_k(P))^{\otimes 2}-\Sigma_k(P)}_{H(k)^{\otimes 2}} \nonumber\\
&~~~~~
+\norm{\frac{1}{m}\sum_{j=1}^{m}(k(\cdot,Y_j)-\mu_k(P))^{\otimes 2}-\Sigma_k(P)}_{H(k)^{\otimes 2}}
\Bigg\}  \norm{A}_{H(k)^{\otimes 2}} \nonumber\\
&=\sqrt{n+m} \left\{
O_p\left(\frac{1}{\sqrt{n+m}}\right) +O_p\left(\frac{1}{\sqrt{n+m}}\right) 
\right\} \cdot O_p\left(\frac{1}{\sqrt{n+m}}\right) \nonumber\\
&=O_p\left(\frac{1}{\sqrt{n+m}}\right). \label{Eq_Theorem3_innner_product_order}
\end{align}

By using (\ref{Eq_A_order}) and (\ref{Eq_Theorem3_innner_product_order}) to (\ref{Eq_Theorem3_h}), we obtain
\begin{align*}
&(n+m)\norm{\widetilde{\Sigma}_k(P)-\widetilde{\Sigma}_k(Q_{nm})}^2_{H(k)^{\otimes 2}} \nonumber\\
&=(n+m)\Bigg\|\frac{1}{n}\sum_{i=1}^{n}(k(\cdot,X_i)-\mu_k(P))^{\otimes 2}-\Sigma_k(P)-\left\{
\frac{1}{m}\sum_{j=1}^{m}(k(\cdot,Y_j)-\mu_k(P))^{\otimes 2}-\Sigma_k(P)
\right\}\Bigg\|^2_{H(k)^{\otimes 2}}\\
&~~~~~+O_p\left(\frac{1}{\sqrt{n+m}}\right)\\
&=(n+m) \left<\frac{1}{n}\sum_{i=1}^{n}(k(\cdot,X_i)-\mu_k(P))^{\otimes 2}-\Sigma_k(P)-\left\{
\frac{1}{m}\sum_{j=1}^{m}(k(\cdot,Y_j)-\mu_k(P))^{\otimes 2}-\Sigma_k(P)
\right\},\right.\\
&~~~~~\left.\frac{1}{n}\sum_{s=1}^{n}(k(\cdot,X_s)-\mu_k(P))^{\otimes 2}-\Sigma_k(P)-\left\{
\frac{1}{m}\sum_{t=1}^{m}(k(\cdot,Y_t)-\mu_k(P))^{\otimes 2}-\Sigma_k(P)
\right\}\right>_{H(k)^{\otimes 2}}\\
&~~~~~+O_p\left(\frac{1}{\sqrt{n+m}}\right)\\
&=\frac{n+m}{n^2m^2} \sum_{i,s=1}^{n} \sum_{j,t=1}^{m}\left<(k(\cdot,X_i)-\mu_k(P))^{\otimes 2}-\Sigma_k(P)-\left\{
(k(\cdot,Y_j)-\mu_k(P))^{\otimes 2}-\Sigma_k(P)
\right\},\right.\\
&~~~~~\left.(k(\cdot,X_s)-\mu_k(P))^{\otimes 2}-\Sigma_k(P)-\left\{
(k(\cdot,Y_t)-\mu_k(P))^{\otimes 2}-\Sigma_k(P)
\right\}\right>_{H(k)^{\otimes 2}}\\
&~~~~~+O_p\left(\frac{1}{\sqrt{n+m}}\right)\\
&=\frac{n+m}{n^2} \sum_{i,s=1}^{n}\left<(k(\cdot,X_i)-\mu_k(P))^{\otimes 2}-\Sigma_k(P),(k(\cdot,X_s)-\mu_k(P))^{\otimes 2}-\Sigma_k(P)\right>_{H(k)^{\otimes 2}}\\
&~~~~~+\frac{n+m}{m^2}\sum_{j,t=1}^{m} \left<(k(\cdot,Y_j)-\mu_k(P))^{\otimes 2}-\Sigma_k(P),(k(\cdot,Y_t)-\mu_k(P))^{\otimes 2}-\Sigma_k(P)\right>_{H(k)^{\otimes 2}}\\
&~~~~~-\frac{2(n+m)}{nm} \sum_{i=1}^{n} \sum_{j=1}^{m} \left<(k(\cdot,X_i)-\mu_k(P))^{\otimes 2}-\Sigma_k(P),(k(\cdot,Y_j)-\mu_k(P))^{\otimes 2}-\Sigma_k(P)\right>_{H(k)^{\otimes 2}}\\
&~~~~~+O_p\left(\frac{1}{\sqrt{n+m}}\right)\\
&=\frac{n+m}{n^2} \sum_{i,s=1}^{n} h(X_i,X_s)
+\frac{n+m}{m^2} \sum_{j,t=1}^{m} h(Y_j,Y_t)
-\frac{2(n+m)}{nm} \sum_{i=1}^{n} \sum_{j=1}^{m} h(X_i,Y_j)+O_p\left(\frac{1}{\sqrt{n+m}}\right),
\end{align*}
where $h(x,y)$ is in (\ref{Eq_h(x,y)}).
Since $\mathbb{E}_{X \sim P} [k(X,X)^2] < \infty$, $S_k$ in (\ref{Eq_S_k}) is a Hilbert--Schmidt operator by Theorem VI.22 in \cite{reed1981functional}.
Therefore, from Theorem 1 in \cite{Minh2006},  we have
\[
h(x,y)=\sum_{\ell=1}^{\infty} \theta_{\ell}\Psi_{\ell}(x) \Psi_{\ell}(y),
\]
where $\theta_{\ell}$ is eigenvalue of $S_k$ and $\Psi_{\ell}$ is eigenfunction corresponding to $\theta_{\ell}$, each satisfies
\[
\int_{\mathcal{H}}h(x,y) \Psi_{\ell}(y) dP(y) =\theta_{\ell} \Psi_{\ell} (x) ~~\text{and}~~\int_{\mathcal{H}} \Psi_i(y) \Psi_{j}(y)dP(y)=\delta_{ij}
\]
with $\delta_{ij}$ Kronecker's delta.
From $\int_{\mathcal{H}}h(x,y) \Psi_{\ell}(y) dP(y) =\theta_{\ell} \Psi_{\ell} (x)$, we  have
\[
\mathbb{E}_{X \sim P}[\Phi_{\ell}(X)]
=\frac{1}{\theta_{\ell}} \int_{\mathcal{H}} \mathbb{E}_{X \sim P}[h(X,y)] \Psi_{\ell}(y) dP(y)=0
\]
and
\begin{align*}
\mathbb{E}_{Y \sim Q_{nm}}[\Psi_{\ell}(Y)]
&=\frac{1}{\theta_{\ell}} \int_{\mathcal{H}} \mathbb{E}_{Y \sim Q_{nm}}[h(Y,y)] \Psi_{\ell}(y) dP(y).
\end{align*}
Since direct calculation, we get
\begin{align*}
&\mathbb{E}_{Y \sim Q_{nm}}[(k(\cdot,Y)-\mu_k(P))^{\otimes 2}-\Sigma_k(P)]\\
&=\mathbb{E}_{Y \sim Q_{nm}}[(k(\cdot,Y)-\mu_k(P))^{\otimes 2}]-\Sigma_k(P)\\
&=\int_{\mathcal{H}} (k(\cdot,y)-\mu_k(P))^{\otimes 2}dQ_{nm}(y) -\Sigma_k(P)\\
&=\int_{\mathcal{H}} (k(\cdot,y)-\mu_k(P))^{\otimes 2}d\left\{
\left(1-\frac{1}{\sqrt{n+m}}\right)P+\frac{1}{\sqrt{n+m}}Q
\right\}(y) -\Sigma_k(P)\\
&=\left(
1-\frac{1}{\sqrt{n+m}} 
\right)\mathbb{E}_{Y \sim P} [(k(\cdot,Y)-\mu_k(P))^{\otimes 2}]
+\frac{1}{\sqrt{n+m}} \mathbb{E}_{Y \sim Q}[(k(\cdot,Y)-\mu_k(P))^{\otimes 2}]
-\Sigma_k(P)\\
&=\left(
1-\frac{1}{\sqrt{n+m}}
\right) \Sigma_k(P)+\frac{1}{\sqrt{n+m}} \left(\Sigma_k(Q) +(\mu_k(Q)-\mu_k(P))^{\otimes 2} \right)-\Sigma_k(P)\\
&=\frac{1}{\sqrt{n+m}}(\Sigma_k(Q)-\Sigma_k(P) +(\mu_k(Q)-\mu_k(P))^{\otimes 2})\\
&=:\frac{1}{\sqrt{n+m}}\zeta(P,Q)
\end{align*}
and 
\begin{align*}
&\theta_{\ell} \mathbb{E}_{Y \sim Q_{nm}}[\Psi_{\ell}(Y)]\\
&=\int_{\mathcal{H}} \mathbb{E}_{Y \sim Q_{nm}}[h(Y,y)] \Psi_{\ell}(y) dP(y)\\
&=\int_{\mathcal{H}} \mathbb{E}_{Y \sim Q_{nm}}[\left<(k(\cdot,Y)-\mu_k(P))^{\otimes 2}-\Sigma_k(P),(k(\cdot,y)-\mu_k(P))^{\otimes 2}-\Sigma_k(P)\right>_{H(k)^{\otimes 2}}] \Psi_{\ell}(y) dP(y)\\
&=\int_{\mathcal{H}} \left<\mathbb{E}_{Y \sim Q_{nm}}[(k(\cdot,Y)-\mu_k(P))^{\otimes 2}-\Sigma_k(P)],(k(\cdot,y)-\mu_k(P))^{\otimes 2}-\Sigma_k(P)\right>_{H(k)^{\otimes 2}} \Psi_{\ell}(y) dP(y)\\
&=\frac{1}{\sqrt{n+m}} \int_{\mathcal{H}} \left<\zeta(P,Q), h(\cdot,y)\right>_{H(k)^{\otimes 2}} \Psi_{\ell}(y) dP(y)\\
&=:\frac{1}{\sqrt{n+m}}\zeta_{\ell}(P,Q).
\end{align*}
Hence
\[
\mathbb{E}_{Y \sim Q_{nm}}[\Psi_{\ell}(Y)] =\frac{1}{\sqrt{n+m}} \frac{\zeta_{\ell}(P,Q)}{\theta_{\ell}}
\]
and
\begin{align*}
&V_{Y \sim Q_{nm}}[\Psi_{\ell}(Y)]\\
&=\mathbb{E}_{Y \sim Q_{nm}}[\Psi_{\ell}(Y)^2] -\left\{
\mathbb{E}_{Y \sim Q_{nm}}[\Psi_{\ell}(Y)]
\right\}^2\\
&=\int_{\mathcal{H}} \Psi_{\ell}(y)^2 dQ_{nm}(y)-\frac{1}{n+m}\frac{\zeta_{\ell}(P,Q)^2}{\theta_{\ell}^2}\\
&=\int_{\mathcal{H}} \Psi_{\ell}(y)^2 dP(y) +\frac{1}{\sqrt{n+m}} \int_{\mathcal{H}} \Psi_{\ell}(y)^2 d(Q-P)(y)-\frac{1}{n+m} \frac{\zeta_{\ell}(P,Q)^2}{\theta_{\ell}^2}\\
&=1+\frac{1}{\sqrt{n+m}}\tau_{\ell \ell} -\frac{1}{n+m}\frac{\zeta_{\ell}(P,Q)^2}{\theta_{\ell}^2},
\end{align*}
where
\[
\tau_{\ell s} =\int_{\mathcal{H}} \Psi_{\ell} (y) \Psi_{s} (y) d(Q-P)(y).
\]
Therefore, by using the central limit theorem for $\Psi_{\ell}$s, we get
\begin{align*}
&(n+m)\norm{\widetilde{\Sigma}_k(P)-\widetilde{\Sigma}_k(Q)}^2_{H(k)^{\otimes 2}}\\
&=\frac{n+m}{n^2} \sum_{i,s=1}^{n} h(X_i,X_s)
+\frac{n+m}{m^2} \sum_{j,t=1}^{m} h(Y_j,Y_t)
-\frac{2(n+m)}{nm} \sum_{i=1}^{n} \sum_{j=1}^{m} h(X_i,Y_j) +O_p\left(\frac{1}{\sqrt{n+m}}\right)\\
&=\frac{n+m}{n^2} \sum_{i,s=1}^{n} \sum_{\ell=1}^{\infty} \theta_{\ell} \Psi_{\ell}(X_i) \Psi_{\ell}(X_s)
+\frac{n+m}{m^2} \sum_{j,t=1}^{m} \sum_{\ell=1}^{\infty} \theta_{\ell} \Psi_{\ell}(Y_j) \Psi_{\ell} (Y_t)\\
&~~~~~-\frac{2(n+m)}{nm} \sum_{i=1}^{n} \sum_{j=1}^{m} \sum_{\ell=1}^{\infty} \theta_{\ell} \Psi_{\ell}(X_i) \Psi_{\ell}(Y_j)\\
&=\frac{n+m}{n} \sum_{\ell=1}^{\infty} \theta_{\ell} \left(
\frac{1}{\sqrt{n}}\sum_{i=1}^{n}\Psi_{\ell}(X_i)
\right)^2
+\frac{n+m}{m}\sum_{\ell=1}^{\infty} \theta_{\ell} \left(
\frac{1}{\sqrt{m}} \sum_{j=1}^{m} \Psi_{\ell}(Y_j)
\right)^2\\
&~~~~~
-\frac{2(n+m)}{\sqrt{nm}} \sum_{\ell=1}^{\infty} \theta_{\ell} \left(
\frac{1}{\sqrt{n}} \sum_{i=1}^{n} \Psi_{\ell}(X_i)
\right)\left(
\frac{1}{\sqrt{m}} \sum_{j=1}^{m} \Psi_{\ell}(Y_j)
\right)\\
&\xrightarrow{\mathcal{D}}
\frac{1}{\rho(1-\rho)} \sum_{\ell=1}^{\infty} \theta_{\ell} Z_{\ell}^2,~~Z_{\ell} \sim N\left(\sqrt{\rho (1-\rho)}\cdot \frac{\zeta_{\ell}(P,Q)}{\theta_{\ell}},1\right) ,
\end{align*}
as $n,m \to \infty$.
\subsection{Proof of Lemma \ref{Lemma_three_evaluate}}
(i) First, we prove $\norm{\mu_k(\widehat{Q}_{nm}) -\mu_k(Q_{nm})}_{H(k)}=O_P(1/\sqrt{n+m})$.
For all $\delta>0$, there exists $N_1 \in \mathbb{N}$ such that 
\[
\frac{1}{\sqrt{n+m}} |\mathbb{E}_{Y \sim P} [\norm{k(\cdot,Y)-\mu_k(P)}^2_{H(k)}]+\mathbb{E}_{Y \sim Q}[\norm{k(\cdot,Y)-\mu_k(P)}^2_{H(k)}]| <\frac{\delta}{2}
\]
for all $n,m > N_1$.
In addition, there exists $N_2 \in \mathbb{N}$ such that 
\[
\frac{1}{\sqrt{n+m}} |\mathbb{E}_{Y \sim P} [\norm{k(\cdot,Y)-\mu_k(P)}^2_{H(k)}]+\mathbb{E}_{Y \sim Q}[\norm{k(\cdot,Y)-\mu_k(P)}^2_{H(k)}]| <\frac{\delta}{2}
\]
for all $n,m > N_2$.
We put
\[
M_\delta=\sqrt{\frac{1}{\delta}\left(\frac{1}{1-\rho}+\frac{\delta}{2}\right)\left(\mathbb{E}_{Y \sim P}[\norm{k(\cdot,Y)-\mu_k(P)}^2_{H(k)}]+\frac{\delta}{2}\right)},
\]
and $N_{\delta}=\max\{N_1,N_2\}$.
Since
\[
\mu_k(Q_{nm}) =\left(1-\frac{1}{\sqrt{n+m}}\right) \mu_k(P) +\frac{1}{\sqrt{n+m}}\mu_k(Q)=\mu_k(P)-\frac{1}{\sqrt{n+m}}(\mu_k(P)-\mu_k(Q)),
\]
for all $n,m>N_{\delta}$, we get 
\begin{align*}
&\mathbb{P}\left(\sqrt{n+m}\norm{\mu_k(\widehat{Q}_{nm})-\mu_k(Q_{nm})}_{H(k)}>M_{\delta} \right)\\
& \leq \frac{\mathbb{E}_{Q_{nm}}\left[(n+m)\norm{\mu_k(\widehat{Q}_{nm})-\mu_k(Q_{nm})}^2_{H(k)} \right]}{M_{\delta^2}}\\
&= \frac{(n+m)\mathbb{E}_{Q_{nm}}\left[\norm{\frac{1}{m}\sum_{j=1}^{m} (k(\cdot,Y_j)-\mu_k(Q_{nm}))}^2_{H(k)}\right]}{M_{\delta}^2}\\
&=\frac{t\mathbb{E}_{Y \sim Q_{nm}} [\norm{k(\cdot,Y)-\mu_k(Q_{nm})}^2_{H(k)}]}{mM_{\delta}^2}\\
&=\frac{n+m}{m M_{\delta}^2} \mathbb{E}_{Y \sim Q_{nm}} \left[\norm{k(\cdot,Y)-\mu_k(P)+\frac{1}{\sqrt{n+m}}(\mu_k(P)-\mu_k(Q))}^2_{H(k)} \right]\\
&=\frac{n+m}{mM{\delta}^2} \Big(
\mathbb{E}_{Y \sim Q_{nm}} [\norm{k(\cdot,Y)-\mu_k(P)}^2_{H(k)}] +\frac{2}{\sqrt{n+m}} \mathbb{E}_{Y \sim Q_{nm}} [\left<k(\cdot,Y)-\mu_k(P),\mu_k(P)-\mu_k(Q)\right>_{H(k)}]\\
&~~~~~
+\frac{1}{n+m} \norm{\mu_k(P)-\mu_k(Q)}^2_{H(k)}
\Big)\\
&=\frac{n+m}{mM_{\delta}^2} \left(
\mathbb{E}_{Y \sim Q_{nm}} [\norm{k(\cdot,Y)-\mu_k(P)}^2_{H(k)}] -\frac{1}{n+m} \norm{\mu_k(P)-\mu_k(Q)}^2_{H(k)}
\right)\\
&\leq \frac{n+m}{m M_{\delta}^2} \mathbb{E}_{Y \sim Q_{nm}} [\norm{k(\cdot,Y)-\mu_k(P)}^2_{H(k)}]\\
&= \frac{t}{m M_{\delta}^2}  \left(
\left(1-\frac{1}{\sqrt{n+m}}\right) \mathbb{E}_{X \sim P} [\norm{k(\cdot,Y)-\mu_k(P)}^2_{H(k)}]
+\frac{1}{n+m} \mathbb{E}_{Y \sim Q} [\norm{k(\cdot,Y)-\mu_k(P)}^2_{H(k)}]
\right)\\
&=\frac{n+m}{mM_{\delta}^2} \Big(
\mathbb{E}_{Y \sim P}[\norm{k(\cdot,Y)-\mu_k(P)}^2_{H(k)}] \\
&~~~~~
-\frac{1}{\sqrt{n+m}} \left(
\mathbb{E}_{Y \sim P}[\norm{k(\cdot,Y)-\mu_k(P)}]-\mathbb{E}_{Y \sim Q}[\norm{k(\cdot,Y)-\mu_k(P)}^2_{H(k)}]
\right)
\Big)\\
&< \frac{1}{M_{\delta}^2} \left(\frac{1}{1-\rho}+\frac{\delta}{2}\right) \left(
\mathbb{E}_{Y \sim P} [\norm{k(\cdot,Y)-\mu_k(P)}^2_{H(k)}] +\frac{\delta}{2}
\right)\\
&=\delta.
\end{align*}
Therefore, we obtain $\norm{\mu_k(\widehat{Q}_{nm})-\mu_k(Q_{nm})}_{H(k)}=O_p(1/\sqrt{n+m})$.

(ii) Next, we prove $\norm{\widetilde{\Sigma}_k(Q_{nm})-\Sigma_k(Q_{nm})}_{H(k)^{\otimes 2}}$.
For all $\delta > 0$, we put
\[
M_{\delta}=\sqrt{\frac{1}{\delta}\left(\frac{1}{1-\rho}+\frac{\delta}{2}\right)(\mathbb{E}_{Y \sim P}[\norm{(k(\cdot,Y)-\mu_k(P))^{\otimes 2} -\Sigma_k(P)}^2_{H(k)^{\otimes 2}}] +\delta)}
\]
and	
\begin{align*}
\frac{1}{\sqrt{n+m}}A(n,m,Y):=&\frac{1}{\sqrt{n+m}} \Big\{(k(\cdot,Y)-\mu_k(P)) \otimes (\mu_k(P)-\mu_k(Q)) \\
&~~~~~
+(\mu_k(P)-\mu_k(Q))\otimes (k(\cdot,Y)-\mu_k(P))
+\frac{1}{\sqrt{n+m}} (\mu_k(P)-\mu_k(Q))^{\otimes 2}
\Big\}.
\end{align*}
Then, there exists $N_1 \in \mathbb{N}$ such that
\begin{align*}
&\Big|
\frac{1}{n+m} \mathbb{E}_{Y \sim Q_{nm}} [\norm{A(n,m,Y)}^2_{H(k)^{\otimes 2}}]\\
&~~~~~-\frac{1}{\sqrt{n+m}} \mathbb{E}_{Y \sim Q_{nm}} \left[
\left<(k(\cdot,Y)-\mu_k(P))^{\otimes 2} -\Sigma_k(P),A(t,Y)\right>_{H(k)^{\otimes 2}} 
\right]
\Big| < \frac{\delta}{2}
\end{align*}
for all $n,m > N_1$.
In addition, there exists $N_2 \in \mathbb{N}$ such that 
\begin{align*}
&\Big|
\frac{1}{\sqrt{n+m}}\Big(
\mathbb{E}_{Y \sim P} [\norm{(k(\cdot,Y)-\mu_k(P))^{\otimes 2} -\Sigma_k(P)}^2_{H(k)^{\otimes 2}}] \\
&~~~~~-\mathbb{E}_{Y \sim Q}  [\norm{(k(\cdot,Y)-\mu_k(P))^{\otimes 2} -\Sigma_k(P)}^2_{H(k)^{\otimes 2}}]
\Big)
\Big| < \frac{\delta}{2}
\end{align*}
for all $n,m > N_2$.
and there exists $N_3 \in \mathbb{N}$ such that 
\[
\left|\frac{t}{m} -\frac{1}{1-\rho} \right| < \frac{\delta}{2}
\]
for all $n,m > N_3$.

Let $N_{\delta}=\max\{N_1,N_2,N_3\}$.
For all $n,m > N_{\delta}$, we obtain that
\begin{align*}
&\mathbb{P}(\sqrt{n+m} \norm{\widetilde{\Sigma}_k(Q_{nm})-\Sigma_{k}(Q_{nm})}_{H(k)^{\otimes 2}} > M_{\delta}) \\
& \leq \frac{\mathbb{E}_{Y \sim Q_{nm}}\left[t \norm{\widetilde{\Sigma}_k(Q_{nm})-\Sigma_{k}(Q_{nm})}_{H(k)^{\otimes 2}}^2 \right]}{M_{\delta}^2}\\
&= \frac{t \mathbb{E}_{Y \sim Q_{nm}}[\norm{(k(\cdot,Y)-\mu_k(Q_{nm}))^{\otimes 2} -\Sigma_k(Q_{nm})}^2_{H(k)^{\otimes 2}}]}{mM_{\delta}^2}\\
&=\frac{n+m}{mM_{\delta}^2} \mathbb{E}_{Y \sim Q_{nm}} \Big[
\Big\|
(k(\cdot,Y)-\mu_k(Q_{nm}))^{\otimes 2} -\Sigma_k(P) +\frac{1}{\sqrt{n+m}} (\Sigma_k(P)-\Sigma_k(Q))\\
&~~~~~ -\frac{1}{\sqrt{n+m}} \left(1-\frac{1}{\sqrt{n+m}}\right) (\mu_k(P)-\mu_k(Q))^{\otimes 2} 
\Big\|^2_{H(k)^{\otimes 2}}
\Big]\\
&=\frac{n+m}{mM_{\delta}^2} \mathbb{E}_{Y \sim Q_{nm}} \left[
\norm{(k(\cdot,Y)-\mu_k(Q_{nm}))^{\otimes 2}-\Sigma_k(P)}^2_{H(k)^{\otimes 2}} 
\right]\\
&~~~~~-\frac{1}{n+m} \norm{\Sigma_k(P)-\Sigma_k(Q)+\left(1-\frac{1}{\sqrt{n+m}}\right)(\mu_k(P)-\mu_k(Q))^{\otimes 2}}^2_{H(k)^{\otimes 2}} \\
&\leq \frac{n+m}{mM_{\delta}^2} \mathbb{E}_{Y \sim Q_{nm}} \left[
\norm{(k(\cdot,Y)-\mu_k(Q_{nm}))^{\otimes 2}-\Sigma_k(P)}^2_{H(k)^{\otimes 2}} 
\right]\\
&=\frac{n+m}{mM_{\delta}^2} \mathbb{E}_{Y \sim Q_{nm}} \left[
\norm{\left\{
	k(\cdot,Y)-\mu_k(P)+\frac{1}{\sqrt{n+m}} (\mu_k(P)-\mu_k(Q))
	\right\}^{\otimes 2} -\Sigma_k(P) }^2_{H(k)^{\otimes 2}}
\right]\\
&=\frac{n+m}{mM_{\delta}^2} \mathbb{E}_{Y \sim Q_{nm}} \left[
\norm{(k(\cdot,Y)-\mu_k(P))^{\otimes 2} -\Sigma_k(P) +\frac{1}{\sqrt{n+m}}A(t,Y) }^2_{H(k)^{\otimes 2}}
\right]\\
&=\frac{n+m}{m M_{\delta}^2} \mathbb{E}_{Y \sim Q_{nm}} \Big[
\norm{(k(\cdot,Y)-\mu_k(P))^{\otimes 2} -\Sigma_k(P)}^2_{H(k)^{\otimes 2}} +\frac{1}{n+m} \norm{A(t,Y)}^2_{H(k)^{\otimes 2}} \\
&~~~~~-\frac{2}{\sqrt{n+m}} \left<(k(\cdot,Y)-\mu_k(P))^{\otimes 2} -\Sigma_k(P),A(t,Y)\right>_{H(k)^{\otimes 2}}
\Big]\\
&=\frac{n+m}{m M_{\delta}^2}  \Big\{\mathbb{E}_{Y \sim Q_{nm}} \Big[
\norm{(k(\cdot,Y)-\mu_k(P))^{\otimes 2} -\Sigma_k(P)}^2_{H(k)^{\otimes 2}} \Big] 
+ \frac{1}{n+m} \mathbb{E}_{Y \sim Q_{nm}} \Big[ \norm{A(t,Y)}^2_{H(k)^{\otimes 2}} \Big]\\
&~~~~~-\frac{2}{\sqrt{n+m}} \mathbb{E}_{Y \sim Q_{nm}} \Big[  \left<(k(\cdot,Y)-\mu_k(P))^{\otimes 2} -\Sigma_k(P),A(t,Y)\right>_{H(k)^{\otimes 2}}
\Big] \Big\}\\
&< \frac{n+m}{m M_{\delta}^2} \Big\{
\mathbb{E}_{Y \sim P} \left[
\norm{(k(\cdot,Y)-\mu_k(P))^{\otimes 2} -\Sigma_k(P)}^2_{H(k)^{\otimes 2}}
\right]\\
&~~~~~~
-\frac{1}{\sqrt{n+m}} \Big(
\mathbb{E}_{Y \sim P} \left[
\norm{(k(\cdot,Y)-\mu_k(P))^{\otimes 2} -\Sigma_k(P)}^2_{H(k)^{\otimes 2}}
\right]\\
&~~~~~
-\mathbb{E}_{Y \sim Q} \left[
\norm{(k(\cdot,Y)-\mu_k(P))^{\otimes 2} -\Sigma_k(P)}^2_{H(k)^{\otimes 2}}
\right]
\Big)+\frac{\delta}{2} 
\Big\}\\
&< \frac{1}{M_{\delta}^2} \left(\frac{1}{1-\rho} +\frac{\delta}{2} \right) \left\{
\mathbb{E}_{Y \sim P} [\norm{(k(\cdot,Y)-\mu_k(P))^{\otimes 2} -\Sigma_k(P)}^2_{H(k)^{\otimes 2}} +\delta]
\right\}\\
&=\delta.
\end{align*}
Therefore, $\norm{\widetilde{\Sigma}_k(Q_{nm})-\Sigma_k(Q_{nm})}_{H(k)^{\otimes 2}}=O_p(1/\sqrt{n+m})$ is proved.

(iii) Finally, we prove $\norm{\Sigma_k(\widehat{Q}_{nm})-\Sigma_k(Q_{nm})}_{H(k)^{\otimes 2}}=O_p(1/\sqrt{n+m})$.
We get 
\begin{align*}
&\Sigma_k(\widehat{Q}_{nm})-\Sigma_k(Q_{nm})\\ 
&=\frac{1}{m} \sum_{j=1}^{m}(k(\cdot,Y_j)-\mu_k(\widehat{Q}_{nm}))^{\otimes 2} -\Sigma_k(Q_{nm}) \\
&=\frac{1}{m} \sum_{j=1}^{m} (k(\cdot,Y_j)-\mu_k(Q_{nm})+\mu_k(Q_{nm})-\mu_k(\widehat{Q}_{nm}))^{\otimes 2} -\Sigma_k(Q_{nm})\\
&=\frac{1}{m} \sum_{j=1}^{m} (k(\cdot,Y_j)-\mu_k(Q_{nm}))^{\otimes 2} -\Sigma_k(Q_{nm})\\
&~~~~~+\frac{1}{m} \sum_{j=1}^{m} (k(\cdot,Y_j)-\mu_k(Q_{nm})) \otimes (\mu_k(Q_{nm})-\mu_k(\widehat{Q}_{nm}))\\
&~~~~~+\frac{1}{m} \sum_{j=1}^{m} (\mu_k(Q_{nm})-\mu_k(\widehat{Q}_{nm})) \otimes (k(\cdot,Y_j)-\mu_k(Q_{nm}))
+(\mu_k(Q_{nm}) -\mu_k(\widehat{Q}_{nm}))^{\otimes 2}\\
&=\widetilde{\Sigma}_k(Q_{nm}) -\Sigma_k(Q_{nm}) -(\mu_k(Q_{nm})-\mu_k(\widehat{Q}_{nm}))^{\otimes 2}
\end{align*}
by a expansion of $\widetilde{\Sigma}_k(Q_{nm})$.
Using (i) and (ii) leads to the following:
\begin{align*}
&\norm{\Sigma_k(\widehat{Q}_{nm})-\Sigma_k(Q_{nm})}_{H(k)^{\otimes 2}}\\
&\leq \norm{\widetilde{\Sigma}_k(Q_{nm})-\Sigma_k(Q_{nm})}_{H(k)^{\otimes
		2}} +\norm{\mu_k(Q_{nm})-\mu_k(\widehat{Q}_{nm})}^2_{H(k)^{\otimes 2}}\\
&=O_p\left(\frac{1}{\sqrt{n+m}}\right) +O_p\left(\frac{1}{n+m}\right)\\
&=O_p\left(\frac{1}{\sqrt{n+m}}\right).
\end{align*}
\subsection{Proof of Lemma \ref{asymptotic_equivalent_contiguous}}
First we have
\begin{align}
	&\Sigma_k(Q_{nm}) \nonumber\\
	&=\mathbb{E}_{Y \sim Q_{nm}} \left[
	(k(\cdot,Y)-\mu_k(Q_{nm}))^{\otimes 2}
	\right] \nonumber\\
	&=\left(1-\frac{1}{\sqrt{n+m}}\right)\mathbb{E}_{Y \sim P} \left[
	(k(\cdot,Y)-\mu_k(Q_{nm}))^{\otimes 2}
	\right]
	+\frac{1}{\sqrt{n+m}} \mathbb{E}_{Y \sim Q} \left[
	(k(\cdot,Y)-\mu_k(Q_{nm}))^{\otimes 2}
	\right] \nonumber\\
	&=\left(1-\frac{1}{\sqrt{n+m}}\right)\mathbb{E}_{Y \sim P} \left[
	(k(\cdot,Y)-\mu_k(P)+\mu_k(P)-\mu_k(Q_{nm}))^{\otimes 2}
	\right] \nonumber\\
	&~~~~~+\frac{1}{\sqrt{n+m}} \mathbb{E}_{Y \sim Q} \left[
	(k(\cdot,Y)-\mu_k(Q)+\mu_k(Q)-\mu_k(Q_{nm}))^{\otimes 2}
	\right] \nonumber\\
	&=\left(1-\frac{1}{\sqrt{n+m}}\right) \left\{\mathbb{E}_{Y \sim P} \left[
	(k(\cdot,Y)-\mu_k(P))^{\otimes 2} 
	\right] +(\mu_k(P)-\mu_k(Q_{nm}))^{\otimes 2} \right\} \nonumber\\
	&~~~~~+\frac{1}{\sqrt{n+m}} \left\{ \mathbb{E}_{Y \sim Q} \left[
	(k(\cdot,Y)-\mu_k(Q))^{\otimes 2} 
	\right] +(\mu_k(Q)-\mu_k(Q_{nm}))^{\otimes 2} \right\} \nonumber\\
	&=\left(1-\frac{1}{\sqrt{n+m}}\right) \Sigma_k(P) +\frac{1}{\sqrt{n+m}} \Sigma_k(Q) 
	+\left(1-\frac{1}{\sqrt{n+m}}\right) \frac{1}{n+m} (\mu_k(P)-\mu_k(Q))^{\otimes 2} \nonumber\\
	&~~~~~
	+\frac{1}{\sqrt{n+m}} \left(1-\frac{1}{\sqrt{n+m}}\right)^2 (\mu_k(P)-\mu_k(Q))^{\otimes 2} \nonumber\\
	&=\left(1-\frac{1}{\sqrt{n+m}}\right) \Sigma_k(P) +\frac{1}{\sqrt{n+m}} \Sigma_k(Q) +\frac{1}{\sqrt{n+m}} \left(1-\frac{1}{\sqrt{n+m}}\right) (\mu_k(P)-\mu_k(Q))^{\otimes 2}. \label{Eq_Lemm_SigQnm}
\end{align}
The result (\ref{Eq_Lemm_SigQnm}) leads to
\begin{align*}
	&\Sigma_k(P)-\Sigma_k(Q_{nm})\\
	&=\Sigma_k(P)-\left(1-\frac{1}{\sqrt{n+m}}\right) \Sigma_k(P)-\frac{1}{\sqrt{n+m}} \Sigma_k(Q)\\
	&~~~~~-\frac{1}{\sqrt{n+m}} \left(1-\frac{1}{\sqrt{n+m}}\right) (\mu_k(P)-\mu_k(Q))^{\otimes 2}\\
	&=\frac{1}{\sqrt{n+m}}\Sigma_k(P)-\frac{1}{\sqrt{n+m}}\Sigma_k(Q) -\frac{1}{\sqrt{n+m}} \left(1-\frac{1}{\sqrt{n+m}}\right) (\mu_k(P)-\mu_k(Q))^{\otimes 2}\\
	&=\frac{1}{\sqrt{n+m}}(\Sigma_k(P)-\Sigma_k(Q))- \frac{1}{\sqrt{n+m}}\left(1-\frac{1}{\sqrt{n+m}}\right) (\mu_k(P)-\mu_k(Q))^{\otimes 2}.
\end{align*}
Hence
\begin{align*}
	&(n+m)\left(
	\norm{\Sigma_k(\widehat{P})-\Sigma_k(\widehat{Q}_{nm})}^2_{H(k)^{\otimes 2}}-\norm{\widetilde{\Sigma}_k(P)-\widetilde{\Sigma}_k(Q_{nm})}^2_{H(k)^{\otimes 2}}
	\right)\\
	&\leq (n+m)\left\{
	\norm{\Sigma_k(\widehat{P})-\Sigma_k(\widehat{Q}_{nm})}_{H(k)^{\otimes 2}}+\norm{\widetilde{\Sigma}_k(P)-\widetilde{\Sigma}_k(Q_{nm})}_{H(k)^{\otimes 2}}
	\right\} \\
	&~~~~~\times \norm{(\mu_k(P)-\mu_k(\widehat{P}))^{\otimes 2}-(\mu_k(Q_{nm})-\mu_k(\widehat{Q}_{nm}))^{\otimes 2}}_{H(k)^{\otimes 2}}\\
	&\leq \sqrt{n+m}\left\{
	\norm{\Sigma_k(\widehat{P})-\Sigma_k(\widehat{Q}_{nm})}_{H(k)^{\otimes 2}}+\norm{\widetilde{\Sigma}_k(P)-\widetilde{\Sigma}_k(Q_{nm})}_{H(k)^{\otimes 2}}
	\right\} \\
	&~~~~~\times\left(
	\frac{\sqrt{n+m}}{n}n\norm{\mu_k(P)-\mu_k(\widehat{P})}^2_{H(k)} +\frac{\sqrt{n+m}}{m}m\norm{\mu_k(Q_{nm})-\mu_k(\widehat{Q}_{nm})}^2_{H(k)} 
	\right)\\
	&=\sqrt{n+m}\left\{
	\norm{\Sigma_k(\widehat{P})-\Sigma_k(P)-(\Sigma_k(\widehat{Q}_{nm})-\Sigma_k(Q_{nm}))+\Sigma_k(P)-\Sigma_k(Q_{nm})}_{H(k)^{\otimes 2}} \right.\\
	&\left.~~~~~+\norm{\widetilde{\Sigma}_k(P)-\Sigma_k(P)-(\widetilde{\Sigma}_k(Q_{nm})-\Sigma_k(Q_{nm}))+\Sigma_k(P)-\Sigma_k(Q_{nm})}_{H(k)^{\otimes 2}}
	\right\} \\
	&~~~~~\times\left(
	\frac{\sqrt{n+m}}{n}n\norm{\mu_k(P)-\mu_k(\widehat{P})}^2_{H(k)} +\frac{\sqrt{n+m}}{m}m\norm{\mu_k(Q_{nm})-\mu_k(\widehat{Q}_{nm})}^2_{H(k)} 
	\right)\\
	&=\sqrt{n+m}\left\{
	\Bigg\|\Sigma_k(\widehat{P})-\Sigma_k(P)-(\Sigma_k(\widehat{Q}_{nm})-\Sigma_k(Q_{nm}))  \right.\\
	&~~~~~
	+\frac{1}{\sqrt{n+m}}(\Sigma_k(P)-\Sigma_k(Q)) -\frac{1}{\sqrt{n+m}} \left(1-\frac{1}{\sqrt{n+m}}\right) (\mu_k(P)-\mu_k(Q))^{\otimes 2}\Bigg\|_{H(k)^{\otimes 2}} \\
	&~~~~~+\Bigg\|\widetilde{\Sigma}_k(P)-\Sigma_k(P)-(\widetilde{\Sigma}_k(Q_{nm})-\Sigma_k(Q_{nm}))\\
	&\left.~~~~~+\frac{1}{\sqrt{n+m}}(\Sigma_k(P)-\Sigma_k(Q))-+\frac{1}{\sqrt{n+m}} \left(1-\frac{1}{\sqrt{n+m}}\right) (\mu_k(P)-\mu_k(Q))^{\otimes 2} \Bigg\|_{H(k)^{\otimes 2}} 
	\right\} \\
	&~~~~~\times\left(
	\frac{\sqrt{n+m}}{n}n\norm{\mu_k(P)-\mu_k(\widehat{P})}^2_{H(k)} +\sqrt{n+m}\norm{\mu_k(Q_{nm})-\mu_k(\widehat{Q}_{nm})}^2_{H(k)} 
	\right)\\
	&\leq \Bigg\{
	\sqrt{\frac{n+m}{n}} \sqrt{n}\norm{\Sigma_k(\widehat{P})-\Sigma_k(P)}_{H(k)^{\otimes 2}}+\sqrt{n+m} \norm{\Sigma_k(\widehat{Q}_{nm})-\Sigma_k(Q_{nm})}_{H(k)^{\otimes 2}} \\
	&~~~~~+\sqrt{\frac{n+m}{n}}\sqrt{n}\norm{\widetilde{\Sigma}_k(P)-\Sigma_k(P)}_{H(k)^{\otimes 2}}+\sqrt{n+m}\norm{\widetilde{\Sigma}_k(Q_{nm})-\Sigma_k(Q_{nm})}_{H(k)^{\otimes 2}}\\ 
	&~~~~~+2\norm{\Sigma_k(P)-\Sigma_k(Q)}_{H(k)^{\otimes 2}}
	+2\left(1-\frac{1}{\sqrt{n+m}}\right) \norm{\mu_k(P)-\mu_k(Q)}^2_{H(k)} 
	\Bigg\}  \\
	&~~~~~\times\left(
	\frac{\sqrt{n+m}}{n}n\norm{\mu_k(P)-\mu_k(\widehat{P})}^2_{H(k)} +\sqrt{n+m}\norm{\mu_k(Q_{nm})-\mu_k(\widehat{Q}_{nm})}^2_{H(k)} 
	\right).
\end{align*}
Therefore, we obtain 
\[
(n+m)\left(
\norm{\Sigma_k(\widehat{P})-\Sigma_k(\widehat{Q}_{nm})}^2_{H(k)^{\otimes 2}}-\norm{\widetilde{\Sigma}_k(P)-\widetilde{\Sigma}_k(Q_{nm})}^2_{H(k)^{\otimes 2}}
\right) \xrightarrow{P}0,
\]
as $n,m \to \infty$ by Lemma \ref{Lemma_three_evaluate}, which completes the proof of Lemma \ref{asymptotic_equivalent_contiguous}.
\subsection{Proof of Proposition \ref{Prop_eigenvalue_T1_and_T3}}
Let $\Upsilon= V[(k(\cdot,x)-\mu_k(P))^{\otimes 2}]$ be the operator with spectral representation $\Upsilon=\sum_{\ell=1}^{\infty} \theta_{\ell} \phi^{\otimes 2}_{\ell}$, and  $T: H(k)^{\otimes 2} \to
L_2(\mathcal{H},P)$ be the operator
\[
(T(A))(x) =\left<\Upsilon^{-1/2} \{(k(\cdot,x)-\mu_k(P))^{\otimes 2} -\Sigma_k(P)\},A\right>_{H(k)^{\otimes 2}}
\]
for all $A \in H(k)^{\otimes 2}$ and $x \in \mathcal{H}$.
We consider the adjoint operator of this $T$,
\begin{align*}
	\left<T^* g,A\right>_{H(k)^{\otimes 2}}
	&=\left<g, TA\right>_{L_2(\mathcal{H},P)}\\
	&=\int_{\mathcal{H}} (T(A))(y) g(y) dP(y)\\
	&=\int_{\mathcal{H}} \left<\Upsilon^{-1/2} \{(k(\cdot,y)-\mu_k(P))^{\otimes 2} -\Sigma_k(P)\},A\right>_{H(k)^{\otimes 2}} g(y) dP(y)\\
	&=\left<\int_{\mathcal{H}}  (\Upsilon^{-1/2} \{(k(\cdot,y)-\mu_k(P))^{\otimes 2} -\Sigma_k(P)\}) g(y) dP(y),A\right>_{H(k)^{\otimes 2}}
\end{align*}
for all $g \in L_2(\mathcal{H},P)$ and $A \in H(k)^{\otimes 2}$, hence we get that adjoint operator $T^*$ of $T$ is
\[
T^* g=\int_{\mathcal{H}}  (\Upsilon^{-1/2} \{(k(\cdot,y)-\mu_k(P))^{\otimes 2} -\Sigma_k(P)\}) g(y) dP(y)
\]
for all $g \in L_2(\mathcal{H},P)$.
Furthermore, since
\begin{align*}
	T^* (T(A))
	&=\int_{\mathcal{H}}  (\Upsilon^{-1/2} \{(k(\cdot,y)-\mu_k(P))^{\otimes 2} -\Sigma_k(P)\}) (T(A))(y) dP(y)\\
	&=\int_{\mathcal{H}}  (\Upsilon^{-1/2} \{(k(\cdot,y)-\mu_k(P))^{\otimes 2} -\Sigma_k(P)\}) \left<\Upsilon^{-1/2} \{(k(\cdot,y)-\mu_k(P))^{\otimes 2} -\Sigma_k(P)\},A\right>_{H(k)^{\otimes 2}} dP(y)\\
	&=\int_{\mathcal{H}}  (\Upsilon^{-1/2} \{(k(\cdot,y)-\mu_k(P))^{\otimes 2} -\Sigma_k(P)\})^{\otimes 2} A dP(y)\\
	&=\Upsilon^{-1/2} \int_{\mathcal{H}} ( \{(k(\cdot,y)-\mu_k(P))^{\otimes 2} -\Sigma_k(P)\})^{\otimes 2}dP(y) \Upsilon^{-1/2} A \\
	&=\Upsilon^{-1/2} \Upsilon \Upsilon^{-1/2} A\\
	&=A
\end{align*}
for all $A \in H(k)^{\otimes 2}$, $T^*T$ is the identity operator from $H(k)^{\otimes 2}$ to $H(k)^{\otimes 2}$.
It follows from direct calculations for $T \Upsilon T^* : L_2(\mathcal{H},P) \to L_2(\mathcal{H},P)$ that 
\begin{align*}
	&(T\Upsilon T^*g) (x)\\
	&=\left<\Upsilon^{-1/2} \{(k(\cdot,x)-\mu_k(P))^{\otimes 2} -\Sigma_k(P)\}, \Upsilon T^* g\right>_{H(k)^{\otimes 2}}\\
	&=\left<\Upsilon^{1/2}\{(k(\cdot,x)-\mu_k(P))^{\otimes 2} -\Sigma_k(P)\},  T^* g\right>_{H(k)^{\otimes 2}}\\
	&=\left<T\Upsilon^{1/2}\{(k(\cdot,x)-\mu_k(P))^{\otimes 2} -\Sigma_k(P)\},   g\right>_{L_2(\mathcal{H},P)}\\
	&=\int_{\mathcal{H}}  [T\Upsilon^{1/2}\{(k(\cdot,x)-\mu_k(P))^{\otimes 2} -\Sigma_k(P)\}](y) g(y) dP(y)\\
	&=\int_{\mathcal{H}} \left<\Upsilon^{-1/2} \{(k(\cdot,y)-\mu_k(P))^{\otimes 2} -\Sigma_k(P)\},\Upsilon^{1/2}\{(k(\cdot,x)-\mu_k(P))^{\otimes 2} -\Sigma_k(P)\}\right>_{H(k)^{\otimes 2}} g(y)  dP(y)\\
	&=\int_{\mathcal{H}} \left<(k(\cdot,y)-\mu_k(P))^{\otimes 2} -\Sigma_k(P),(k(\cdot,x)-\mu_k(P))^{\otimes 2} -\Sigma_k(P)\right>_{H(k)^{\otimes 2}} g(y)  dP(y)\\
	&=(S_kg)(x)
\end{align*}
for all $g \in L_2(\mathcal{H},P)$ an $x \in \mathcal{H}$, thus we see that  $T\Upsilon T^*=S_k$.
Therefore, $\theta_{\ell}$ and $T\phi_{\ell}$ are the eigenvalue and eigenfunction of $S_k$, by the following that
\begin{align*}
	S_kg
	&= T\Upsilon T^* g\\
	&=T \sum_{\ell=1}^{\infty} \theta_{\ell} \phi_{\ell}^{\otimes 2} T^* g\\
	&=\sum_{\ell=1}^{\infty} \theta_{\ell} \left<\phi_{\ell},T^* g\right>_{H(k)^{\otimes 2}}T\phi_{\ell}\\
	&=\sum_{\ell=1}^{\infty} \theta_{\ell} \left<T\phi_{\ell}, g\right>_{H(k)^{\otimes 2}}T\phi_{\ell}\\
	&=\sum_{\ell=1}^{\infty} \theta_{\ell} (T\phi_{\ell})^{\otimes 2} g
\end{align*}
and $\{T\phi_{\ell}\}_{\ell=1}^{\infty}$ is an orthonormal system in $L_2(\mathcal{H},P)$ which holds
\[
\left<T \phi_{\ell},T\phi_{s}\right>_{L_2(\mathcal{H},P)}=\left<T^*T\phi_{\ell}, \phi_{s}\right>_{H(k)^{\otimes 2}}=\left<\phi_{\ell}, \phi_{s}\right>_{H(k)^{\otimes 2}}=\delta_{\ell s}
\]
for all $\ell, s \in \mathbb{N}$.
In fact,
\[
S_k (T \phi_{\ell})
=T \Upsilon T^* T \phi_{\ell}
=T \Upsilon \phi_{\ell}
=T \theta_{\ell}  \phi_{\ell}
=\theta_{\ell} (T \phi_{\ell})
\]
shows that $\theta_{\ell}$ and $T\theta_{\ell}$ are eigenvalue and eigenfunction of $S_k$.

\subsection{Proof of Theorem \ref{Asymptotic_distribution_Gram_matrix_spectrum}}
Let $W'_n=\sum_{\ell=1}^{\infty}({\widehat{\lambda}}^{(n)}_{\ell}-\lambda_{\ell})Z_{\ell}^2$.
Then
\begin{align}
&\mathbb{E}[W'_n]=\mathbb{E}\left[\sum_{\ell=1}^{\infty}({\widehat{\lambda}}^{(n)}_{\ell}-\lambda_{\ell})\right], \nonumber\\
&V[W'_n]=2 \mathbb{E}\left[\sum_{\ell=1}^{\infty}({\widehat{\lambda}}^{(n)}_{\ell}-\lambda_{\ell})^2\right]+Cov\left(\sum_{\ell=1}^{\infty}({\widehat{\lambda}}^{(n)}_{\ell}-\lambda_{\ell}),\sum_{\ell'=1}^{\infty}(\widehat{\lambda}_{\ell'}-\lambda_{\ell'})\right). \label{Eq_VW_Theorem4}
\end{align}
By the definition of trace of the Hilbert--Schmidt operator, we get the following inequality
\begin{equation}\label{Ineq_lambda_Upsilon}
\sum_{\ell=1}^{\infty}({\widehat{\lambda}}^{(n)}_{\ell}-\lambda_{\ell})
=\text{tr}[\widehat{\Upsilon}^{(n)}]-\text{tr}[\Upsilon]
=\left<\widehat{\Upsilon}^{(n)}-\Upsilon,I\right>_{H(k)^{\otimes 4}}
\leq \norm{\widehat{\Upsilon}^{(n)}-\Upsilon}_{H(k)^{\otimes 4}}.
\end{equation}
By using a notation
\[
B(X_1)=(k(\cdot,X_1)-\mu_k(P)) \otimes (\mu_k(P)-\mu_k(\widehat{P}))
+(\mu_k(P)-\mu_k(\widehat{P}))\otimes (k(\cdot,X_1)-\mu_k(P))
+(\mu_k(P)-\mu_k(\widehat{P}))^{\otimes 2},
\]
we have
\begin{align}
\widehat{\Upsilon}^{(n)} \nonumber
&=\frac{1}{n}\sum_{i=1}^{n}\left\{
(k(\cdot,X_i)-\mu_k(\widehat{P}))^{\otimes 2} -\Sigma_k(\widehat{P})
\right\}^{\otimes 2}\nonumber\\
&=\frac{1}{n}\sum_{i=1}^{n}\left\{(k(\cdot,X_i)-\mu_k(P))^{\otimes 2}-\Sigma_k(P) +B(X_i)+\Sigma_k(P)-\Sigma_k(\widehat{P})\right\}^{\otimes 2}\nonumber\\
&=\frac{1}{n}\sum_{i=1}^{n}\left\{(k(\cdot,X_i)-\mu_k(P))^{\otimes 2}-\Sigma_k(P)\right\}^{\otimes 2}\nonumber\\
&~~~~~+\frac{1}{n}\sum_{i=1}^{n}\left\{(k(\cdot,X_i)-\mu_k(P))^{\otimes 2}-\Sigma_k(P)\right\} \otimes \left\{
B(X_i)+\Sigma_k(P)-\Sigma_k(\widehat{P})
\right\}\nonumber\\
&~~~~~+\frac{1}{n}\sum_{i=1}^{n}\left\{
B(X_i)+\Sigma_k(P)-\Sigma_k(\widehat{P})
\right\} \otimes \left\{(k(\cdot,X_i)-\mu_k(P))^{\otimes 2}-\Sigma_k(P)\right\}\nonumber\\
&~~~~~+\frac{1}{n}\sum_{i=1}^{n}\left\{
B(X_i)+\Sigma_k(P)-\Sigma_k(\widehat{P})
\right\}^{\otimes 2} \nonumber\\
&=:I_1+I_2+I_3+I_4. \label{Upsilonhat_equality_1}
\end{align}
Since direct calculation, we get the following three expressions,
\begin{align}
I_2
&=\frac{1}{n}\sum_{i=1}^{n}\left\{(k(\cdot,X_i)-\mu_k(P))^{\otimes 2}-\Sigma_k(P)\right\} \otimes \left\{
B(X_i)+\Sigma_k(P)-\Sigma_k(\widehat{P})
\right\}\nonumber\\
&=\frac{1}{n}\sum_{i=1}^{n}\left\{(k(\cdot,X_i)-\mu_k(P))^{\otimes 2}-\Sigma_k(P)\right\} \otimes B(X_i)\nonumber\\
&~~~~~
+\frac{1}{n}\sum_{i=1}^{n}\left\{(k(\cdot,X_i)-\mu_k(P))^{\otimes 2}-\Sigma_k(P)\right\} \otimes \left(\Sigma_k(P)-\Sigma_k(\widehat{P})
\right)\nonumber\\
&=\frac{1}{n}\sum_{i=1}^{n}\left\{(k(\cdot,X_i)-\mu_k(P))^{\otimes 2}-\Sigma_k(P)\right\} \otimes B(X_i)
-\left(\Sigma_k(P)-\Sigma_k(\widehat{P})
\right)^{\otimes 2}, \label{Upsilonhat_equality_2}
\end{align}
\begin{align}
I_3&=\frac{1}{n}\sum_{i=1}^{n}\left\{
B(X_i)+\Sigma_k(P)-\Sigma_k(\widehat{P})
\right\} \otimes \left\{
(k(\cdot,X_i)-\mu_k(P))^{\otimes 2}-\Sigma_k(P)
\right\}\nonumber\\
&=\frac{1}{n}\sum_{i=1}^{n}B(X_i) \otimes \left\{
(k(\cdot,X_i)-\mu_k(P))^{\otimes 2}-\Sigma_k(P)
\right\}
-\left(
\Sigma_k(P)-\Sigma_k(\widehat{P})
\right)^{\otimes 2} \label{Upsilonhat_equality_3}
\end{align}	
and
\begin{align}
I_4=&\frac{1}{n}\sum_{i=1}^{n}\left\{
B(X_i)+\Sigma_k(P)-\Sigma_k(\widehat{P})
\right\}^{\otimes 2}\nonumber\\
&=\frac{1}{n}\sum_{i=1}^{n}B(X_i)^{\otimes 2}
+\frac{1}{n}\sum_{i=1}^{n}B(X_i) \otimes \left(\Sigma_k(P)-\Sigma_k(\widehat{P})\right)\nonumber\\
&~~~~~
+\frac{1}{n}\sum_{i=1}^{n}(\Sigma_k(P)-\Sigma_k(\widehat{P})) \otimes B(X_i)
+(\Sigma_k(P)-\Sigma_k(\widehat{P}))^{\otimes 2}\nonumber\\
&=\frac{1}{n}\sum_{i=1}^{n}B(X_i)^{\otimes 2}
-(\mu_k(P)-\mu_k(\widehat{P}))^{\otimes 2} \otimes (\Sigma_k(P)-\Sigma_k(\widehat{P}))\nonumber\\
&~~~~~
-(\Sigma_k(P)-\Sigma_k(\widehat{P})) \otimes (\mu_k(P)-\mu_k(\widehat{P}))^{\otimes 2}
+(\Sigma_k(P)-\Sigma_k(\widehat{P}))^{\otimes 2}. \label{Upsilonhat_equality_4}
\end{align}
By (\ref{Upsilonhat_equality_2}), (\ref{Upsilonhat_equality_3}) and (\ref{Upsilonhat_equality_4}) are combined into (\ref{Upsilonhat_equality_1}), we have an expression
\begin{align*}
\widehat{\Upsilon}^{(n)}
&=\frac{1}{n}\sum_{i=1}^{n}\left\{(k(\cdot,X_i)-\mu_k(P))^{\otimes 2} -\Sigma_k(P)\right\}^{\otimes 2}
-(\Sigma_k(P)-\Sigma_k(\widehat{P}))^{\otimes 2}
\\
&~~~~~
-(\mu_k(P)-\mu_k(\widehat{P}))^{\otimes 2} \otimes (\Sigma_k(P)-\Sigma_k(\widehat{P}))
-(\Sigma_k(P)-\Sigma_k(\widehat{P})) \otimes (\mu_k(P)-\mu_k(\widehat{P}))^{\otimes 2}
\\
&~~~~~
+\frac{1}{n}\sum_{i=1}^{n}B(X_i)^{\otimes 2}
+\frac{1}{n}\sum_{i=1}^{n}B(X_i) \otimes \left\{
(k(\cdot,X_i)-\mu_k(P))^{\otimes 2}-\Sigma_k(P)
\right\}\\
&~~~~~
+\frac{1}{n}\sum_{i=1}^{n} \left\{
(k(\cdot,X_i)-\mu_k(P))^{\otimes 2}-\Sigma_k(P)
\right\}\otimes B(X_i).
\end{align*}
Therefore, 
\begin{align}
&\norm{\widehat{\Upsilon}^{(n)}-\Upsilon}_{H(k)^{\otimes 4}} \nonumber\\
&\leq \norm{\frac{1}{n}\sum_{i=1}^{n}\{(k(\cdot,X_i)-\mu_k(P))^{\otimes 2} -\Sigma_k(P)\}^{\otimes 2} -\Upsilon}_{H(k)^{\otimes 4}}
+\norm{(\Sigma_k(P)-\Sigma_k(\widehat{P}))^{\otimes 2}}_{H(k)^{\otimes 4}}\nonumber\\
&~~~~+\norm{(\mu_k(P)-\mu_k(\widehat{P}))^{\otimes 2}\otimes (\Sigma_k(P)-\Sigma_k(\widehat{P}))}_{H(k)^{\otimes 4}}
+\norm{(\Sigma_k(P)-\Sigma_k(\widehat{P}))\otimes (\mu_k(P)-\mu_k(\widehat{P}))^{\otimes 2}}_{H(k)^{\otimes 4}}\nonumber\\
&~~~~~+\norm{\frac{1}{n}\sum_{i=1}^{n}B(X_i)^{\otimes 2}}_{H(k)^{\otimes 4}}
+\norm{\frac{1}{n}\sum_{i=1}^{n}B(X_i)\otimes \{(k(\cdot,X_i)-\mu_k(P))^{\otimes 2}-\Sigma_k(P)\}}_{H(k)^{\otimes 4}}\nonumber\\
&~~~~~
+\norm{\frac{1}{n}\sum_{i=1}^{n} \{(k(\cdot,X_i)-\mu_k(P))^{\otimes 2}-\Sigma_k(P)\} \otimes B(X_i) }_{H(k)^{\otimes 4}}\nonumber\\
&=:I_1+I_2+I_3+I_4+I_5+I_6+I_7. \label{Upsilonhat-Upsilon}
\end{align}
By the law of large number in Hilbert spaces (see \cite{hoffmann-jorgensen1976}), we have
\begin{equation} \label{I1_Theorem4}
I_1=\norm{\frac{1}{n}\sum_{i=1}^{n}\{(k(\cdot,X_i)-\mu_k(P))^{\otimes 2} -\Sigma_k(P)\}^{\otimes 2} -\Upsilon}_{H(k)^{\otimes 4}}=O_p\left(\frac{1}{\sqrt{n}}\right),
\end{equation}
and we get
\begin{equation}\label{I2_Theorem4}
I_2=\norm{\Sigma_k(P)-\Sigma_k(\widehat{P})}^2_{H(k)^{\otimes 2}}=O_p\left(\frac{1}{n}\right),
\end{equation}
\begin{equation}\label{I3_Theorem4}
I_3=\norm{\mu_k(P)-\mu_k(\widehat{P})}^2_{H(k)} \norm{\Sigma_k(P)-\Sigma_k(\widehat{P})}_{H(k)^{\otimes 2}}=O_p\left(\frac{1}{n\sqrt{n}}\right),
\end{equation}
and
\begin{equation}\label{I4_Theorem4}
I_4=\norm{\mu_k(P)-\mu_k(\widehat{P})}^2_{H(k)} \norm{\Sigma_k(P)-\Sigma_k(\widehat{P})}_{H(k)^{\otimes 2}}=O_p\left(\frac{1}{n\sqrt{n}}\right).
\end{equation}
Next our focus goes to $I_5$. 
We have by direct computations that
\begin{align*}
I_5
&\leq \frac{1}{n}\sum_{i=1}^{n} \norm{B(X_i)}^2_{H(k)^{\otimes 2}} \\
&=\frac{1}{n}\sum_{i=1}^{n} \left\{
\norm{k(\cdot,X_i)-\mu_k(P)}_{H(k)} \norm{\mu_k(P)-\mu_k(\widehat{P})}_{H(k)}
+\norm{\mu_k(P)-\mu_k(\widehat{P})}_{H(k)}
\norm{k(\cdot,X_i)-\mu_k(P)}_{H(k)} \right.\\
&~~~~~~~~~~\left.
+\norm{\mu_k(P)-\mu_k(\widehat{P})}^2_{H(k)}
\right\}^2\\
&=\frac{1}{n}\sum_{i=1}^{n} \left\{ \norm{\mu_k(P)-\mu_k(\widehat{P})}^2_{H(k)}
+2\norm{k(\cdot,X_i)-\mu_k(P)}_{H(k)} \norm{\mu_k(P)-\mu_k(\widehat{P})}_{H(k)}	\right\}^2\\
&=\norm{\mu_k(P)-\mu_k(\widehat{P})}^4_{H(k)} +\left\{\frac{4}{n}\sum_{i=1}^{n} \norm{k(\cdot,X_i)-\mu_k(P)}_{H(k)} \right\} \norm{\mu_k(P)-\mu_k(\widehat{P})}^3_{H(k)}\\
&~~~~~+\left\{\frac{4}{n}\sum_{i=1}^{n} \norm{k(\cdot,X_i)-\mu_k(P)}_{H(k)}^2 \right\} \norm{\mu_k(P)-\mu_k(\widehat{P})}_{H(k)}^2. 
\end{align*}
Since Condition is assumed, the followings  hold
\[
\frac{1}{n}\sum_{i=1}^{n} \norm{k(\cdot,X_i)-\mu_k(P)}_{H(k)}^2 \xrightarrow{P} \mathbb{E}[\norm{k(\cdot,X_i)-\mu_k(P)}_{H(k)}^2]
\]
and
\[
\frac{1}{n}\sum_{i=1}^{n} \norm{k(\cdot,X_i)-\mu_k(P)}_{H(k)} \xrightarrow{P} \mathbb{E}[\norm{k(\cdot,X_i)-\mu_k(P)}_{H(k)}]
\]
as $n \to \infty$ by the law of large numbers.
Hence, we get
\begin{align}
I_5
&=O_p\left(\frac{1}{n^2}\right) +\frac{1}{n}\sum_{i=1}^{n}\norm{k(\cdot,X_i)-\mu_k(P)}_{H(k)} O_p\left(\frac{1}{n\sqrt{n}}\right) +\frac{1}{n}\sum_{i=1}^{n}\norm{k(\cdot,X_i)-\mu_k(P)}_{H(k)}^2 O_p\left(\frac{1}{n}\right) \nonumber\\
&=O_p\left(\frac{1}{n}\right) \label{I5_Theorem4}.
\end{align}
Also, we see that
\begin{align*}
I_6
&=\norm{\frac{1}{n}\sum_{i=1}^{n}B(X_i)\otimes \{(k(\cdot,X_i)-\mu_k(P))^{\otimes 2}-\Sigma_k(P)\}}_{H(k)^{\otimes 4}}\\
&\leq \frac{1}{n}\sum_{i=1}^{n} \norm{B(X_i)}_{H(k)^{\otimes 2}}\norm{(k(\cdot,X_i)-\mu_k(P))^{\otimes 2}-\Sigma_k(P)}_{H(k)^{\otimes 2}}\\
&\leq \left(
\frac{1}{n}\sum_{i=1}^{n} \norm{B(X_i)}_{H(k)^{\otimes 2}}^2
\right)^{1/2} 
\left(
\frac{1}{n}\sum_{i=1}^{n} \norm{(k(\cdot,X_i)-\mu_k(P))^{\otimes 2}-\Sigma_k(P)}_{H(k)^{\otimes 2}}^2
\right)^{1/2}
\end{align*}
by Cauchy--Schwartz inequality, and we get
\[
\frac{1}{n}\sum_{i=1}^{n} \norm{(k(\cdot,X_i)-\mu_k(P))^{\otimes 2}-\Sigma_k(P)}_{H(k)^{\otimes 2}}^2 \xrightarrow{P} \mathbb{E}[\norm{(k(\cdot,X_i)-\mu_k(P))^{\otimes 2}-\Sigma_k(P)}_{H(k)^{\otimes 2}}^2]
\]
as $n \to \infty$.
Hence, we get
\begin{equation}\label{I6_Theorem4}
I_6=\left(\frac{1}{n}\sum_{i=1}^{n} \norm{(k(\cdot,X_i)-\mu_k(P))^{\otimes 2}-\Sigma_k(P)}_{H(k)^{\otimes 2}}^2\right)^{1/2} O_p\left(\frac{1}{\sqrt{n}}\right)=O_p\left(\frac{1}{\sqrt{n}}\right)
\end{equation}
since $I_5 = O_p(1/n)$.
By the same argument as for $I_6=O_p(1/\sqrt{n})$ in (\ref{I6_Theorem4}), we have
\begin{align}
I_7
&=\norm{\frac{1}{n}\sum_{i=1}^{n} \{(k(\cdot,X_i)-\mu_k(P))^{\otimes 2}-\Sigma_k(P)\} \otimes B(X_i) }_{H(k)^{\otimes 4}} \nonumber\\
&\leq \frac{1}{n} \norm{B(X_i)}_{H(k)^{\otimes 2}}\norm{(k(\cdot,X_i)-\mu_k(P))^{\otimes 2}-\Sigma_k(P)}_{H(k)^{\otimes 2}}\nonumber\\
&=O_p\left(\frac{1}{\sqrt{n}}\right). \label{I7_Theorem4}
\end{align}
The equations (\ref{I1_Theorem4}), (\ref{I2_Theorem4}), (\ref{I3_Theorem4}), (\ref{I4_Theorem4}), (\ref{I5_Theorem4}), (\ref{I6_Theorem4}) and (\ref{I7_Theorem4}) are combined into (\ref{Upsilonhat-Upsilon}), which leads to
\begin{equation}\label{Upsilon_hat}
\norm{\widehat{\Upsilon}^{(n)}-\Upsilon}_{H(k)^{\otimes 4}}=O_p\left(\frac{1}{\sqrt{n}}\right).
\end{equation}
Therefore, we get
\begin{equation}\label{Eq_order_lambda}
\sum_{\ell=1}^{\infty}(\widehat{\lambda}_{\ell}^{(n)}-\lambda_{\ell})=O_p\left(\frac{1}{\sqrt{n}}\right)
\end{equation}
by (\ref{Ineq_lambda_Upsilon}) and (\ref{Upsilon_hat}), that is $\mathbb{E}[W'] \to 0$, as $n \to \infty$.

Next, we consider $V[W']$.
By $\widehat{\Upsilon}$ and $\Upsilon$ are compact Hermitian operators and (28) of \cite{Bhatia1994},
\[
\sum_{\ell=1}^{\infty} (\widehat{\lambda}_{\ell}^{(n)}-\lambda_{\ell})^2
\leq \norm{\widehat{\Upsilon}^{(n)}-\Upsilon}_{H(k)^{\otimes 4}}^2.
\]
Futhermore, we have got
$\norm{\widehat{\Upsilon}^{(n)}-\Upsilon}_{H(k)^{\otimes 4}}^2=O_p(1/n)$, hence 
\begin{equation}
\mathbb{E}\left[
\sum_{\ell=1}^{\infty}(\widehat{\lambda}_{\ell}^{(n)}-\lambda_{\ell})^2
\right]
\to 0,
\end{equation}
as $n \to \infty$.
Also, 
\begin{equation}\label{Eq_Cov}
Cov\left(\sum_{\ell=1}^{\infty}({\widehat{\lambda}}^{(n)}_{\ell}-\lambda_{\ell}),\sum_{\ell'=1}^{\infty}(\widehat{\lambda}_{\ell'}-\lambda_{\ell'})\right)
=\mathbb{E}\left[
\left(
\sum_{\ell=1}^{\infty}({\widehat{\lambda}}^{(n)}_{\ell}-\lambda_{\ell})
\right)^2
\right]
+\left\{
\mathbb{E}\left[\sum_{\ell=1}^{\infty}({\widehat{\lambda}}^{(n)}_{\ell}-\lambda_{\ell})\right]
\right\}^2
\end{equation}
and since $\left(
\sum_{\ell=1}^{\infty}({\widehat{\lambda}}^{(n)}_{\ell}-\lambda_{\ell})
\right)^2 =O_p(1/n)$ by (\ref{Eq_order_lambda}), we obtain 
\begin{equation}\label{Eq_lambda2}
\mathbb{E}\left[
\left(
\sum_{\ell=1}^{\infty}({\widehat{\lambda}}^{(n)}_{\ell}-\lambda_{\ell})
\right)^2
\right] \to 0,
\end{equation}
as $n \to \infty$.
In addition, $\mathbb{E}\left[\sum_{\ell=1}^{\infty}({\widehat{\lambda}}^{(n)}_{\ell}-\lambda_{\ell})\right] \to 0,$ as $n \to \infty$ by (\ref{Eq_order_lambda}), we get  
\[
Cov\left(\sum_{\ell=1}^{\infty}({\widehat{\lambda}}^{(n)}_{\ell}-\lambda_{\ell}),\sum_{\ell'=1}^{\infty}(\widehat{\lambda}_{\ell'}-\lambda_{\ell'})\right) \to 0,
\]
as $n \to \infty$ by (\ref{Eq_Cov}) and (\ref{Eq_lambda2}).
Therefore, $V[W'] \to 0,~~n \to \infty$.

Finally, we shall show $W' \xrightarrow{P}0$, as $n \to \infty$.
From Chebyshev's inequality, for any $\varepsilon>0$,
\begin{align*}
P(|W'| \geq \varepsilon)
&=P(|W'| -|\mathbb{E}[W']|\geq \varepsilon-|\mathbb{E}[W']|)\\
&\leq P(|W' -\mathbb{E}[W']|\geq \varepsilon-|\mathbb{E}[W']|)\\
&\leq \frac{V[W']}{(\varepsilon -|\mathbb{E}[W']|)^2}\\
& \to 0,
\end{align*}
as $n \to \infty$.
Therefore, $W' \xrightarrow{P} 0$ as $n \to \infty$.

\subsection{Proof of Proposition \ref{Lem_eigenvaue_same}}
Let $\widetilde{h}:\mathcal{H} \times \mathcal{H} \to \mathbb{R}$ be a kernel defined as,
\[
\widetilde{h}(x,y)=\left<(k(\cdot,x)-\mu_k(\widehat{P}))^{\otimes 2} -\Sigma_k(\widehat{P}),(k(\cdot,y)-\mu_k(\widehat{P}))^{\otimes 2} -\Sigma_k(\widehat{P}) \right>_{H(k)^{\otimes 2}},~~x,y \in \mathcal{H}
\]
and the associated function $\widetilde{h}(\cdot,x)$ represent $\widetilde{h}(\cdot,x)=(k(\cdot,x)-\mu_k(\widehat{P}))^{\otimes 2} -\Sigma_k(\widehat{P})$  for all $x \in \mathcal{H}$.
Then, $T: H(k)^{\otimes 2} \to \mathbb{R}^n$ is defined by
\[
T(A)=
\begin{bmatrix}
\left<A,\widetilde{h}(\cdot,X_1)\right>_{H(k)^{\otimes 2}}\\ 
\vdots\\ 
\left<A,\widetilde{h}(\cdot,X_n)\right>_{H(k)^{\otimes 2}}
\end{bmatrix} 
\]
for any $A \in H(k)^{\otimes 2}$.
The conjugate operator $T^*$ of $T$ (see Section VI.2 in \cite{reed1981functional} for details) is obtained as $T^* \underline{a}=\sum_{i=1}^{n}a_i \widetilde{h}(\cdot,X_i)$  for all $\underline{a}=\begin{bmatrix}
a_1& \cdots &a_n 
\end{bmatrix}^T  \in  \mathbb{R}^n$, because for all $A \in H(k)^{\otimes 2}$,
\begin{align*}
&\left<T^* \underline{a} , A\right>_{H(k)^{\otimes 2}}\\
&=\underline{a}^T(TA)\\
&=\underline{a}^T 
\begin{bmatrix}
\left<\widetilde{h}(\cdot,X_1),A\right>_{H(k)^{\otimes 2}}\\ 
\vdots\\ 
\left<\widetilde{h}(\cdot,X_1),A\right>_{H(k)^{\otimes 2}}
\end{bmatrix} \\
&=a_1 \left<\widetilde{h}(\cdot,X_1),A\right>_{H(k)^{\otimes 2}}+\cdots+a_n \left<\widetilde{h}(\cdot,X_1),A\right>_{H(k)^{\otimes 2}}\\
&=\left<\sum_{i=1}^{n}a_i \widetilde{h}(\cdot,X_i),A\right>_{H(k)^{\otimes 2}}.
\end{align*}
Let $\lambda$ and $A$ be the eigenvalue and eigenvector of $\widehat{\Upsilon}^{(n)}$, respectevely. 
Then, it is holds that
\[
\frac{1}{n}\sum_{j=1}^{n}\left<h(\cdot,X_j),A\right>_{H(k)^{\otimes 2}} h(\cdot,X_j)=\frac{1}{n}\sum_{j=1}^{n}\{\widetilde{h}(\cdot,X_j)\}^{\otimes 2} A= \lambda A
\]
from the definition of eigenvalue and eigenvector.
By mapping both sides with $T$,
\begin{align*}
&\frac{1}{n} 
\begin{bmatrix}
\widetilde{h}(X_1,X_1)&\cdots  &\widetilde{h}(X_1,X_n)  \\ 
\vdots& \ddots & \vdots  \\ 
\widetilde{h}(X_n,X_1)&\cdots  &\widetilde{h}(X_n,X_n) 
\end{bmatrix} 
\begin{bmatrix}
\left<A,\widetilde{h}(\cdot,X_1)\right>_{H(k)^{\otimes 2}}\\ 
\vdots\\ 
\left<A,\widetilde{h}(\cdot,X_1)\right>_{H(k)^{\otimes 2}}
\end{bmatrix}\\ 
&=
\frac{1}{n}\sum_{j=1}^{n}\left<\widetilde{h}(\cdot,X_j),A\right>_{H(k)^{\otimes 2}}
\begin{bmatrix}
\widetilde{h}(X_1,X_j)\\ 
\vdots\\ 
\widetilde{h}(X_n,X_j)
\end{bmatrix} 	\\
&=\frac{1}{n}\sum_{j=1}^{n}\left<\widetilde{h}(\cdot,X_j),A\right>_{H(k)^{\otimes 2}} T(\widetilde{h}(\cdot,X_j)) \\
&=\lambda T(A)\\
&=\lambda\begin{bmatrix}
\left<A,\widetilde{h}(\cdot,X_1)\right>_{H(k)^{\otimes 2}}\\ 
\vdots\\ 
\left<A,\widetilde{h}(\cdot,X_n)\right>_{H(k)^{\otimes 2}}
\end{bmatrix}. 
\end{align*} 
Hence, the eigenvalues of $\widehat{\Upsilon}^{(n)}$ are that of $H/n$.

Conversely, let $\tau$ and 
$\underline{u}=
\begin{bmatrix}
u_1&\cdots&u_n 
\end{bmatrix}^T$ be the eigenvalue and correspondent eigenvector of $H/n$, then
\[
\frac{1}{n}\sum_{j=1}^{n}u_j 
\begin{bmatrix}
\widetilde{h}(X_1,X_j)\\ 
\vdots\\ 
\widetilde{h}(X_n,X_j)
\end{bmatrix} =\frac{1}{n}H\underline{u}=\lambda \underline{u},
\]
and 
\begin{align*}
&\widehat{\Upsilon}^{(n)}\left\{ \sum_{j=1}^{n} u_j \widetilde{h}(\cdot,X_j)\right\}\\
&=\frac{1}{n} \sum_{i=1}^{n} \left\{\widetilde{h}(\cdot,X_i)\right\}^{\otimes 2} \left\{ \sum_{j=1}^{n} u_j \widetilde{h}(\cdot,X_j)\right\}\\
&=\frac{1}{n}\sum_{i=1}^{n} \left<\widetilde{h}(\cdot,X_i), \sum_{j=1}^{n} u_j\widetilde{h}(\cdot,X_j)\right>_{H(k)^{\otimes 2}} \widetilde{h}(\cdot,X_i)\\
&=\frac{1}{n}\sum_{i,j=1}^{n} u_j \left<\widetilde{h}(\cdot,X_i),\widetilde{h}(\cdot,X_j)\right>_{H(k)^{\otimes 2}} \widetilde{h}(\cdot,X_i)\\
&=\frac{1}{n}\sum_{j=1}^{n} u_j \sum_{i=1}^{n} \widetilde{h}(X_i,X_j) \widetilde{h}(\cdot,X_i)\\
&=\frac{1}{n}\sum_{j=1}^{n}u_j 
T^* \left(\begin{bmatrix}
\widetilde{h}(X_1,X_j)\\ 
\vdots\\ 
\widetilde{h}(X_n,X_j)
\end{bmatrix}\right)\\
&=\lambda T^*(\underline{u})\\
&= \lambda \sum_{i=1}^{n}u_i \widehat{h}(\cdot,X_i)	
\end{align*}
form mapping both sides with $T^*$, hence the eigenvalue of $H/n$ are that of $\widehat{\Upsilon}^{(n)}$.

\subsection{Proof of Proposition \ref{Variance_MVD_P=Q}}
We see that
\begin{align*}
	&\mathbb{E}[(n+m)\widehat{T}^2_{n,m}]\\
	&=\frac{n+m}{n^2}\sum_{i,s=1}^{n}\mathbb{E}[ h(X_i,X_s)]
	+\frac{n+m}{m^2} \sum_{j,t=1}^{m}\mathbb{E}[h(Y_j,Y_t)]
	-\frac{2(n+m)}{nm} \sum_{i=1}^{n} \sum_{j=1}^{m} \mathbb{E}[h(X_i,Y_j)],
\end{align*}
where $h(x,y)$ is in (\ref{Eq_h(x,y)}).
Since we have
\begin{align*}
	\mathbb{E}[h(X_1,X_2)]
	&=\mathbb{E}\left[\left<(k(\cdot,X_1)-\mu_k(P))^{\otimes 2}-\Sigma_k(P),(k(\cdot,X_2)-\mu_k(P))^{\otimes 2}-\Sigma_k(P)\right>_{H(k)^{\otimes 2}} \right]\\
	&=\left<\mathbb{E}[(k(\cdot,X_1)-\mu_k(P))^{\otimes 2}-\Sigma_k(P)],\mathbb{E}[(k(\cdot,X_2)-\mu_k(P))^{\otimes 2}-\Sigma_k(P)]\right>_{H(k)^{\otimes 2}}\\
	&=0
\end{align*}
and
\begin{align*}
	\mathbb{E}[h(X_1,X_1)]
	&=\mathbb{E}\left[\left<(k(\cdot,X_1)-\mu_k(P))^{\otimes 2}-\Sigma_k(P),(k(\cdot,X_1)-\mu_k(P))^{\otimes 2}-\Sigma_k(P)\right>_{H(k)^{\otimes 2}}\right]\\
	&=\mathbb{E}\left[\left<\left\{(k(\cdot,X_1)-\mu_k(P))^{\otimes 2}-\Sigma_k(P)\right\}^{\otimes 2}, I\right>_{H(k)^{\otimes 4}}\right]\\
	&=\left<\mathbb{E}\left[\left\{(k(\cdot,X_1)-\mu_k(P))^{\otimes 2}-\Sigma_k(P)\right\}^{\otimes 2}\right], I\right>_{H(k)^{\otimes 4}}\\
	&=\left<\Upsilon,I\right>_{H(k)^{\otimes 4}},
\end{align*}
under $P=Q$, it follows that
\begin{equation}\label{Eq_That_Expactation}
\mathbb{E}[(n+m)\widehat{T}^2_{n,m}]
=\frac{n+m}{n^2} n \left<\Upsilon,I\right>_{H(k)^{\otimes 4}}
+\frac{n+m}{m^2} m \left<\Upsilon,I\right>_{H(k)^{\otimes 4}}\\
=\frac{(n+m)^2}{nm}\left<\Upsilon,I\right>_{H(k)^{\otimes 4}}.
\end{equation}
Next, we consider $\mathbb{E}\left[\left\{(n+m)\widehat{T}_{n,m}^2\right\}^2 \right]$.
It follows from direct calculations that
\begin{align*}
	&\mathbb{E}\left[\left\{(n+m)\widehat{T}_{n,m}^2\right\}^2\right]\\
	&=\mathbb{E}\left[
	\left\{
	\frac{n+m}{n^2}\sum_{i,s=1}^{n}h(X_i,X_s)+\frac{n+m}{m^2} \sum_{j,t=1}^{m} h(Y_j,Y_t) -\frac{2(n+m)}{nm} \sum_{i=1}^{n} \sum_{j=1}^{m} h(X_i,Y_j)
	\right\}^2
	\right]\\
	&=\frac{(n+m)^2}{n^4}\mathbb{E}\left[
	\left\{\sum_{i,s=1}^{n} h(X_i,X_s)\right\}^2
	\right]
	+\frac{(n+m)^2}{m^4} \mathbb{E}\left[
	\left\{
	\sum_{j,t=1}^{m}h(Y_j,Y_t)
	\right\}^2
	\right]\\
	&~~~~~
	+\frac{4(n+m)^2}{n^2m^2} \mathbb{E}\left[
	\left\{
	\sum_{i=1}^{n} \sum_{j=1}^{m} h(X_i,Y_j)
	\right\}^2
	\right] 
	+\frac{2(n+m)^2}{n^2m^2}\mathbb{E}\left[
	\left\{
	\sum_{i,s=1}^{n} h(X_i,X_s)
	\right\}
	\left\{
	\sum_{j,t=1}^{m} h(Y_j,Y_t)
	\right\}
	\right]\\
	&~~~~~
	-\frac{4(n+m)^2}{n^3m}\mathbb{E}\left[
	\sum_{i,s,\ell=1}^{n} \sum_{j=1}^{m}h(X_i,X_s) h(X_\ell,Y_j)
	\right]
	-\frac{4(n+m)^2}{nm^3} \mathbb{E}\left[
	\sum_{i=1}^{n} \sum_{j,t,k=1}^{m} h(Y_j,Y_t)h(X_i,Y_k)
	\right].
\end{align*}
A straightforward but lengthy computation yields that 
\begin{equation}\label{Eq_That_Expectation_2}
\mathbb{E}\left[
\left\{\sum_{i,s=1}^{n} h(X_i,X_s)\right\}^2
\right]
=n\left<A,I\right>_{H(k)^{\otimes 8}}
+2n(n-1)\norm{\Upsilon}^2_{H(k)^{\otimes 4}}
+n(n-1)\left<\Upsilon,I\right>_{H(k)^{\otimes 4}}^2,
\end{equation}
where $A=\mathbb{E}[\{(k(\cdot,X_1)-\mu_k(P))^{\otimes 2}-\Sigma_k(P)\}^{\otimes 4}]$.
In addition, we obtain from direct calcuration that
\begin{align*}
	&\mathbb{E}\left[
	\left\{
	\sum_{i=1}^{n} \sum_{j=1}^{m} h(X_i,Y_j)
	\right\}^2
	\right]
	=nm\norm{\Upsilon}^2_{H(k)^{\otimes 4}},\\
	&\mathbb{E}\left[
	\left\{
	\sum_{i,s=1}^{n} h(X_i,X_s)
	\right\}
	\left\{
	\sum_{j,t=1}^{m} h(Y_j,Y_t)
	\right\}
	\right]
	=nm \left<\Upsilon,I\right>_{H(k)^{\otimes 2}}^2.
\end{align*}
Therefore, using (\ref{Eq_That_Expactation}) and (\ref{Eq_That_Expectation_2})
\begin{align*}
	&V[(n+m)\widehat{T}_{n,m}]\\
	&=\mathbb{E}\left[\left\{(n+m)\widehat{T}^2_{n,m}\right\}^2\right] -\{\mathbb{E}[(n+m)\widehat{T}^2_{n,m}]\}^2\\
	&=\frac{(n+m)^2}{n^4}\left(n\left<A,I\right>_{H(k)^{\otimes 8}}+2n(n-1)\norm{\Upsilon}^2_{H(k)^{\otimes 4}}+n(n-1)\left<\Upsilon,I\right>_{H(k)^{\otimes 4}}^2 \right)\\
	&~~~~~+\frac{(n+m)^2}{m^4}\left(m\left<A,I\right>_{H(k)^{\otimes 8}}+2m(m-1)\norm{\Upsilon}^2_{H(k)^{\otimes 4}}+m(m-1)\left<\Upsilon,I\right>_{H(k)^{\otimes 4}}^2 \right)\\
	&~~~~~+\frac{4(n+m)^2}{n^2m^2} nm \norm{\Upsilon}^2_{H(k)^{\otimes 4}}
	+\frac{2(n+m)^2}{n^2m^2} nm \left<\Upsilon,I\right>^2_{H(k)^{\otimes 2}}
	-\frac{(n+m)^4}{n^2m^2} \left<\Upsilon,I\right>^2_{H(k)^{\otimes 4}}\\
	&=\frac{2(n+m)^4}{n^2m^2}\norm{\Upsilon}^2_{H(k)^{\otimes 4}}+O\left(\frac{1}{n}\right)+O\left(\frac{1}{m}\right).
\end{align*}
\subsection{Proof of Proposition \ref{Variance_MMD_P=Q}}
Since
\[
(n+m)\widehat{\Delta}^2_{n,m}
=\frac{n+m}{n^2} \sum_{i,s=1}^{n}k(X_i,X_s)
+\frac{n+m}{m^2}\sum_{j,t=1}^{m} k(Y_j,Y_t)
-\frac{2(n+m)}{nm} \sum_{i=1}^{n} \sum_{j=1}^{m} k(X_i,Y_j),
\]
first, we need to calculate
\[
\mathbb{E}[(n+m)\widehat{\Delta}^2_{n,m}]
=\frac{n+m}{n^2} \sum_{i,s=1}^{n} \mathbb{E}[k(X_i,X_s)]
+\frac{n+m}{m^2} \sum_{j,t=1}^{m} \mathbb{E}[k(Y_j,Y_t)]
-\frac{2(n+m)}{nm} \sum_{i=1}^{n} \sum_{j=1}^{m} \mathbb{E}[k(X_i,Y_j)].
\]

From the expected values of each term are obtained as
\begin{align*}
	&\mathbb{E}[k(X_1,X_2)]=\norm{\mu_k(P)}^2_{H(k)},\\
	&\mathbb{E}[k(X_1,X_1)]
	=\left<\Sigma_k(P),I\right>_{H(k)^{\otimes 2}} +\norm{\mu_k(P)}_{H(k)}^2
\end{align*}
we get
\begin{align}
	\mathbb{E}[(n+m)\widehat{\Delta}^2_{n,m}]
	=\frac{(n+m)^2}{nm} \left<\Sigma_k(P),I\right>_{H(k)^{\otimes 2}} \label{MMD_firstE}
\end{align}
under $P = Q$.

Next, the second moment of $(n+m)\widehat{\Delta}^2_{n,m}$ is 
\begin{align}
	&\mathbb{E}[\{(n+m)\widehat{\Delta}^2_{n,m}\}^2] \nonumber\\
	&=\mathbb{E}\left[
	\left\{
	\frac{n+m}{n^2} \sum_{i,s=1}^{n} k(X_i,X_s)
	+\frac{n+m}{m^2} \sum_{j,t=1}^{m} k(Y_j,Y_t)
	-\frac{2(n+m)}{nm} \sum_{i=1}^{n} \sum_{j=1}^{m} k(X_i,Y_j)
	\right\}^2
	\right]\nonumber\\
	&=\frac{(n+m)^2}{n^4}\mathbb{E} \left[
	\left\{
	\sum_{i,s=1}^{n}k(X_i,X_s)
	\right\}^2
	\right]
	+\frac{(n+m)^2}{m^4}\mathbb{E}\left[
	\left\{
	\sum_{j,t=1}^{m} k(Y_j,Y_t)
	\right\}^2
	\right]\nonumber\\
	&~~~~~
	+\frac{4(n+m)^2}{n^2m^2} \mathbb{E} \left[
	\left\{
	\sum_{i=1}^{n} \sum_{j=1}^{m} k(X_i,Y_j)
	\right\}^2
	\right]
	+\frac{2(n+m)^2}{n^2m^2} \mathbb{E}\left[
	\left\{
	\sum_{i,s=1}^{n} k(X_i,X_s)
	\right\}
	\left\{
	\sum_{j,t=1}^{m} k(Y_j,Y_t)
	\right\}
	\right]\nonumber\\
	&~~~~~
	-\frac{4(n+m)^2}{n^3m} \mathbb{E} \left[
	\sum_{i,s,\ell=1}^{n} \sum_{j=1}^{m} k(X_i,X_s) k(X_{\ell},Y_j)
	\right]
	-\frac{4(n+m)^2}{nm^3} \mathbb{E} \left[
	\sum_{i=1}^{n} \sum_{j,t,k=1}^{m} k(Y_j,Y_t) k(X_i,Y_k)
	\right]. \label{MMD_secondE}
\end{align}
These expectations are obtained as 
\begin{align}
	&\mathbb{E}\left[\left\{
	\sum_{i,s=1}^{n}k(X_i,X_s)
	\right\}^2
	\right] \nonumber\\
	&=n \mathbb{E}[k(X_1,X_1)^2]
	+4n(n-1) \left<\mathbb{E}[k(\cdot,X_1)^{\otimes 2} k(\cdot,X_1)],\mu_k(P)\right>_{H(k)}
	+2n(n-1) \norm{\Sigma_k(P)}^2_{H(k)^{\otimes 2}} \nonumber\\
	&~~~~~
	+n(n-1) \left<\Sigma_k(P),I\right>_{H(k)^{\otimes 2}}^2
	+2n(n-1)^2\left<\Sigma_k(P),I\right>_{H(k)^{\otimes 2}} \norm{\mu_k(P)}^2_{H(k)} \nonumber\\
	&~~~~~
	+4n(n-1)^2 \left<\Sigma_k(P),\mu_k(P)^{\otimes 2}\right>_{H(k)^{\otimes 2}}
	+n(n-1)(n^2+n-3) \norm{\mu_k(P)}^4_{H(k)}, \label{MMD_E1}\\
	&\mathbb{E} \left[
	\left\{
	\sum_{i=1}^{n} \sum_{j=1}^{m}k(X_i,Y_j)
	\right\}^2
	\right] \nonumber\\
	&=nm\norm{\Sigma_k(P)}^2_{H(k)^{\otimes 2}} 
	+nm(n+m) \left<\Sigma_k(P),\mu_k(P)^{\otimes 2}\right>_{H(k)^{\otimes 2}}
	+n^2m^2 \norm{\mu_k(P)}^4_{H(k)}, \label{MMD_E2}\\
	&\mathbb{E}\left[
	\left\{
	\sum_{i,s=1}^{n} k(X_i,X_s)
	\right\}
	\left\{
	\sum_{j,t=1}^{m} k(Y_j,Y_t)
	\right\}
	\right]\nonumber\\
	&=nm \left<\Sigma_k(P),I\right>_{H(k)^{\otimes 2}}^2 
	+nm(n+m) \left<\Sigma_k(P),I\right>_{H(k)} \norm{\mu_k(P)}^2_{H(k)}
	+n^2m^2 \norm{\mu_k(P)}^4, \label{MMD_E3}\\
	&\mathbb{E} \left[
	\sum_{i,s,\ell=1}^{n} \sum_{j=1}^{m} k(X_i,X_s) k(X_{\ell},Y_j)
	\right]\nonumber\\
	&=nm \left<\mathbb{E}[k(\cdot,X_1)^{\otimes 2} k(\cdot,X_1)],\mu_k(P)\right>_{H(k)}
	+n(n-1)m \left<\Sigma_k(P),I\right>_{H(k)^{\otimes 2}} \norm{\mu_k(P)}^2_{H(k)}\nonumber\\
	&~~~~~
	+2n(n-1)m \left<\Sigma_k(P),\mu_k(P)^{\otimes 2}\right>_{H(k)^{\otimes 2}}
	+n(n-1)(n+1)m \norm{\mu_k(P)}^4_{H(k)},\label{MMD_E4}\\
	&\mathbb{E} \left[
	\sum_{i=1}^{n} \sum_{j,t,k=1}^{m} k(Y_j,Y_t) k(X_i,Y_k)
	\right] \nonumber\\
	&=nm \left<\mathbb{E}[k(\cdot,X_1)^{\otimes 2} k(\cdot,X_1)],\mu_k(P)\right>_{H(k)}
	+m(m-1)n \left<\Sigma_k(P),I\right>_{H(k)^{\otimes 2}} \norm{\mu_k(P)}^2_{H(k)}\nonumber\\
	&~~~~~
	+2m(m-1)n \left<\Sigma_k(P),\mu_k(P)^{\otimes 2}\right>_{H(k)^{\otimes 2}}
	+m(m-1)(m+1)n \norm{\mu_k(P)}^4_{H(k)} \label{MMD_E5}.
\end{align}
The combining (\ref{MMD_secondE}) and (\ref{MMD_E1})-(\ref{MMD_E5}) provides that
\begin{align}
	&\mathbb{E}[\{(n+m)\widehat{\Delta}^2_{n,m}\}^2] \nonumber\\
	&=\frac{2(n+m)^4}{n^2m^2} \norm{\Sigma_k(P)}^2_{H(k)^{\otimes 2}}
	+\frac{(n+m)^4}{n^2m^2} \left<\Sigma_k(P),I\right>_{H(k)^{\otimes 2}}^2
	+O\left(\frac{1}{n}\right) +O\left(\frac{1}{m}\right). \label{MMD_secondEE}
\end{align}
Therefore, from (\ref{MMD_firstE}) and (\ref{MMD_secondEE}), the variance of $(n+m)\widehat{\Delta}^2_{n,m}$ is 
\begin{align*}
	V[(n+m)\widehat{\Delta}^2_{n,m}]
	&=\mathbb{E}\left[\left\{(n+m)\widehat{\Delta}^2_{n,m}\right\}^2\right]-\left\{\mathbb{E}[(n+m)\widehat{\Delta}^2_{n,m}]\right\}^2\\
	&=\frac{2(n+m)^4}{n^2m^2} \norm{\Sigma_k(P)}^2_{H(k)^{\otimes 2}}
	+O\left(\frac{1}{n}\right) +O\left(\frac{1}{m}\right).
\end{align*}
\subsection{Proof of (\ref{H_matrix_Gram})}
The $(i,j)$-th element of the matrix $H$ is
\begin{align*}
H_{ij}
&=\left<(k(\cdot,X_i)-{\mu}_k(\widehat{P}))^{\otimes 2} -{\Sigma}_k(\widehat{P}),(k(\cdot,X_j)-{\mu}_k(\widehat{P}))^{\otimes 2} -{\Sigma}_k(\widehat{P})\right>_{H(k)^{\otimes 2}}\\
&=\left<(k(\cdot,X_i)-{\mu}_k(\widehat{P}))^{\otimes 2} ,(k(\cdot,X_j)-{\mu}_k(\widehat{P}))^{\otimes 2} \right>_{H(k)^{\otimes 2}}
-\left<{\Sigma}_k(\widehat{P}),(k(\cdot,X_j)-{\mu}_k(\widehat{P}))^{\otimes 2} \right>_{H(k)^{\otimes 2}}\\
&~~~~~
-\left<(k(\cdot,X_i)-{\mu}_k(\widehat{P}))^{\otimes 2} ,{\Sigma}_k(\widehat{P})\right>_{H(k)^{\otimes 2}}
+\left<{\Sigma}_k(\widehat{P}),{\Sigma}_k(\widehat{P})\right>_{H(k)^{\otimes 2}}.
\end{align*}
Each term of this $H_{ij}$ can be expressed as
\begin{align*}
&\left<(k(\cdot,X_i)-{\mu}_k(\widehat{P}))^{\otimes 2} ,(k(\cdot,X_j)-{\mu}_k(\widehat{P}))^{\otimes 2} \right>_{H(k)^{\otimes 2}}\\
&=\left<k(\cdot,X_i)-{\mu}_k(\widehat{P}) ,k(\cdot,X_j)-{\mu}_k(\widehat{P}) \right>_{H(k)}^2\\
&=\left\{
k(X_i,X_j)-\mu_k(\widehat{P})(X_i)-\mu_k(\widehat{P})(X_j)+\left<\mu_k(\widehat{P}),\mu_k(\widehat{P})\right>_{H(k)}
\right\}^2\\
&=\left\{
k(X_i,X_j) -\frac{1}{n}\sum_{s=1}^{n}k(X_j,X_s)-\frac{1}{n}\sum_{\ell=1}^{n}k(X_i,X_{\ell})+\frac{1}{n^2}\sum_{s,\ell=1}^{n}k(X_s,X_{\ell})
\right\}^2\\
&=\left(\widetilde{K}_{ij}\right)^2\\
&=\left(\widetilde{K} \odot \widetilde{K} \right)_{ij},	\\
&\left<\Sigma_k(\widehat{P}),(k(\cdot,X_i)-{\mu}_k(\widehat{P}))^{\otimes 2} \right>_{H(k)^{\otimes 2}}\\
&=\left<\frac{1}{n}\sum_{s=1}^{n}(k(\cdot,X_s)-{\mu}_k(\widehat{P}))^{\otimes 2},(k(\cdot,X_i)-{\mu}_k(\widehat{P}))^{\otimes 2} \right>_{H(k)^{\otimes 2}}\\
&=\frac{1}{n}\sum_{s=1}^{n}\left<(k(\cdot,X_s)-{\mu}_k(\widehat{P}))^{\otimes 2},(k(\cdot,X_i)-{\mu}_k(\widehat{P}))^{\otimes 2} \right>_{H(k)^{\otimes 2}}\\
&=\frac{1}{n}\sum_{s=1}^{n}\left<k(\cdot,X_s)-{\mu}_k(\widehat{P}),k(\cdot,X_i)-{\mu}_k(\widehat{P})\right>_{H(k)}^2\\
&=\frac{1}{n}\sum_{s=1}^{n}\left(\widetilde{K}_{sj}\right)^2\\
&=\frac{1}{n}\sum_{s=1}^{n} \left(\widetilde{K} \odot \widetilde{K}\right)_{sj}
\end{align*}
and
\begin{align*}
&\left<\Sigma_k(\widehat{P}),\Sigma_k(\widehat{P})\right>_{H(k)^{\otimes 2}}\\
&=\left<\frac{1}{n}\sum_{s=1}^{n}(k(\cdot,X_s)-{\mu}_k(\widehat{P}))^{\otimes 2},\frac{1}{n}\sum_{\ell=1}^{n}(k(\cdot,X_{\ell})-{\mu}_k(\widehat{P}))^{\otimes 2} \right>_{H(k)^{\otimes 2}}\\
&=\frac{1}{n^2} \sum_{s,\ell=1}^{n} \left<(k(\cdot,X_s)-{\mu}_k(\widehat{P}))^{\otimes 2},(k(\cdot,X_{\ell})-{\mu}_k(\widehat{P}))^{\otimes 2} \right>_{H(k)^{\otimes 2}}\\
&=\frac{1}{n^2} \sum_{s,\ell=1}^{n} \left<k(\cdot,X_s)-{\mu}_k(\widehat{P}),k(\cdot,X_{\ell})-{\mu}_k(\widehat{P})\right>_{H(k)}^2\\
&=\frac{1}{n^2} \sum_{s,\ell=1}^{n} \left(\widetilde{K}_{s\ell}\right)^2\\
&=\frac{1}{n^2} \sum_{s,\ell=1}^{n} \left(\widetilde{K} \odot \widetilde{K}\right)_{s\ell}.
\end{align*}
Therefore, 
\[
H_{ij}=\left(\widetilde{K} \odot \widetilde{K}\right)_{ij} 
-\frac{1}{n} \sum_{s=1}^{n} \left(\widetilde{K} \odot \widetilde{K}\right)_{sj} 
-\frac{1}{n} \sum_{s=1}^{n} \left(\widetilde{K} \odot \widetilde{K}\right)_{si}
+\frac{1}{n^2} \sum_{s,\ell=1}^{n} \left(\widetilde{K} \odot \widetilde{K}\right)_{s \ell},
\]
which gives the expression (\ref{H_matrix_Gram}).

\bibliography{Mylibrary.bib}

\begin{thebibliography}{10}

\bibitem{Alon}
U.~Alon, N.~Barkai, D.~A. Notterman, K.~Gish, S.~Ybarra, D.~Mack, and A.~J.
  Levine.
\newblock {Broad patterns of gene expression revealed by clustering analysis of
  tumor and normal colon tissues probed by oligonucleotide arrays}.
\newblock {\em Proceedings of the National Academy of Sciences},
  96(12):6745--6750, 1999.

\bibitem{Aronszajn}
N.~Aronszajn.
\newblock {Theory of reproducing kernels}.
\newblock {\em Transactions of the American Mathematical Society},
  68(3):337--404, 1950.

\bibitem{Bhatia1994}
R.~Bhatia and L.~Elsner.
\newblock {The Hoffman-Wielandt inequality in infinite dimensions}.
\newblock {\em Proceedings of the Indian Academy of Sciences - Mathematical
  Sciences}, 104(3):483--494, 1994.

\bibitem{Boente2018}
G.~Boente, D.~Rodriguez, and M.~Sued.
\newblock {Testing equality between several populations covariance operators}.
\newblock {\em Annals of the Institute of Statistical Mathematics},
  70(4):919--950, 2018.

\bibitem{Fukumizu}
K.~Fukumizu, A.~Gretton, X.~Sun, and B.~Sch{\"{o}}lkopf.
\newblock {Kernel measures of conditional dependence}.
\newblock {\em Advances in Neural Information Processing Systems 20 -
  Proceedings of the 2007 Conference}, pages 1--13, 2009.

\bibitem{Gretton2007}
A.~Gretton, K.~M. Borgwardt, M.~Rasch, B.~Sch{\"{o}}lkopf, and A.~J. Smola.
\newblock {A kernel method for the two-sample-problem}.
\newblock In B.~Sch{\"{o}}lkopf, J.~C. Platt, and T.~Hoffman, editors, {\em
  Advances in Neural Information Processing Systems}, volume~19, pages
  513--520. MIT Press, 2007.

\bibitem{Gretton2009a}
A.~Gretton, K.~Fukumizu, Z.~Harchaoui, and B.~Sriperumbudur.
\newblock {A fast, consistent kernel two-sample test}.
\newblock {\em Advances in Neural Information Processing Systems}, pages
  673--681, 2009.

\bibitem{hoffmann-jorgensen1976}
J.~Hoffmann-Jorgensen and G.~Pisier.
\newblock {The law of large numbers and the central limit theorem in Banach
  spaces}.
\newblock {\em The Annals of Probability}, 4(4):587--599, 1976.

\bibitem{Kellner_and_Celisse}
J.~Kellner and A.~Celisse.
\newblock {A one-sample test for normality with kernel methods}.
\newblock {\em Bernoulli}, 25(3):1816--1837, 2019.

\bibitem{Minh2006}
H.~Q. Minh, P.~Niyogi, and Y.~Yao.
\newblock {Mercer's theorem, feature maps, and smoothing}.
\newblock In {\em Proceedings of the 19th Annual Conference on Learning
  Theory}, COLT'06, pages 154--168, Berlin, Heidelberg, 2006. Springer-Verlag.

\bibitem{politis1999}
D.~Politis, D.~Wolf, J.~Romano, M.~Wolf, P.~Bickel, P.~Diggle, and S.~Fienberg.
\newblock {\em Subsampling}.
\newblock Springer Series in Statistics. Springer New York, 1999.

\bibitem{reed1981functional}
M.~Reed and B.~Simon.
\newblock {\em {Functional Analysis}}.
\newblock Methods of Modern Mathematical Physics. Elsevier Science, 1981.

\bibitem{Wand_and_Jones}
M.~P. Wand and M.~C. Jones.
\newblock {\em Kernel Smoothing}.
\newblock Chapman \& Hall, New York, 1994.

\end{thebibliography}
%
%
%
%
%
\end{document}